\documentclass[reqno]{amsart}
 \usepackage{amsbsy,amssymb,amscd,amsfonts,latexsym,amstext,delarray,
 amsmath,graphicx, color}
 \usepackage[dvips]{epsfig}
\usepackage{hyperref}
\input xypic

\definecolor{Green}{rgb}{0,1,0}
\definecolor{Blue}{RGB}{0,0,191}
\definecolor{mathmodecolor}{RGB}{0,102,0}
\definecolor{keywordcolor}{RGB}{0,51,151}
\definecolor{sourcebackgroundcolor}{RGB}{255,247,223}
\definecolor{unixagred}{RGB}{255,0,0}
\definecolor{lightgray}{RGB}{191,191,191}
\definecolor{green}{RGB}{1,191,191}

\newtheorem{thm}{Theorem}[section]
\newtheorem{prop}[thm]{Proposition}
\newtheorem{cor}[thm]{Corollary}
\newtheorem{lem}[thm]{Lemma}

\newtheorem{defn}[thm]{Definition}
\newtheorem{rem}[thm]{Remark}
\newtheorem{example}[thm]{Example}

\def\qqq{\,,\quad~\forall}

\def\Aut{{\rm Aut}}

\def\End{{\rm End}}

\def\Hom{{\rm Hom}}

\def\Spec{{\rm Spec\,}}
\def\Sp{{\rm Spec}\,}

\def\B{{\mathbb B}}
\def\C{{\mathbb C}}
\def\F{{\mathbb F}}
\def\K{{\mathbb K}}

\def\N{{\mathbb N}}
\def\P{{\mathbb P}}
\def\Q{{\mathbb Q}}
\def\R{{\mathbb R}}
\def\Z{{\mathbb Z}}

\def\cD{{\mathcal D}}

\def\cF{{\mathcal F}}

\def\cO{{\mathcal O}}

\def\cU{{\mathcal U}}

\newcommand{\ie}{{\it i.e.\/}\ }
\newcommand{\eg}{{\it e.g.\/}\ }
\newcommand{\cf}{{\it cf.\/}\ }
\newcommand{\opcit}{{\it op.cit.\/}\ }

\def\sin{{{\rm sin}}}
\def\cos{{{\rm cos}}}

\def\Hom {{\mbox{Hom}}}

\def\End{{\mbox{End}}}

\def\ffp{\mathfrak{p}}

\def\Mo{\mathfrak{Mo}}
\def\An{\mathfrak{ Ring}}
\def\Ab{\mathfrak{ Ab}}
\def\Mr{\mathfrak{ MR}}

\def\Sch{\mathfrak{  Sch}}

\def\Se{\mathfrak{ Sets}}

\def\GS{\mathfrak G\mathcal S}

\def\ssp{{\mathfrak{ spec}\,}}

\def\spad{{\P^1_{\F_1}}}

\def\Omm{\Omega}

\def\dio{{idempotent semi-ring  }}
\def\plus  {\uplus}
\def\dioid{{semi-ring of characteristic $1$ }}
\def\dioids{{semi-rings of characteristic $1$}}
\def\mul{{\mathbb G}}
\def\rmax{\R_+^{\rm max}}

\newcommand{\nil}[1]{}

\title
{Characteristic $1$, entropy and the absolute point}

\author[Connes]{Alain Connes}
\author[Consani]{Caterina Consani}
\address{A.~Connes: Coll\`ege de France \\
3, rue d'Ulm \\ Paris, F-75005 France
\\ I.H.E.S. and Vanderbilt
University} \email{alain\@@connes.org}
\address{C.~Consani: Mathematics Department \\ Johns Hopkins
University \\ Baltimore, MD 21218 USA} \email{kc\@@math.jhu.edu}
\thanks{The authors are partially supported by the NSF grant
DMS-FRG-0652164. The first author thanks C.~Breuil for pointing out the reference to M. Krasner \cite{Krasner}.
 Both authors acknowledge the paper \cite{Lescot} of Paul Lescot, that brings up the role of the idempotent semi-field $\B$}


\begin{document}
\begin{abstract} We show that the mathematical meaning of working in characteristic one  is directly connected to the fields of idempotent analysis
and tropical algebraic geometry  and  we relate this idea to the notion of the absolute point $\Spec\F_1$. After introducing the notion of ``perfect" semi-ring of characteristic one,  we explain how to adapt the construction of the Witt ring in characteristic $p>1$ to the limit case of characteristic one. This construction also unveils an interesting  connection with entropy and thermodynamics, while shedding a
new light on the classical Witt construction itself. We simplify our earlier construction of the geometric
realization of an $\F_1$-scheme and extend our earlier computations of the zeta function to cover the case of $\F_1$-schemes with torsion.  Then, we  show that the study of the additive structures on monoids provides a natural map $M\mapsto A(M)$ from monoids to sets which comes close to
fulfill the requirements for the hypothetical curve $\overline {\Sp\,\Z}$ over the
 absolute point $\Spec\F_1$. Finally, we test the computation of the zeta function on elliptic curves over the rational numbers.
\end{abstract}

\maketitle

 \tableofcontents

\section{Introduction}

The main goal of this article is to explore the mathematical meaning of working in characteristic one and to relate this idea with the notion of the absolute point $\Spec\F_1$.
In the first part of the paper, we explain that an already well established branch of mathematics supplies a satisfactory answer to the question of the meaning of  mathematics in characteristic one. Our starting point is the investigation of the structure of fields of positive characteristic  $p>1$, in the degenerate case $p=1$. The  outcome is that this limit case is directly connected to the fields of idempotent analysis (\cf\cite{Maslov})
and tropical algebraic geometry (\cf\cite{Mikhalkin}, \cite{Gathmann}). In parallel with the development of classical mathematics over fields, its ``shadow'' idempotent version appears in characteristic one

\begin{gather}\label{Shidiagr}
\raisetag{-24pt} \,\hspace{30pt}\begin{xy}
(0,20) *+++{Real \ Analysis} *\frm{-,} ;
(0,0) *+++{Idempotent\  Analysis} *\frm{-,}**\dir{-};
(50,20) *+++{Geometry \ over \  \C} *\frm{-,} ;
(50,0) *+++{Tropical \ Geometry} *\frm{-,}**\dir{-};
\end{xy}\hspace{40pt}
\end{gather}
\bigskip

 To perform idempotent or tropical mathematics   means to perform analysis and/or geometry after replacing  the Archimedean local field $\R$ with the locally compact {\em semi-field} $\rmax$. This latter structure is defined as the locally compact space $\R_+=[0,\infty)$ endowed with the ordinary product and a new idempotent addition $x\plus y=\max\{ x ,y \}$ that appears as the limit case of the conjugates of ordinary addition by the map $x\mapsto x^h$, when $h\to 0$. The semi-field $\rmax$ is trivially isomorphic to the ``schedule algebra" (or ``max-plus algebra'') $(\R\cup\{-\infty\},\max,+)$, by means of the map $x\mapsto \log x$. The reason why we prefer to work with $\rmax$ in this paper (rather than with the ``max-plus'' algebra) is that the structure of $\rmax$  makes the role of the ``Frobenius"  more transparent. On a field of characteristic $p>1$, the arithmetic Frobenius
is given by the ``additive" map $x \mapsto x^p$.  On $\rmax$, the Frobenius flow has the analogous description given by the action $x \mapsto x^\lambda$ ($\lambda \in \R_+^*$)
which is also ``additive" since it is monotonic and hence compatible with $\max\{ x ,y \}$.\vspace{.02in}

It is well known that the classification of local fields of positive characteristic reduces  to the classification of finite fields of positive characteristic. A local field of characteristic $p>1$ is isomorphic to the field
$\F_q((T))$ of formal power series, over a finite extension $\F_q$ of $\F_p$, with finite order pole at $0$.  There is a strong relation between the $p$-adic field $\Q_p$ of characteristic zero and the field $\F_p((T))$ of characteristic $p$. The connection is described by the Ax-Kochen theorem \cite{AK} (that depends upon the continuum hypothesis) which states the following isomorphism of ultraproducts
\begin{equation}\label{ultra}
   \prod_\omega\,\Q_p \sim \prod_\omega\,\F_p((T))
\end{equation}
for any non-trivial ultrafilter $\omega$ on the integers.
 \vspace{.02in}

The local field $\F_p((T))$ is the unique local field of characteristic $p$ with residue field $\F_p$. Then, the Ax-Kochen theorem essentially states that, for sufficiently large $p$, the field $\Q_p$ resembles its simplification  $\F_p((T))$, in which  addition forgets the carry over rule in the description of a $p$-adic number as a power series in $p^n$. The process that allows one to view $\F_p((T))$ as a limit of fields of characteristic zero obtained by adjoining roots of $p$ to the $p$-adic field $\Q_p$ was described in \cite{Krasner}.
The role of the Witt vectors is that to define a process going backwards from the simplified rule of $\F_p((T))$ to the original algebraic law holding for $p$-adic numbers.

In characteristic one, the role of the local field $\F_p((T))$ is played by the local semi-field $\rmax$.

\bigskip
\begin{center}
\begin{tabular}{|c|c|c|}
\hline & & \\
 Characteristic $0$&$\Q_p$ , \ $p> 1$  &\ Archimedean \ $\R$  \\
 & &\\ \hline & & \\
 Characteristic $\neq 0$& $\F_p((T))$ &\ $p=1$,\ \ $\rmax$  \\
& &\\ \hline
\end{tabular}
\end{center}\bigskip

The  process that allows one to view $\rmax$ as a result of a deformation of the usual algebraic structure on real numbers is known as ``dequantization'' and can be described as a semi-classical limit in the following way.  First of all note that, for $\R\ni h>0$  one has
\[
 h\ln(e^{w_1/h}+e^{w_2/h}) \to \max\{w_1,w_2\}\quad\text{as}~h\to 0.
\]
Thus,  the natural map
\[
D_h: \R_+ \to \R_+\qquad D_h(u) = u^h
\]
satisfies the equation
$$
\lim_{h\to 0}D_h(D_h^{-1}(u_1)+D_h^{-1}(u_2))=\max \{u_1,u_2\}.
$$

In the limit $h\to 0$, the usual algebraic rules on $\R_+$ deform to become those of $\rmax$. Moreover, in the limit one also obtains a one-parameter group of automorphisms of $\rmax$
\begin{equation}\label{frobone}
    \vartheta_\lambda(u)=u^\lambda,\qquad \lambda\in\rmax
\end{equation}
which corresponds to the arithmetic Frobenius in characteristic $p>1$.\vspace{.05in}

After introducing the notion of ``perfect" \dioid (\S \ref{subperfect}), our first main result (\cf \S\ref{entropy}) states that one can adapt the construction of the Witt ring to this limit case of characteristic one. This fact also unveils an interesting deep connection with entropy and thermodynamics, while shedding a
new light on the classical Witt construction itself.

 Our starting point is the basic physics formula expressing the free energy in terms of entropy from a variational principle. In its simplest mathematical form, it states that
\begin{equation}\label{free}
   \log(e^a+e^b)=\sup_{x\in [0,1]}S(x)+x a+ (1-x)b
\end{equation}
where $S(x)=- x\log x-(1-x)\log(1-x)$ is the entropy of the partition of $1$ as $1=x+(1-x)$.

In Theorem \ref{deformcharone}, we explain how formula \eqref{free} allows one to reconstruct addition from its
 characteristic one degenerate limit. This process has an analogue in the construction of the Witt ring which reconstructs the $p$-adic numbers from their degenerate characteristic $p$ limit $\F_p((t))$. Interestingly, this result also leads to a reformulation of the Witt construction in complete analogy with the case of characteristic one (\cf Proposition \ref{deform}).
These developments, jointly with the existing body of results in idempotent analysis and tropical algebraic geometry supply a strong and motivated expectation on the existence of a meaningful notion of ``local" mathematics in characteristic one. In fact, they also suggest the
existence of a notion of ``global" mathematics in characteristic one, pointing out to the possibility that the discrete, co-compact embedding of a global field into a locally compact semi-simple non-discrete ring may have an extended version for semi-fields in \dioids.\vspace{.1in}

The second part of the paper reconsiders the study of the ``absolute point" initiated in \cite{Kapranov}. It is important to explain the distinction between $\Spec\B$, where $\B$ is the initial object in the category of \dioids~(\cf\cite{Golan}) and the sought for ``absolute point" which has been traditionally denoted by $\Spec\F_1$. It is generally expected that $\Sp(\F_1)$ sits under $\Spec\Z$ in the sense of \cite{TV}

\begin{gather}
\label{functCFTmap}
 \,\hspace{100pt}\raisetag{-47pt} \xymatrix@C=25pt@R=25pt{
 \Spec\Z\ar[d] &
  \Spec\B\ar[dl]& \\
\Spec\F_1   & \\
}\hspace{25pt}
\end{gather}

 In order to reasonably fulfill the property to be ``absolute", $\Spec\F_1$ should also sit under $\Spec\B$. Here, it is important to stress the fact that $\Spec\B$ does not qualify to be the ``absolute point", since there is no morphism from $\B$ to $\Z$. We refer to \cite{Lescot} for an interesting approach to algebraic geometry over $\B$.

 Theorem \ref{thmsym} singles out the conditions on a bijection $s$ of the set $\K=H\cup \{0\}$ onto itself ($H=$ abelian group), so that $s$ becomes the addition of $1$ in a field whose multiplicative structure is that of the monoid $\K$. These conditions are:\vspace{.05in}

$\bullet$~ $s$ commutes with its conjugates
  for the action of $H$ by multiplication on the monoid $\K$.\vspace{.05in}

$\bullet$~$s(0)\neq 0$.\vspace{.05in}

The degenerate case of characteristic one appears when one drops the requirement that $s$ is a bijection and replaces it with the idempotent condition $s\circ s =s$, so that $s$ becomes a retraction onto its range. The degenerate structure of $\F_1$ arises when one drops the condition $s(0)\neq 0$. The trivial solution given by the identity map $s={\rm id}$, then yields such degenerate structure in which, except for the addition of $0$, the addition is indeterminate as $0/0$. This fact simply means that, except for $0$, one forgets completely the additive structure.
This outcome agrees with the point of view developed in the earlier papers \cite{deit}, \cite{deit1} and \cite{TV}.

 In \cite{announc3}, we have introduced a  geometric theory of algebraic schemes (of finite type) over $\F_1$ which unifies an earlier construction developed in \cite{Soule} with the approach initiated in our paper \cite{ak}. This theory is reconsidered and fully developed in the present article. The geometric objects covered by our construction are rational varieties which admit a suitable cell decomposition in toric varieties. Typical examples  include Chevalley groups schemes, projective spaces etc. Our initial interest for an algebraic geometry over $\F_1$ originated with the study of Chevalley schemes as (affine) algebraic varieties over $\F_1$. Incidentally, we note that these structures are not toric varieties in general and this shows that the point of view of \cite{deit} and \cite{deit1} is a bit too restrictive. Theorem 2.2 of \cite{Lopez} gives a precise characterization of the schemes studied by our theory  in terms of general torifications. This result also points out to a deep relation with tropical geometry that we believe is worth to be explored in details.

 On the other hand, one should  not forget the fact that the class of varieties covered even by this extended definition of schemes over $\F_1$ is still extremely restrictive and, due to the rationality requirement, excludes curves of positive genus. Thus, at first sight, one seems still quite far from the real objective which is that to describe the ``curve" (of infinite genus) $C=\overline {\Sp\,\Z}$ over $\F_1$,  whose zeta function $\zeta_C(s)$ is the complete Riemann zeta function.
However, a thorough development of the theory of schemes over $\F_1$ as initiated in \cite{announc3} shows that any such a scheme $X$  is described (among other data) by a covariant functor $\underline X$  from the category $\Mo$ of commutative monoids to sets. The functor $\underline X$ associates to an object $M$ of $\Mo$ the set of points of $X$ which are ``defined over $M$". The ad\`ele class space of a global field $\K$ is a particularly important example of a monoid and this relevant description allows one to evaluate any scheme over $\F_1$ on such monoid. In \cite{announc3}, we have shown that using the simplest example of a scheme over $\F_1$, namely the curve $\spad$, one constructs a perfect set-up to understand conceptually the spectral realization of the zeros of the Riemann zeta function (and more generally of Hecke $L$-functions over $\K$). The spectral realization appears naturally from the (sheaf) cohomology $H^1(\spad,\Omm)$ of very simply defined sheaves $\Omm$ on the geometric realization of the scheme $\spad$. Such sheaves are obtained from the functions on the points of $\spad$ defined over the monoid $M$ of ad\`ele classes. Moreover, the computation of the cohomology $H^0(\spad,\Omm)$ provides a description of the graph of the Fourier transform which thus appears conceptually in this picture.\vspace{.05in}

 The above construction makes use of a peculiar property of the geometric realization of an $\Mo$-functor which has no analogue for $\Z$-schemes (or more generally for $\Z$-functors). Indeed, the geometric realization of an $\Mo$-functor $\underline X$ coincides, as a set, with the set of  points which are defined over $\F_1$, \ie with the value $\underline X(\F_1)$ of the functor on the most trivial monoid. After a thorough development of the general theory of $\Mo$-schemes
(an essential ingredient in the finer notion of an $\F_1$-scheme), we investigate in full details their geometric realization, which turn out to be  schemes in the sense of \cite{deit} and \cite{deit1}. Then, we show that both the topology and the structure sheaf of an $\F_1$-scheme  can be obtained in a natural manner on the set $\underline X(\F_1)$ of its $\F_1$-points.\vspace{.05in}

At the beginning of \S \ref{repfunfunctor}, we shortly recall our construction of an $\F_1$-scheme and the description of the associated zeta function, under the (local) condition of  no-torsion on the scheme (\cf \cite{announc3}).
We then remove this hypothesis and compute, in particular, the zeta function of the extensions $\F_{1^n}$ of $\F_1$.
 The general result is stated in Theorem \ref{dthmfonesch1} and in Corollary \ref{corzeta} which present a description of
 the zeta function of a Noetherian $\F_1$-scheme $X$
as the product
$$
    \zeta(X,s)=e^{h(s)}\prod_{j=0}^n(s-j)^{\alpha_j}
$$
of the exponential of an entire function by a finite product of fractional powers of simple monomials. The exponents $\alpha_j$ are rational numbers defined explicitly, in terms of the structure sheaf $\cO_X$ in monoids, as follows
$$
    \alpha_j=(-1)^{j+1}\sum_{x\in X}(-1)^{n(x)}{n(x)\choose j}\sum_{d}\frac 1d\nu(d,\cO_{X,x}^\times)
$$
where  $\nu(d,\cO_{X,x}^\times)$ is the number of  injective homomorphisms from $\Z/d\Z$ to the group $\cO_{X,x}^\times$ of invertible elements of the monoid $\cO_{X,x}$.
In order to establish this result, we need to study the case of a  counting function
 \begin{equation}\label{count1}
\# \,X(\F_{1^n})=N(n+1)
\end{equation}
no longer  polynomial in the integer $n\in \N$. In \cite{announc3} we showed that the limit formula that was used in \cite{Soule} to define the zeta function of an algebraic variety over $\F_1$, can be replaced by an equivalent integral formula which determines the equation
\begin{equation}\label{integralrep}
\frac{\partial_s\zeta_N(s)}{\zeta_N(s)}=-\int_1^\infty N(u)u^{-s-1}du
\end{equation}
 describing the logarithmic derivative of the zeta function $\zeta_N(s)$ associated to the counting function $N(q)$. We use this result to treat the case of the counting function of
a Noetherian $\F_1$-scheme and Nevanlinna theory to uniquely extend the counting function $N(n)$ to arbitrary complex arguments $z\in\C$ and finally we  compute the corresponding integrals.

 In \S\ref{ellipticsect} we show that the replacement, in  formula \eqref{integralrep}, of the integral by the discrete sum
$$\frac{\partial_s\zeta_N^{\rm disc}(s)}{\zeta_N^{\rm disc}(s)}=-\sum_1^\infty N(u)u^{-s-1}$$
only modifies the zeta function of a Noetherian $\F_1$-scheme by an exponential factor of an entire function. This observation leads to the Definition~\ref{modifiedzeta} of a modified zeta function over $\F_1$ whose main advantage  is that of being applicable to the case of an arbitrary counting function with polynomial growth.

In \cite{announc3}, we determined the counting distribution $N(q)$, defined for $q\in [1,\infty)$,   such that the  zeta function  \eqref{integralrep} gives the complete Riemann zeta function
$\zeta_\Q(s)=\pi^{-s/2}\Gamma(s/2)\zeta(s)$: \ie such that the following equation holds
\begin{equation}\label{thelogzetabisnew}
    \frac{\partial_s\zeta_\Q(s)}{\zeta_\Q(s)}=-\int_1^\infty  N(u)\, u^{-s}d^*u.
\end{equation}
In \S \ref{riemannzeta},
 we use the modified version of zeta function (as described above) to study a simplified form of \eqref{thelogzetabisnew} \ie
\begin{equation}\label{simplerzeta1}
   \frac{\partial_s\zeta(s)}{\zeta(s)}=-\sum_1^\infty  N(n)\, n^{-s-1}.
\end{equation}
This gives, for the counting function $N(n)$, the formula
$
    N(n)=n\Lambda(n)\,,
$
where  $\Lambda(n)$ is the von-Mangoldt function\footnote{with value $\log p$ for powers $p^\ell$ of  primes and zero otherwise}.
Using the relation \eqref{count1}, this shows that, for a suitably extended notion of ``scheme over $\F_1$", the  hypothetical ``curve" $X=\overline {\Sp\,\Z}$ should fulfill the following requirement
\begin{equation}\label{idealequ}
   \#(X(\F_{1^n}))=\left\{
                              \begin{array}{ll}
                                0 & \hbox{if}\ n+1 \ \hbox{is not a prime power}\\
                                (n+1)\log p, & \hbox{if}\ n=p^\ell-1,\ p \ \hbox{prime.}
                              \end{array}
                            \right.
\end{equation}
This expression is neither functorial nor integer valued, but by reconsidering the theory of additive structures previously developed in \S \ref{addi}, we show in Corollary \ref{corfine}  that the natural construction which assigns to an object $M$ of $\Mo$ the set $A(M)$ of maps $s:M\to M$ such that\vspace{.05in}

$\bullet$~$s$ commutes with its conjugates by multiplication by elements of $M^\times$\vspace{.05in}

$\bullet$~$s(0)=1$\vspace{.05in}

comes close to solving the requirements of \eqref{idealequ} since it gives (for $n>1$, $\varphi$ the Euler totient function)
\begin{equation}\label{cardxmequ}
     \#(A(\F_{1^n}))=\left\{
                              \begin{array}{ll}
                                0 & \hbox{if}\ n+1 \ \hbox{is not a prime power}\\
                                \frac{\varphi(n)}{\log (n+1)}\,\log p, & \hbox{if}\ n=p^\ell-1,\ p \ \hbox{prime.}
                              \end{array}
                            \right.
\end{equation}

In \S\ref{ellipticsubsect} we experiment with elliptic curves $E$ over $\Q$, by computing the zeta function associated to a {\em specific} counting function on $E$. The specific  function $N(q,E)$ is uniquely determined by the following two conditions:\vspace{.05in}

-~For any prime power $q=p^\ell$, the value of $N(q,E)$ is the number\footnote{including the singular point} of points of the reduction of $E$ modulo $p$ in the finite field $\F_q$.\vspace{.05in}

-~The function $t(n)$ occurring in the equation $N(n,E)=n+1-t(n)$ is multiplicative.\vspace{.05in}

Then, we prove that
the obtained zeta function $\zeta_N^{\rm disc}(s)$ of $E$  fulfills the equation
$$
\frac{\partial_s\zeta_N^{\rm disc}(s)}{\zeta_N^{\rm disc}(s)}=-\zeta(s+1)-\zeta(s)+\frac{L(s+1,E)}{\zeta(2s+1)M(s+1)}
$$
where $L(s,E)$ is the $L$-function of the elliptic curve, $\zeta(s)$ is the Riemann zeta function
and
$$
    M(s)=\prod_{p\in S}(1-p^{1-2s})
$$
is a finite product of local factors indexed on the set $S$ of primes at which $E$ has  bad reduction. In Example \ref{concrete} we exhibit the singularities of $\zeta_E(s)$ in a concrete case.

\medskip

\section{Working in characteristic $1$: Entropy and Witt construction}\label{wittsect}

In this  section we investigate the degeneration of the structure of fields of prime positive characteristic $p$, in the limit case $p=1$. A thorough investigation of this process and the algebraic structures of characteristic one, points out to a very interesting link between the fields of idempotent analysis \cite{Maslov} and tropical geometry \cite{Mikhalkin}, with the algebraic structure of a degenerate geometry (of characteristic one) obtained as the limit case of geometries over finite fields $\F_p$.
Our main result is that of adapting the construction of the Witt ring to this limit case of characteristic $p=1$ while unraveling also a deep connection of the Witt construction with entropy and thermodynamics. Remarkably, we find that this process sheds also a new light on the classical Witt construction itself.\vspace{.05in}

\subsection{Additive structure}\hfill\label{addi}\vspace{.05in}

The multiplicative structure of a field $\K$ is obtained by adjoining an absorbing element $0$ to its multiplicative group $H=\K^\times$. Our goal here, is to understand the additional structure on the multiplicative monoid $H\cup \{0\}$ (the multiplicative monoid underlying $\K$) corresponding to the addition. For simplicity, we shall still denote by $\K$ such monoid.

We first notice that to define an additive structure on $\K$ it is enough to know how to define the operation
\begin{equation}\label{ops}
s_o: \K\to\K\qquad  s_o(x)=x+1 \qquad \forall~ x\in \K
\end{equation}
since then the definition of the addition follows by setting
\begin{equation}\label{adds}
+: \K\times\K\to\K\qquad    x+y=x\,s_o(yx^{-1}) \qquad \forall~ x, y\in \K.
\end{equation}
Moreover, with this definition given, one has (using the commutativity of the monoid $\K$)
$$
x(y+z)=x\,y\,s_o(zy^{-1})= xy\,s_o(xzy^{-1}x^{-1})=xy\,s_o(xz(xy)^{-1}).
$$
Thus, the distributivity follows automatically
\begin{equation}\label{distr}
    x(y+z)=xy+xz\qquad\forall~ x,y,z\in \K.
\end{equation}
The following result characterizes the bijections $s$ of the monoid $\K$ onto itself such that, when enriched by the addition law \eqref{adds}, the monoid $\K$ becomes a field.

\begin{thm} \label{thmsym} Let  $H$ be an abelian group.
Let $s$ be a bijection  of the set $\K=H\cup \{0\}$ onto itself that commutes with its conjugates
  for the action of $H$ by multiplication on the monoid $\K$. Then, if $s(0)\neq 0$, the operation
\begin{equation}\label{mainadddefn}
    x+y:=\begin{cases}  y&\text{if}~x=0\\
s(0)^{-1}xs(s(0)yx^{-1})&\text{if}~x\neq 0\end{cases}
\end{equation}
defines a commutative group law on $\K$. With this law as addition, the monoid   $\K$ becomes a commutative field. Moreover the field $\K$ is of characteristic $p$ if and only if
\begin{equation}\label{iter}
   s^p= s\circ s\circ \ldots \circ s={\rm id}
\end{equation}

\end{thm}

\proof By replacing $s$ with its conjugate, one can assume, since $s(0)\neq 0$, that $s(0)=1$.
Let $A(x,y)=x+y$ be given as in \eqref{mainadddefn}. By
definition, one has $A(0,y)=y$. For $x\neq 0$, $A(x,0)=xs(0)=x$. Thus $0$ is
a neutral element for $\K$. The commutation of $s$ with its conjugate for the action of $H$ on $\K$ by multiplication with the element $x\neq 0$ means that
\begin{equation}\label{asshyp}
    s(x\,s(yx^{-1}))=x\,s(s(y)x^{-1}) \qquad \forall  y\in \K.
\end{equation}
Taking $y=0$ (and using $s(0)=1$), the above equation gives
\begin{equation}\label{comcons}
s(x)=x\, s(x^{-1})\qquad\forall x\neq 0.
\end{equation}
Assume now that $x\neq 0$ and $y\neq 0$, then from \eqref{comcons} one gets
$$
A(x,y)=xs(yx^{-1})=x(yx^{-1})s(xy^{-1})=y\,s(xy^{-1})=A(y,x).
$$
Thus $A$ is a commutative law of composition. The associativity of $A$ follows from the commutation of left and right addition which is a
 consequence of the commutation of the conjugates of $s$. More precisely, assuming first that all elements involved are non zero, one has
\begin{align*}
A(A(x,y),z)&=A(xs(yx^{-1}),z)=A(z,xs(yx^{-1}))=zs(xs(yx^{-1})z^{-1})\\&=zs(xz^{-1}\,s(yx^{-1})).
\end{align*}
Using \eqref{asshyp} one obtains
$$
zs(xz^{-1}\,s(yx^{-1}))=z\,xz^{-1}\,s(s(yz^{-1})zx^{-1})=x\,s(s(yz^{-1})zx^{-1})
$$
$$
A(x,A(y,z))=A(x, s(yz^{-1})z)=x\,s(s(yz^{-1})zx^{-1})
$$
which yields the required equality. It remains  to verify a few special cases. If
either $x$ or $y$ or $z$ is zero, the equality follows since $0$ is a neutral element.
Since we never had to divide by $s(a)$, the above argument applies without
restriction. Finally let $\theta=s^{-1}(0)$, then for any $x\neq 0$ one has $A(x,\theta x)=xs(\theta)=0$. This shows
that $\theta x$ is the inverse of $x$ for the law $A$. We have thus proven that this law
defines an abelian group structure on $\K$.

Finally, we claim that the distributive law holds. This means that for any $a\in \K$ one has
$$
A(ax,ay)=aA(x,y)\qqq x,y\in \K.
$$
We can assume that all elements involved are $\neq 0$. Then, one obtains
$$
A(ax,ay)=ax\, s(ax(ay)^{-1})=ax\, s(xy^{-1})=a\, A(x,y).
$$
This suffices to show that $\K$ is a field.
\endproof

\smallskip
As a corollary, one  obtains the following general uniqueness result

\begin{cor}\label{uniqueness} Let $H$ be a finite commutative group and let $s_j$ ($j=1,2$) be two bijections of the monoid $\K=H\cup \{0\}$ fulfilling the conditions of Theorem \ref{thmsym}. Then $H$ is a cyclic group of order $m=p^\ell-1$ for some prime $p$ and there exists  a group automorphism $\alpha\in \Aut(H)$ and an element $g\in H$ such that
$$
 s_2   =T\circ s_1\circ T^{-1}\,, \ T(x)=g\alpha(x)\qqq x\in H\,, \ T(0)=0.
$$
\end{cor}

\proof By replacing $s_j$ with its conjugate by $g_j=s_j(0)$, one can assume that $s_j(0)=1$.
 Each $s_j$ defines on the monoid $\K=H\cup \{0\}$ an additive structure which, in view of  Theorem~\ref{thmsym}, turns this set into a finite field $\F_q$, for some prime power $q=p^\ell$. If $m$ is the order of $H$, then
 $m=p^\ell-1=q-1$ and $H$ is the cyclic group of order $m$. The field $\F_q$ with $q=p^\ell$ elements is unique up to isomorphism.  Thus there exists a bijection $\alpha$ from $\K$ to
itself, which is an isomorphism with respect to the multiplicative  and the
additive structures given by $s_j$. This map sends $0$ to $0$ and $1$ to $1$,
thus transforms the addition of $1$, for the first structure, into the
addition of $1$ for the second one. Since this map also respects the multiplication, it is necessarily given by  a group automorphism $\alpha\in \Aut(H)$.
\endproof\vspace{.05in}

\subsection{Characteristic $p=2$}\hfill \label{chartwo}\vspace{.05in}

 We apply the above discussion to a concrete case. In this subsection we work in characteristic two, thus we consider an algebraic closure $\bar\F_2$.
 It is well known that
the multiplicative group $(\bar\F_2)^\times$ is {\em non-canonically} isomorphic to the group $\mul^{\rm odd}$ of roots of unity in $\C$
of odd order.  We consider the  monoid $\K=\mul^{\rm odd} \cup \{0\}$. The product between non zero elements in $\K$ is the same as in  $\mul^{\rm odd}$ and  $0$ is an absorbing element (\ie $0\cdot x=x\cdot 0=0$, $\forall~x\in \K$). For each positive integer $\ell$, there is a unique subgroup $\mu_{(2^{\ell}-1)}<\mul^{\rm odd}$ made by the roots of $1$ of order $2^\ell-1$; this is also the subgroup of the invertible elements of the finite subfield $\F_{2^\ell}\subset \bar\F_2$.  We know how
to multiply two elements in $\K$ but we do not know how to add them.

Since we work in characteristic two,  the transformation $s_o$ of \eqref{ops} given by the addition of $1$ fulfills $s_o\circ s_o=id$, \ie  $s_o$ is an involution.\vspace{.05in}

Let $\Sigma$ be the set of {\em involutions} $s$ of $\K$ which
commute with all their conjugates by rotations with elements of $\mul^{\rm odd}$,
and which fulfill the condition $s(0)=1$. Thus an element $s\in \Sigma$ is an involution  such
that $s(0)=1$ and it satisfies the identity $s\circ(R\circ s\circ R^{-1})=(R\circ s\circ R^{-1})\circ s$, for any rotation $R: \K\to\K$.

\begin{prop}\label{car2}
 For any  choice of a group isomorphism
\begin{equation}\label{rootembed}
    j:(\bar \F_2)^\times\stackrel{\sim}{\longrightarrow}\mul^{\rm odd},
\end{equation}
the following map defines an element $s\in \Sigma$
\begin{equation}\label{rootinv}
s:\K\to\K\qquad s(x)=\begin{cases}j(j^{-1}(x)+1)& \text{if $x\neq 1$,}
\\0& \text{if $x=1$.}
\end{cases}
\end{equation}
Moreover, two pairs $(\bar \F_2,j)$ are isomorphic if and only if the associated symmetries
  $ s\in \Sigma$ are the same.
Each element of $\Sigma$ corresponds to an uniquely associated  pair $(\bar \F_2,j)$.
\end{prop}

 \begin{figure}
\begin{center}
\includegraphics[scale=0.7]{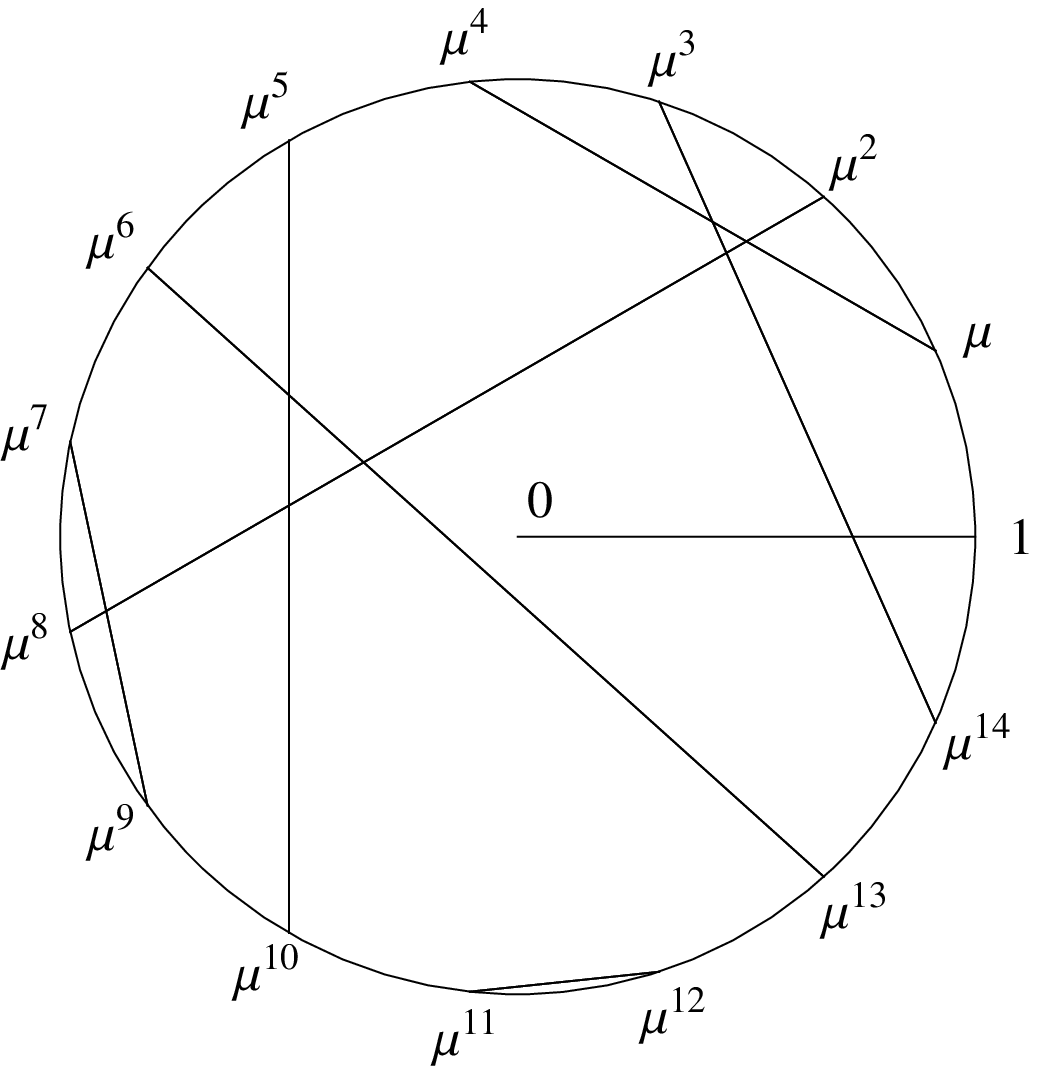}
\end{center}
\caption{The symmetry $s$ for $\F_{16}$ represented by the graph $G_s$
with vertices the $15$-th roots of $1$ and edges $(x,s(x))$.\label{safe7draw} }
\end{figure}

\begin{figure}
\begin{center}
\includegraphics[scale=0.7]{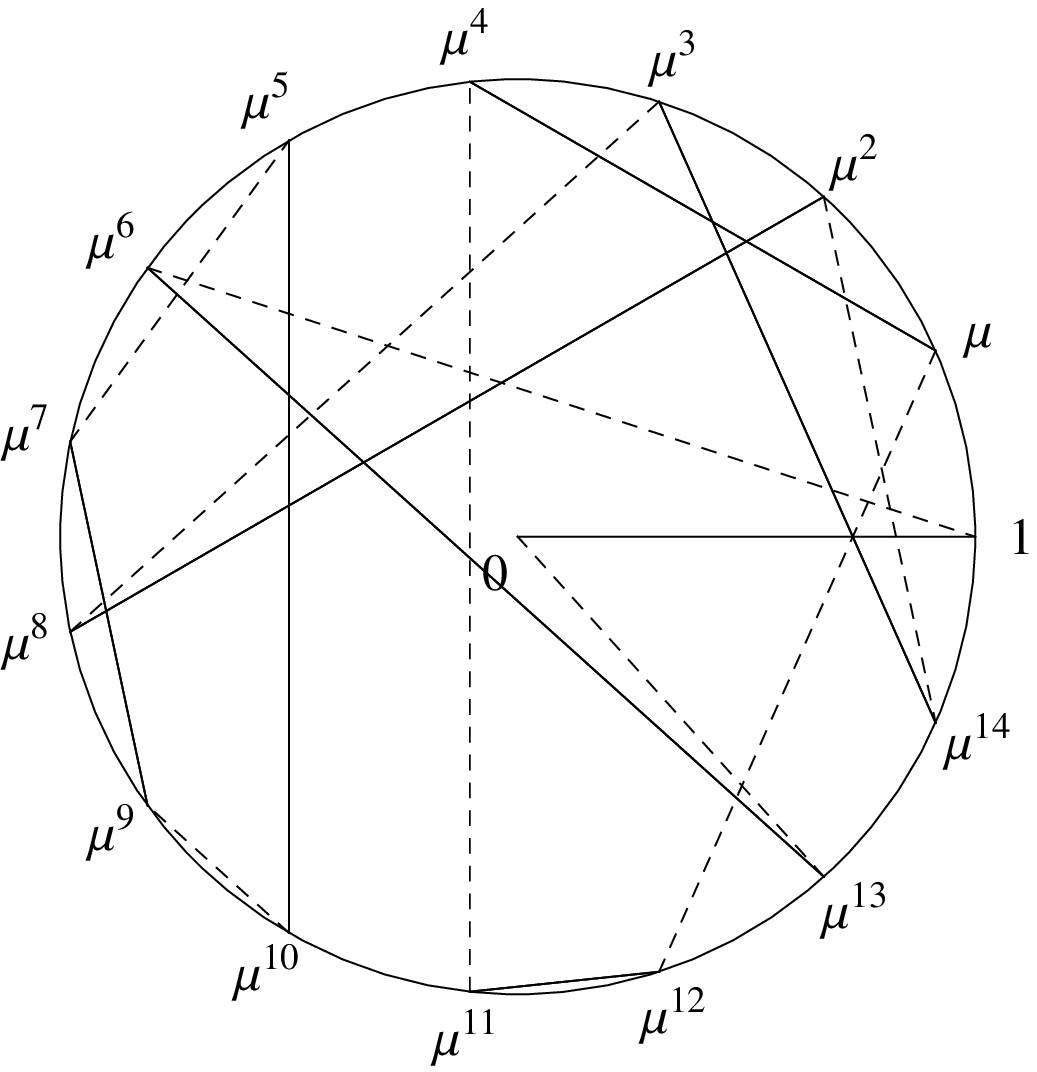}
\end{center}
\caption{The commutativity of $s$ with its conjugates by rotations $R\circ s\circ R^{-1}$ means that the union of the graph $G_s$ with a rotated graph $R(G_s)$ (dashed)  is a union of quadrilaterals.\label{safe8draw} }
\end{figure}
\proof We first show that $s\in \Sigma$. Since ${\rm char}(\bar \F_2)=2$, $s$ is an involution. To show that $s$ commutes with its conjugates
by rotations with elements of $\mul^{\rm odd}$, it is enough to check that $s$ commutes with its conjugate by
$j(y)$,  $y\in(\bar \F_2)^\times$. One first uses the distributivity of the addition  to see that
$$
j(y)s(xj(y)^{-1})=j(y)j(j^{-1}(xj(y)^{-1})+1)=j(y(j^{-1}(xj(y)^{-1})+1))=j(j^{-1}(x)+y).
$$
The commutativity of the additive structure on $\bar \F_2$ gives the required
commutation.

We now prove the second statement. Given two pairs $(\K_i,j_i)$ ($i=1,2$) as in
\eqref{rootembed}, and an isomorphism $\theta:\K_1\to \K_2$ such that $j_2\circ
\theta=j_1$, one has
$$
s_2(x)=j_2(j_2^{-1}(x)+1)=j_1\circ \theta^{-1}(\theta(j_1^{-1}(x))+1)=s_1(x)
$$
since $\theta(1)=1$.
\endproof

Let $s\in \Sigma$ be an involution of $\K$ satisfying the required properties.  We introduce the following operation in $\K$
\begin{equation}\label{adddefn1}
    x+y=\begin{cases}  y&\text{if}~x=0\\
xs(yx^{-1})&\text{if}~x\neq 0.\end{cases}
\end{equation}
Then, it follows as a corollary of  Theorem~\ref{thmsym} that the above operation \eqref{adddefn1} defines a commutative group law on $\K=\mul^{\rm odd}\cup \{0\}\subset
\C$. This law and the induced multiplication, turn $\K$ into a field of
characteristic $2$.

\begin{rem}\label{conway}{\rm
We avoid to refer to $\bar\F_p$ as to `the' algebraic closure of $\F_p$, since for $q=p^n$ a prime power, the finite
field  $\F_q$ with $q$ elements is well defined only up to isomorphism.
In
computations, such as the construction of tables of modular characters,
one uses an explicit construction of $\F_q$ as the quotient ring  $\F_p[T]/(C_n(T))$,
where $C_n(T)$ is a specific irreducible polynomial, called Conway's
polynomial. This polynomial is of degree $n$ and fulfills simple algebraic conditions and a
normalization property involving the lexicographic ordering (\cf \eg
\cite{Luebeck}).
In the particular case of characteristic $2$, J. Conway was able to
produce {\em canonically}  an algebraic closure $\bar\F_2$ using inductively defined operations on the ordinals\footnote{The first author is grateful to Javier Lopez for pointing out this construction of J. Conway.} less than $\omega^{\omega^\omega}$ (\cf\cite{Conway}). His construction also provides one with a
natural choice of the isomorphism
$$
j:(\bar \F_2)^\times\stackrel{\sim}{\longrightarrow}\mul^{\rm odd}.
$$
 In fact, one can use the well ordering to choose the smallest solution as a root of unity of order $2^\ell-1$
fulfilling compatibility conditions with the previous choices for divisors of $\ell$. One obtains this way a well-defined corresponding symmetry $s$ on $\mul^{\rm odd}\cup \{0\}$. Figure \ref{safe7draw} shows the restriction of this symmetry to $\mul^{\rm 15}\cup \{0\}$, while  Figure \ref{safe8draw} gives the simple geometric meaning of the commutation of $s$ with its conjugates by rotations.
}\end{rem}

\subsection{Characteristic $p=1$ and idempotent analysis}\hfill \label{charone}\vspace{.05in}

 To explain how idempotent analysis and the semifield $\rmax$ appear naturally in the above framework, we shall no longer require that the map $s$, which represents the ``addition of $1$" is a bijection of $\K$, but we shall still retain the property that $s$ is a {\em retraction} \ie it satisfies the idempotent condition
 \begin{equation}\label{retract}
    s\circ s=s.
 \end{equation}
 Before  stating the analogue of Theorem \ref{thmsym} in this  context we will review a few definitions.

 The classical theory of rings has been generalized by the more general theory of semi-rings (\cf \cite{Golan}).

\begin{defn} \label{semi-ring } A semi-ring  $(R,+,\cdot)$ is a non-empty set $R$ endowed with operations of addition and multiplication such that\vspace{.05in}

$(a)$~$(R,+)$ is a commutative monoid with neutral element $0$\vspace{.05in}

$(b)$~$(R,\cdot)$ is a monoid with identity $1$\vspace{.05in}

$(c)$~$\forall r,s,t\in R$:~~$r(s+t) = rs+rt$ and $(s+t)r = st+sr$\vspace{.05in}

$(d)$~$\forall r\in R$:~~$r\cdot 0=0\cdot r =0$\vspace{.05in}

$(e)$~$0\neq 1$.

\end{defn}

To any semi-ring $R$ one associates a characteristic.  The set $\N1_R = \{n1_R|n\in\N\}$ is a commutative subsemi-ring  of $R$ called the prime semi-ring. The prime semi-ring is the smallest semi-ring contained in $R$. If $n\neq m \Rightarrow n1_R\neq m1_R$, then $\N1 \simeq \N$ naturally and the semi-ring $R$ is said to have characteristic zero. On the other hand, if it happens that $k1_R=0_R$, for some $k>1$, then there is a least positive integer $n$ with $n1_R = 0_R$. In this case,  $\N1\simeq\Z/n\Z$ and one shows that $R$ itself is a ring of positive characteristic. Finally, if $\ell 1_R = m1_R$ for some $\ell\neq m$, $\ell, m \ge 1$ but $k1_R \neq 0$ for all $ k\ge 1$, one writes $j$ for $j1_R$ and it follows that  $\N1_R = \{0,1,\ldots,n-1\}$. In this case one has (\cf \cite{Golan} Proposition 9.7)

 \begin{prop}\label{propprime} Let $n$ be the least positive integer such that $n1_R \in\{ i1_R\mid 1\le i\le n-1\}$ and let
 $i\in \{1,n-1\}$ with $n1_R=i1_R$. Then the prime semi-ring  $\N1_R = \{n1_R|n\in\N\}$ is the following semi-ring  $B(n,i)$
\[
B(n,i) := \{0,1,\ldots,n-1\}
\]
with the following operations, where $m = n-i$
\[
x+'y := \begin{cases} x+y &\text{if $x+y\le  n-1$},\\
\ell \,, i\le\ell\le n-1\,, \ell = (x+y)\,{\rm mod}\,m &\text{if $x+y\ge n$}
\end{cases}
\]
\[
x\cdot y := \begin{cases} xy & \text{if $xy\le n-1$},\\
\ell \,, i\le\ell\le n-1\,, \ell = (xy)\,{\rm mod}\,m &\text{if $xy\ge n.$}
\end{cases}
\]
\end{prop}
 The semi-ring  $B(n,i)$ is the homomorphic image of $(\N,+,\cdot)$ by the map $\pi: \N \to R$, $\pi(k) = k$ for $0\le k\le n-1$ and for $k\ge n$, $\pi(k)$ is the unique natural number congruent to $k$ mod $m = n-i$ with $i\le \pi(k)\le n-1$. In this case we say that $R$ has characteristic $(n,i)$.\vspace{.05in}

Note that if $R$ is a semi-field the only possibility for the semi-ring  $B(n,i)$ is when $n=2$ and $i=1$.  Indeed the subset $\{\ell \,, i\le\ell\le n-1\}$ is stable under product and hence should contain $1$ since a finite submonoid of an abelian group is a subgroup. Thus $i=1$, but in that case one gets $(n-1)(n-1)=(n-1)$ which is possible only for $n=2$.
 The smallest finite (prime) semi-ring structure arises when $n=2$ (and  $i=1$). We shall denote this structure by $\B = B(2,1)$. By construction, $\B = \{0,1\}$ with the usual multiplication law and an addition requiring the idempotent rule $1+1=1$.

\begin{defn} \label{charisone} A semi-ring  $R$ is said to have characteristic $1$ when
\begin{equation}\label{idemplus}
1+ 1=1
\end{equation}
in $R$ \ie when $R$ contains $\B$ as the prime sub-semi-ring.
\end{defn}
When $R$ is a semi-ring of characteristic $1$, we denote the addition  in $R$ by the symbol $\plus$.
Then, it follows from distributivity  that
\begin{equation}\label{idemplusbis}
a\plus a=a \qqq a\in R.
\end{equation}
 This justifies the term ``additively idempotent" frequently used in semi-ring theory as a synonymous for ``characteristic one". A semi-ring  $R$ is called a {\em semi-field} when every non-zero element in $R$ has a multiplicative inverse, or equivalently when the set of non-zero elements in $R$ is a (commutative) group for the multiplicative law.

\begin{thm} \label{thmsymnil} Let  $H$ be an abelian group.
Let $s$ be a retraction ($s\circ s =s$) of the set $\K=H\cup \{0\}$  that commutes with its conjugates
  for the action of $H$ by multiplication on the monoid $\K$. Then, if $s(0)\neq 0$, the operation
\begin{equation}\label{adddefn}
    x+y=\begin{cases}  y&\text{if}~x=0\\
s(0)^{-1}xs(s(0)yx^{-1})&\text{if}~x\neq 0\end{cases}
\end{equation}
defines a commutative monoid law on $\K$. With this law as addition, the monoid   $\K$ becomes a commutative  semi-field of characteristic $1$.
\end{thm}

\proof The proof of Theorem \ref{thmsym} applies without modification. Notice that that proof did not use the hypothesis that $s$ is a bijection, except to get the element $\theta=s^{-1}(0)$ which was used to define the additive inverse. The fact that $s$ is a retraction shows that $K$ is of characteristic $1$. \endproof

\begin{example}  \label{rplusstar}{\rm Let $H=\R_+^*$ be the multiplicative group of the positive real numbers.
Let $s$ be the retraction of $\K=H\cup \{0\}=\R_+$ on $[1,\infty)\subset H$  given by
$$
s(x)=1\qqq x\leq 1\,, \ s(x)=x\qqq x \geq 1.
$$
 The conjugate $s^\lambda$ under multiplication (by $\lambda$) is the retraction on $[\lambda,\infty)$ and one easily checks that $s$ commutes with $s^\lambda$. The resulting commutative idempotent semi-field is denoted by $\rmax$. Thus $\rmax$ is the set $\R_+$ endowed with the following two operations\vspace{.05in}

$\bullet$~Addition
  \begin{equation}\label{addr}
x\plus y:={\rm max}(x,y) \qqq x,y\in \R_+
\end{equation}
$\bullet$~Multiplication
   is unchanged.
}\end{example}

\begin{example}  \label{groupsubsets}{\rm  Let $H$ be a group and let $2^H$ be the set of subsets of $H$ endowed with the following two operations\vspace{.05in}

$\bullet$~Addition
  \begin{equation}\label{add}
X\plus Y:=X\cup Y,\qquad \forall X,Y\in 2^H
\end{equation}
$\bullet$~Multiplication
  \begin{equation}\label{mult}
X.Y=\{ab \mid a\in X, b\in Y\} \qqq X,Y\in 2^H.
\end{equation}\vspace{.05in}

Addition is commutative and associative and admits the empty set as a neutral element. Multiplication is
associative and it admits the empty set
as an absorbing element (which we denote by $0$ since it is also the neutral element
for the additive structure). The multiplication is also distributive with respect to the addition.
}\end{example}

\begin{example} \label{monotonicmaps}{\rm Let $(M,\plus )$ be an idempotent semigroup. Then, one endows the set $\End(M)$ of endomorphisms  $h: M\to M$ such that
$$
h(a \plus  b)=h(a)\plus  h(b)\qqq a,b\in M
$$
with the following operations (\cf \cite{Golan}, I Example 1.14),\vspace{.05in}

$\bullet$~Addition
  \begin{equation}\label{add1}
(h\plus  g)(a)=h(a)\plus  g(a)\qqq a\in M,\ \forall h,g\in \End(M)
\end{equation}
$\bullet$ Multiplication
  \begin{equation}\label{mult1}
(h\cdot g)(a)=h(g(a)),\qquad \forall a\in M,\ \forall h,g\in \End(M).
\end{equation}

For instance, if one lets $(M,\plus )$ be the idempotent semigroup given by $(\R,{\rm max})$, the set $\End(M)$ of endomorphisms of $M$ is the set of monotonic maps $\R\to\R$. $\End(M)$ becomes a \dioid with the above operations.

}\end{example}

\begin{example} \label{ccstar}{\rm

Let $\C_{\rm star}$ be the set of finitely generated star shaped subsets of the complex numbers $\C$. Thus an element of $\C_{\rm star}$ is of the form
$$
S=\{0\}\bigcup_{z\in F}\{\lambda z\mid z\in F, \ \lambda \in [0,1]\}
$$
for some finite subset $F\subset \C$. One has a natural injection $\C_{\rm star}\to 2^{\C^*}$ given by
$$
S\mapsto S\cap \C^*.
$$
The image of $\C_{\rm star}$ through this injection is stable under the semi-ring  operations \eqref{add} and \eqref{mult} in the semi-ring  $2^{\C^*}$, where we view $\C^*$ as a multiplicative group. This shows that $\C_{\rm star}$ is an \dio under the following operations\vspace{.05in}

$\bullet$~Addition
  \begin{equation}\label{add2}
X\plus Y:=X\cup Y \qqq X,Y\in \C_{\rm star}
\end{equation}
$\bullet$~Multiplication
  \begin{equation}\label{mult2}
X.Y=\{ab \mid a\in X, b\in Y\} \qqq X,Y\in \C_{\rm star}
\end{equation}
}\end{example}\vspace{.05in}

\subsubsection{Finite semi-field of characteristic $1$}

We quote the following result from \cite{Golan} (Example 4.28 Chapter 4).
\begin{prop}\label{uniquesemi} The semi-field $\B=B(2,1)$ is the only finite semi-field of characteristic $1$.
\end{prop}\vspace{.05in}

The semi-field $\B=B(2,1)$ is called the Boolean semi-field (\cf \cite{Golan}, I Example 1.5).

\subsubsection{Lattices}

The idempotency of addition in semi-rings of characteristic $1$ gives rise to a natural {\em partial order} which differentiates this theory from that of the more general semi-rings. We recall the following well known fact

\begin{prop}\label{lattic} Let $(A,\plus)$ be a commutative semigroup with an idempotent addition. Define
\begin{equation}\label{deforder}
a\preccurlyeq b\quad\Leftrightarrow\quad a\plus b = b
\end{equation}
Then $(A,\preccurlyeq)$ is a sup-semilattice (\ie a semilattice in which any two elements have a supremum). Furthermore
\begin{equation}\label{plus}
max\{a,b\} = a\plus b.
\end{equation}
Conversely, if $(A,\preccurlyeq)$ is a sup-semilattice and $\plus$ is defined to satisfy \eqref{plus}, then $(A,\plus)$ is an idempotent semigroup. These two constructions are inverse to each other.
\end{prop}
\proof We check that \eqref{deforder} defines a partial order $\preccurlyeq$ on $A$. Let $a,b,c\in A$. The reflexive property ($a\le a$) follows from the idempotency of the addition in $A$ \ie $a\plus a = a$. The antisymmetric property ($a\le b$, $b\le a$ $\Rightarrow$ $a=b$) follows from the commutativity of the addition in $A$ \ie $b =a\plus b = b\plus a = a$. Finally, the transitivity property ($a\le b$, $b\le c$ $\Rightarrow$ $a\le c$) follows from the associativity of the binary operation $\plus$ \ie $a\plus c = a\plus(b\plus c) = (a\plus b)\plus c = b\plus c=c$. Thus $(A,\preccurlyeq)$ is a poset and due to the idempotency of the addition, $(A,\preccurlyeq)$ is also a semilattice.  One defines the join of two elements in $A$ as: $a\vee b := a\plus b$ and  due to the closure property of the law $\plus$ in $A$ (\ie $a\plus b\in A$, $\forall a,b\in A$), one concludes that $(A,\preccurlyeq)$ is a sup-semilattice (the supremum$=$ maximum of two elements in $A$ being their join).

The converse statement follows too since the above statements on the idempotency and associativity of the operation $\plus$ hold also in reverse and the closure property derives from \eqref{plus}.\endproof

It follows that any \dioid  has a natural ordering, denoted $\preccurlyeq$. The addition in the semi-ring is a monotonic operation and $0$ is the least element. Distributivity implies that left and right multiplication are semilattice homomorphisms and in particular they are monotonic. Thus, any \dioid may be thought of as a semilattice ordered semigroup.\vspace{.05in}

\subsubsection{Perfect semi-rings of characteristic $1$}\label{subperfect}

 The following notion is the natural extension to semi-rings of the absence of zero divisors in a ring.

\begin{defn}\label{cancell}
A commutative semi-ring is called multiplicatively-cancellative when the multiplication by any non-zero element is injective.
\end{defn}

We recall, from \cite{Golan} Propositions 4.43 and 4.44   the following result which describes the analogue of the Frobenius endomorphism in characteristic $p$.

\begin{prop}\label{pow} Let $R$ be a multiplicatively-cancellative commutative \dioid, then for any integer $n\in \N$, the map $x\mapsto x^n$ is an injective endomorphism of $R$.
\end{prop}

Recall now that in characteristic $p>1$ a ring is called {\em perfect} if the map $x\mapsto x^p$ is surjective. Let $R$ be a multiplicatively-cancellative commutative \dioid, then if the endomorphism $R\to R,~x\mapsto x^n$  is surjective, then it is bijective and one can invert it and construct the fractional powers
$$
\vartheta_\alpha: R \to R,\qquad \vartheta_\alpha(x)=x^\alpha\qqq \alpha\in \Q_+^*.
$$
Then, by construction, the $\vartheta_\alpha$'s are  automorphisms of $R$ and they fulfill the following properties
 $$
 \vartheta_\lambda\circ \vartheta_\mu=\vartheta_{\lambda \mu}\qqq \lambda,\mu \in \Q_+^*
 $$
 and
 $$
 \vartheta_\lambda(x)\vartheta_\mu(x)=\vartheta_{\lambda+\mu}(x)\qqq \lambda,\mu \in \Q_+^*\,, \ x\in R.
 $$
 Moreover, under natural completeness requirements one can extend this construction by passing from $\Q_+^*$ to its completion $\R_+^*$. Thus we shall adopt the following notion

\begin{defn}\label{perfect}
A commutative \dioid is called {\em perfect} when there exists a one parameter group of automorphisms
$\vartheta_\lambda\in \Aut(R)$, for $\lambda \in \R_+^*$, such that\vspace{.05in}

$\bullet$~$\vartheta_n(x)=x^n$ for all $n\in \N$ and $x\in R$.\vspace{.05in}

$\bullet$~$\vartheta_\lambda\circ \vartheta_\mu=\vartheta_{\lambda \mu}$ for all $\lambda,\mu \in \R_+^*$.\vspace{.05in}

$\bullet$~$\vartheta_\lambda(x)\vartheta_\mu(x)=\vartheta_{\lambda+\mu}(x)$ for all $\lambda,\mu \in \R_+^*$ and $x\in R.$
\end{defn}\vspace{.05in}

\subsubsection{Localization}\label{localsemi}

One can localize a commutative \dioid with respect to any multiplicative subset $S$ not containing $0$ (\cf\cite{Golan} Proposition 11.5). Moreover, when $S$ is made by multiplicatively cancellative elements the natural morphism $R\to S^{-1}R$ is injective. We shall apply this construction in the case $S=\{\rho^n\mid n\in \N\}$ for a multiplicatively cancellative element $\rho\in R$.\vspace{.05in}

By analogy with the standard notation adopted for rings, we denote by $R_\rho$ the semi-ring  $S^{-1}R$  for $S=\{\rho^n\mid n\in \N\}$.

The following result is straightforward.

\begin{prop}\label{local} Let $R$ be a perfect commutative \dioid and let $\rho\in R$ be a multiplicatively cancellative element. Then the localization $R_\rho$ is a perfect commutative \dioid.
\end{prop}
\proof It is immediate to check that $R_\rho$ is a commutative semi-ring of characteristic one. To verify that $R_\rho$ is also perfect, it is enough to show that an automorphism $\vartheta_\lambda\in \Aut(R)$ ($\lambda \in \R_+^*$) verifying the conditions of Definition~\ref{perfect} extends uniquely to an automorphism $\vartheta'_\lambda\in \Aut(R_\rho)$, when $\rho\in R$ is a multiplicatively cancellative element. Indeed, one defines $\vartheta'_\lambda(\frac{r}{\rho^n}) = \vartheta_\lambda(r)\cdot\vartheta_\lambda(\rho^n)^{-1}$. It belongs to $R_\rho$ because for $m>n\lambda$ one can replace $\vartheta_\lambda(\rho^n)^{-1}$ by $\vartheta_{m-n\lambda}(\rho)/\rho^m$. It is now straightforward to verify that $\vartheta'_\lambda$ verifies the properties of Definition~\ref{perfect}.
\endproof\vspace{.05in}

\subsection{Witt ring in characteristic $p=1$ and entropy}\hfill \label{entropy}\vspace{.05in}

The places of the global field $\Q$ of the rational numbers fall in two classes: the infinite archimedean place $\infty$ and the finite places which are labeled by the prime integer numbers. The $p$-adic completion of $\Q$ at a finite place $p$ determines the corresponding global field $\Q_p$ of $p$-adic numbers. These local fields have close relatives with simpler structure: the local fields $\F_p((T))$ of formal power series with coefficients in the finite fields $\F_p$ and with finite order pole at $0$. We already explained  in the introduction that the similarity between the structures of $\Q_p$ and of $\F_p((T))$ is embodied in the Ax-Kochen Theorem (\cf\cite{AK}) which states the isomorphism of arbitrary ultraproducts as in \eqref{ultra}. By means of the natural construction of the ring of Witt vectors, one recovers the ring $\Z_p$ of $p$-adic integers from the finite field $\F_p$. This construction (\cf also \S \ref{pbig}) can be interpreted as a deformation of the ring of formal power series $\F_p[[T]]$ to $\Z_p$.

It is then natural to wonder on the existence of a similar phenomenon at the infinite archimedean  place of $\Q$. We have already introduced  the  semi-field of characteristic one $\rmax$ as the degenerate version of the field  of real numbers. In this subsection we shall describe how the basic physics formula for the free energy involving entropy allows one to move canonically, not only from $\rmax$ to $\R$ but in even greater generality from a perfect \dioid to an ordinary algebra over $\R$. This construction is also in perfect analogy with the construction of the Witt ring.\vspace{.05in}

 Let   $K$ be a perfect \dioid, and let us first assume that it contains
 $\rmax$ as a sub semi-ring. Thus the operation of raising to a power $s\in \R_+$ is well defined in $K$ and  determines an automorphism $\vartheta_s: K\to K,~ \vartheta_s(x)=x^s$. We shall assume in this section that we can use  the idempotent integral $\int^\plus$ (\cf \cite{Maslov} I, \S 1.4) to integrate the functions from $[0,1]$ to $K$ involved below and will concentrate on the algebraic aspects leaving the technical aspects aside. The intuitive way of thinking of the idempotent integral $\int^\plus_Y f(s)$ is as the least upper bound of the range $f(Y)\subset K$, this least upper bound is assumed to exist under suitable compactness conditions.
 We can then make sense of the following formula
 \begin{equation}\label{quant}
    x+'y:=\int^\plus_{s\in [0,1]}c(s)\,x^s\,y^{1-s}\qqq x,y\in K
 \end{equation}
 where the function $c(s)\in \R_+$ is defined by
 \begin{equation}\label{cx}
 c:[0,1]\to\R_+,\quad   c(s)=e^{S(s)}=s^{-s}(1-s)^{-(1-s)}.
 \end{equation}

  \begin{figure}
\begin{center}
\includegraphics[scale=0.8]{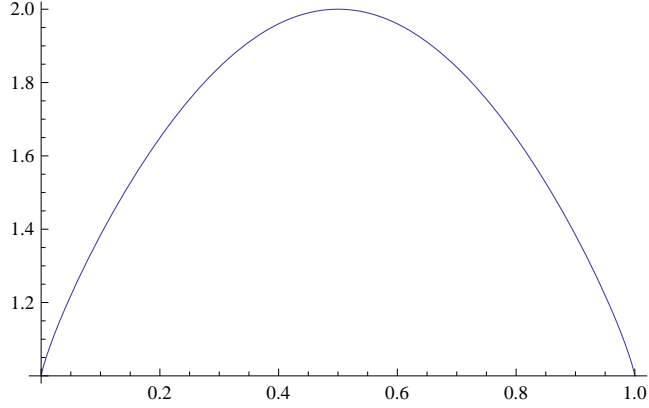}
\end{center}
\caption{Graph of $c(x)$.\label{fbinone} }
\end{figure}

The property of the function $c(s)$ which is at the root of the associativity
of this addition law is the following: for any $s,t\in [0,1]$ the product $c(s)c(t)^s$ only depends upon the
partition of $1$ as $st+s(1-t)+(1-s)$. This fact also implies the functional equation
\begin{equation}\label{funequa}
    c(u)c(v)^u=c(uv)c(w)^{(1-uv)}\,, \ \ w=\frac{u(1-v)}{1-uv}.
\end{equation}
 The function $c$ fulfils the symmetry $c(1-u)=c(u)$ and, by taking it into account, \eqref{funequa} means then that the function
on the simplex $\Sigma_2=\{(s_j)\mid s_j\geq 0,\sum_0^2 s_j=1\}$ defined by
\begin{equation}\label{funequa1}
    c_2(s_0,s_1,s_2)=c(s_2)c(\frac{s_0}{s_0+s_1})^{s_0+s_1}
\end{equation}
is symmetric under all permutations of the $s_j$.
More generally, for any integer $n$ one may define the function on the $n$-simplex
\begin{equation}\label{funequan}
    c_n(s_0,s_1,\ldots, s_n)=\prod \gamma(s_j)^{-1}\,, \ \ \gamma(x)=x^x.
\end{equation}

\begin{lem}\label{freeenerg} Let $x,y\in \R_+$, then one has
\begin{equation}\label{freeenerg1}
    \sup_{s\in [0,1]}c(s)\,x^s\,y^{1-s}=x+y.
\end{equation}
\end{lem}

\proof
Let $x=e^a$ and $y=e^b$. The function $$f(s)=\log(c(s)\,x^s\,y^{1-s})=S(s)+s a+(1-s)b$$
is strictly convex on the interval $[0,1]$ and reaches its unique maximum for
$s=\frac{e^a}{e^a+e^b}\in [0,1]$. Its value at the maximum is $\log(e^a+e^b)$.
\endproof

\begin{thm}\label{deformcharone} Let $K$ be a perfect \dioid. Let $\rho\in K$, $\rho \preccurlyeq 1$ be a multiplicatively cancellative element and let $K_\rho$ be the localized semi-ring (\cf \S \ref{localsemi}). Then the formula
\begin{equation}\label{quantprime}
    x+_\rho y:=\int^\plus_{s\in [0,1]}\rho^{-S(s)}\,x^s\,y^{1-s}\qqq x,y\in K_\rho
 \end{equation}
defines an associative law on  $K_\rho $ with $0$ as neutral element. The multiplication is distributive with respect to this law. These operations turn the Grothendieck group of the additive monoid $(K_\rho,+_\rho)$ into an algebra $W(K,\rho)$ over $\R$ which depends functorially upon $(K,\rho)$.
\end{thm}

\proof Since $\rho \preccurlyeq 1$ one has, using the automorphisms $\vartheta_\lambda$, that
$$
\vartheta_\lambda(\rho)=\rho^\lambda \preccurlyeq 1 \qqq \lambda \in \R_+.
$$
Then it follows that
 $$
 \rho^{\lambda} \preccurlyeq \rho^{\lambda'}\qqq \lambda \geq \lambda'.
 $$
 Thus the following map defines a homomorphism of semi-rings
\begin{equation}\label{homom}
\alpha:(\R_+,{\rm max},\cdot )=\rmax\to K_\rho,\quad  \alpha(\lambda)=\rho^{-\log \lambda}\,, \ \alpha(0)=0\,,
\end{equation}
where we have used the invertibility of $\rho$ in $K_\rho$ to make sense of the negative powers of $\rho$. Let $c(t)=\rho^{-S(t)}$. The associativity of the operation $+_\rho$ follows  from
$$
(x+_\rho y)+_\rho z=\int^\plus_{s\in [0,1]}c(s)\,\vartheta_s\left(\int^\plus_{t\in [0,1]}c(t)\,x^t\,y^{1-t}\right)\,z^{1-s}=
$$
$$
\int^\plus_{s\in [0,1]}\int^\plus_{t\in [0,1]}c(s)\,c(t)^s\,x^{ts}\,y^{(1-t)s}\,z^{1-s}=\int^\plus_{\Sigma_2}c_2(s_0,s_1,s_2)\,x^{s_0}\,y^{s_1}\,z^{s_2}
$$
which is symmetric in $(x,y,z)$ by making use of \eqref{funequa1} and \eqref{funequan}. Moreover, one has by homogeneity the distributivity
$$
(x+_\rho y)z=xz+_\rho  yz\qqq x,y,z\in K\,.
$$
We let $K_\rho $ be the semi-ring with the operations $+_\rho $ and the multiplication unchanged.
If we endow $\R_+$ with its ordinary addition, we have by applying Lemma \ref{freeenerg}, that
$$
\alpha(\lambda_1)+_\rho \alpha(\lambda_2)=\alpha(\lambda_1+\lambda_2).
$$
  Thus $\alpha$ defines a homomorphism
  \begin{equation}\label{homomtilde}
\tilde \alpha:\R_+\to (K_\rho,+_\rho,\cdot),\quad  \tilde\alpha(\lambda)=\rho^{-\log \lambda}\,, \ \alpha(0)=0\,,
\end{equation}
 of the semi-ring $\R_+$ (with ordinary addition) into the semi-ring $(K_\rho,+_\rho,\cdot)$.

Let $G$ be the functor from semi-rings to rings which associates to a semi-ring $R$
the Grothendieck group  of the additive monoid $R$ with the natural extension of the product. One can view the ring $G(R)$ as the quotient of the semi-ring $R[\Z/2\Z]$ by  the equivalence relation
\begin{equation}\label{equiv}
(m_1, m_2)\sim (n_1, n_2)\Leftrightarrow \exists k\in R \,, \ m_1 +   n_2 +   k = m_2 +   n_1 +   k.
\end{equation}
 The addition is given coordinate-wise and the product is defined by
 $$
 (m_1, m_2)\cdot (n_1, n_2)=(m_1n_1+  m_2n_2,m_1n_2+  m_2n_1).
 $$
The equivalence relation is compatible with the product which turns the quotient $G(R)$ into a ring. We define
\begin{equation}\label{defnW}
    W(K,\rho)=G((K_\rho,+_\rho,\cdot))
\end{equation}
One has $G(\R_+)=\R$.
By functoriality of $G$ one thus gets a morphism $G(\tilde \alpha)$ from $\R$ to $W(K,\rho)$. As long as $1\neq 0$ in $G((K_\rho,+_\rho,\cdot))$ this morphism endows $W(K,\rho)$ with the structure of an algebra over $\R$. When $1=0$ in $G((K_\rho,+_\rho,\cdot))$ one gets the degenerate case $W(K,\rho)=\{0\}$.
 \endproof

 We did not discuss conditions on $\rho$ which ensure that $(K_\rho,+_\rho,\cdot)$ injects in $W(K,\rho)$, let alone that $W(K,\rho)\neq\{0\}$. The following example shows that the problem comes from how strict the inequality $\rho \preccurlyeq 1$ is assumed to be, \ie where the function $T(x)$ playing the role of absolute temperature actually vanishes.

\begin{example}\label{beta}{\rm Let $X$ be a compact space and let $R=C(X,[0,1])$ be the space of continuous maps from $X$ to the interval $[0,1]$. We endow $R$ with the operations ${\rm max}$ for addition and the ordinary pointwise product for multiplication. The associated semi-ring is a perfect \dioid. Let then $\rho\in R$ be a non-vanishing map then it is of the form
$$
\rho(x)=e^{-T(x)}\qqq x\in X\,.
$$
Let $\beta(x)=1/T(x)$ if $T(x)>0$ and $\beta(x)=\infty$ for $T(x)=0$. Then, the addition $+_\rho$ in $R$ is given by
\begin{equation}\label{betaform}
   (f+_\rho g)(x)=\left(f(x)^{\beta(x)}+g(x)^{\beta(x)}\right)^{T(x)}\qqq x\in X\,.
\end{equation}
This follows from Lemma \ref{freeenerg} which implies more generally that
$$
\sup_{s\in [0,1]}e^{TS(s)}\,x^s\,y^{1-s}=(x^\beta+y^\beta)^T\,, \ \beta=\frac 1T\,.
$$
Then, provided that $\rho(x)<1$ for all $x\in X$, the algebra $W(R,\rho)$ is isomorphic,
as an algebra over $\R$, with the real $C^*$-algebra $C(X,\R)$ of continuous functions on $X$.
}\end{example}\vspace{.05in}

\subsection{The Witt ring in characteristic $p>1$ revisited}\hfill \label{pbig}\vspace{.05in}

In this subsection we explain in which sense we interpret the formula \eqref{quantprime} as the analogue of the construction of the Witt ring in characteristic one. Let $K$ be a perfect ring of characteristic $p$. We start by reformulating the construction of the Witt ring in characteristic $p>1$. One knows  (\cf \cite{Serre} Chapter II) that there is a strict $p$-ring $R$, unique up to canonical isomorphism, such that its residue ring $R/pR$ is equal to $K$. One also knows that there exists a unique multiplicative section $\tau:K\to R$ of the residue map. For $x\in K$, $\tau(x)\in R$ is called the Teichm\"{u}ller representative of $x$. Every element  $z\in R$ can be  written uniquely as
\begin{equation}\label{binwitt0}
z=\sum_0^\infty \tau(x_n)p^n\,.
\end{equation}
The Witt construction of $R$ is functorial and an easy corollary of its properties is the following

\begin{table}\label{w5}
$$
\begin{array}{cccccccccc}
 \frac{1}{125} & 4 T^3 & \frac{2}{125} & 3 T^3 & \frac{3}{125} & 3 T^3 & \frac{4}{125} & 4 T^3 & \frac{1}{25} & 4 T^2 \\
 \frac{6}{125} & 3 T^3 & \frac{7}{125} & 2 T^3 & \frac{8}{125} & 0 & \frac{9}{125} & T^3 & \frac{2}{25} & 3 T^2 \\
 \frac{11}{125} & 4 T^3 & \frac{12}{125} & 2 T^3 & \frac{13}{125} & 2 T^3 & \frac{14}{125} & 4 T^3 & \frac{3}{25} & 3 T^2+2 T^3 \\
 \frac{16}{125} & 0 & \frac{17}{125} & 0 & \frac{18}{125} & 3 T^3 & \frac{19}{125} & 3 T^3 & \frac{4}{25} & 4 T^2+2 T^3 \\
 \frac{21}{125} & T^3 & \frac{22}{125} & 3 T^3 & \frac{23}{125} & 0 & \frac{24}{125} & 3 T^3 & \frac{1}{5} & 4 T \\
 \frac{26}{125} & 3 T^3 & \frac{27}{125} & 0 & \frac{28}{125} & 2 T^3 & \frac{29}{125} & T^3 & \frac{6}{25} & 3 T^2+3 T^3 \\
 \frac{31}{125} & T^3 & \frac{32}{125} & 2 T^3 & \frac{33}{125} & 4 T^3 & \frac{34}{125} & 0 & \frac{7}{25} & 2 T^2+T^3 \\
 \frac{36}{125} & 3 T^3 & \frac{37}{125} & T^3 & \frac{38}{125} & 0 & \frac{39}{125} & T^3 & \frac{8}{25} & T^3 \\
 \frac{41}{125} & 2 T^3 & \frac{42}{125} & 2 T^3 & \frac{43}{125} & 0 & \frac{44}{125} & 3 T^3 & \frac{9}{25} & T^2+T^3 \\
 \frac{46}{125} & 3 T^3 & \frac{47}{125} & 2 T^3 & \frac{48}{125} & 4 T^3 & \frac{49}{125} & 2 T^3 & \frac{2}{5} & 3 T \\
 \frac{51}{125} & T^3 & \frac{52}{125} & 3 T^3 & \frac{53}{125} & 4 T^3 & \frac{54}{125} & T^3 & \frac{11}{25} & 4 T^3 \\
 \frac{56}{125} & 2 T^3 & \frac{57}{125} & T^3 & \frac{58}{125} & T^3 & \frac{59}{125} & T^3 & \frac{12}{25} & 4 T^2+T^3 \\
 \frac{61}{125} & 4 T^3 & \frac{62}{125} & T^3 & \frac{63}{125} & T^3 & \frac{64}{125} & 4 T^3 & \frac{13}{25} & 4 T^2+T^3 \\
 \frac{66}{125} & T^3 & \frac{67}{125} & T^3 & \frac{68}{125} & T^3 & \frac{69}{125} & 2 T^3 & \frac{14}{25} & 4 T^3 \\
 \frac{71}{125} & T^3 & \frac{72}{125} & 4 T^3 & \frac{73}{125} & 3 T^3 & \frac{74}{125} & T^3 & \frac{3}{5} & 3 T \\
 \frac{76}{125} & 2 T^3 & \frac{77}{125} & 4 T^3 & \frac{78}{125} & 2 T^3 & \frac{79}{125} & 3 T^3 & \frac{16}{25} & T^2+T^3 \\
 \frac{81}{125} & 3 T^3 & \frac{82}{125} & 0 & \frac{83}{125} & 2 T^3 & \frac{84}{125} & 2 T^3 & \frac{17}{25} & T^3 \\
 \frac{86}{125} & T^3 & \frac{87}{125} & 0 & \frac{88}{125} & T^3 & \frac{89}{125} & 3 T^3 & \frac{18}{25} & 2 T^2+T^3 \\
 \frac{91}{125} & 0 & \frac{92}{125} & 4 T^3 & \frac{93}{125} & 2 T^3 & \frac{94}{125} & T^3 & \frac{19}{25} & 3 T^2+3 T^3 \\
 \frac{96}{125} & T^3 & \frac{97}{125} & 2 T^3 & \frac{98}{125} & 0 & \frac{99}{125} & 3 T^3 & \frac{4}{5} & 4 T \\
 \frac{101}{125} & 3 T^3 & \frac{102}{125} & 0 & \frac{103}{125} & 3 T^3 & \frac{104}{125} & T^3 & \frac{21}{25} & 4 T^2+2 T^3 \\
 \frac{106}{125} & 3 T^3 & \frac{107}{125} & 3 T^3 & \frac{108}{125} & 0 & \frac{109}{125} & 0 & \frac{22}{25} & 3 T^2+2 T^3 \\
 \frac{111}{125} & 4 T^3 & \frac{112}{125} & 2 T^3 & \frac{113}{125} & 2 T^3 & \frac{114}{125} & 4 T^3 & \frac{23}{25} & 3 T^2 \\
 \frac{116}{125} & T^3 & \frac{117}{125} & 0 & \frac{118}{125} & 2 T^3 & \frac{119}{125} & 3 T^3 & \frac{24}{25} & 4 T^2 \\
 \frac{121}{125} & 4 T^3 & \frac{122}{125} & 3 T^3 & \frac{123}{125} & 3 T^3 & \frac{124}{125} & 4 T^3 & 1 & 1
\end{array}
$$
\caption{Values of $w_5(\alpha)$ modulo $T^4$.}
\end{table}

\begin{lem} \label{wittcons} For any prime number $p$, there exists a universal sequence
$$
w(p^n,k)\in \Z/p\Z,\quad \ 0<k<p^n
$$
such that
\begin{equation}\label{binwitt1}
    \tau(x)+\tau(y)=\tau(x+y)+\sum_{n=1}^\infty\tau\left(\sum w(p^n,k)\,x^{\frac{k}{p^n}}y^{1-\frac{k}{p^n}}\right )p^n.
\end{equation}
\end{lem}

Note incidentally that the fractional powers such as $x^{\frac{k}{p^n}}$ make sense in a perfect ring such as $K$. The main point here is that  formula \eqref{binwitt1} suffices to reconstruct the whole ring structure on $R$ and allows one to add and multiply series of the form \eqref{binwitt0}.

\proof Formula \eqref{binwitt1} is a special case of the basic formula
$$
s_n=S_n(x_0^{1/p^n},y_0^{1/p^n},x_1^{1/p^{n-1}},y_1^{1/p^{n-1}},\ldots, x_n,y_n)
$$
(\cf\cite{Serre}, Chapter II \S 6) which computes the sum
$$
\sum_0^\infty \tau(x_n)p^n+\sum_0^\infty \tau(y_n)p^n=\sum_0^\infty \tau(s_n)p^n
$$
\endproof

We can now introduce a natural map from the set $I_p$ of rational numbers in $[0,1]$ whose denominator is a power of $p$, to the maximal compact subring of the local field $\F_p((T))$ as follows
\begin{equation}\label{wmap}
w_p: I_p \to \F_p((T)),\quad  w_p(\alpha)=\sum_{\frac{a}{p^n}=\alpha}w(p^n,a)T^n\in \F_p((T)).
\end{equation}
We can then rewrite equation \eqref{binwitt1} in a more suggestive form as a deformation of the addition, by first introducing the map
$$
\tau: K[[T]]\to R,\quad \tau\left(\sum_0^\infty z_n T^n\right )=\sum_0^\infty \tau(z_n)p^n.
$$
This map is an homeomorphism and is multiplicative on monomials \ie
$$
\tau(aT^n Z)=\tau(a)\tau(T)^n\tau(Z)\qqq a\in K,\ Z\in K[[T]].
$$
Since $\tau$ is not additive, one defines a new addition on $K[[T]]$ by setting
$$
X+' Y:=\tau^{-1}(\tau(X)+\tau(Y))\qqq X,Y\in K[[T]].
$$

\begin{prop}\label{deform} Let us view the ring of formal series $K[[T]]$ as a module over $\F_p[[T]]$. Then one has
\begin{equation}\label{binwitt}
    x+' y:=\sum_{\alpha\in I_p} w_p(\alpha)\, x^\alpha\,y^{1-\alpha}\qqq x,y\in K.
\end{equation}

\end{prop}

This formula is  in perfect analogy with formula \eqref{quantprime} of Theorem \ref{deformcharone} and suggests to develop in a deeper way the properties of the $K[[T]]$-valued function $w_p(\alpha)$ in analogy with the classical properties of the entropy function $w_1(\alpha)=e^{S(\alpha)}$ taking its values in $\rmax$.\vspace{.05in}

\subsection{$\B$, $\F_1$ and the ``absolute point"}\hfill \label{degcase}\vspace{.05in}

The \dioids $\,$ provide a natural framework for a mathematics of finite characteristic $p=1$. Moreover, the semi-field  $\B$ is the initial object among \dioids. However, $\Spec \B$ does not fulfill  the requirements of the ``absolute point" $\Spec\F_1$, as defined in \cite{Kapranov}. In particular, one expects that $\Spec\F_1$  sits under $\Spec\Z$. This property does not hold for $\Spec \B$ since there is no unital homomorphism of semi-rings from $\B$ to $\Z$.

We conclude this section by explaining how the ``naive" approach to $\F_1$ emerges in the framework of \S \ref{addi} and \S \ref{chartwo}. First, notice that
Proposition \ref{car2} generalizes to characteristic $p>2$ by implementing the following  modifications\vspace{.05in}

$\bullet$~The group $\mul^{\rm odd}$ is replaced by the group $\mul^{(p)}$ of roots of $1$ in $\C$ of order prime to $p$.\vspace{.05in}

$\bullet$~The involution  $s^2=id$ is replaced by a bijection $s:\K\to\K$ of $\K =\mul^{(p)}\cup\{0\}$, satisfying the condition $s^p=id$ and commuting with its conjugates by rotations.\vspace{.05in}

To treat the degenerate case $p=1$ in \S \ref{charone} we dropped the condition that $s$ is a bijection and we replaced it by the idempotency condition $s\circ s=s$. There is however another trivial possibility which is that to leave unaltered the condition  $s^p=id$ and simply take $s=id$ for $p=1$. The limit case is then obtained  by implementing the following procedure\vspace{.05in}

$\bullet$~The group $\mul^{(p)}$ is replaced by the group $\mul$ of all roots of $1$ in $\C$\vspace{.05in}

$\bullet$~ The map $s$ is the identity map on $\mul\cup\{0\}$.\vspace{.05in}

The additive structure on
\begin{equation}\label{fone}
    \F_{1^\infty}=\mul \cup \{0\}
\end{equation}
  degenerates
since this case corresponds to setting  $s(0)=0$ in  Theorem~\ref{thmsym}, with the bijection $s$ given by the identity map, so that \eqref{adddefn} becomes the indeterminate expression $0/0$.

The multiplicative structure on the monoid \eqref{fone} is the same as the multiplicative structure on the group $\mul$, where we consider $\{0\}$ as an absorbing element.
By construction, for each integer $n$ the group $\mul$ contains a unique cyclic subgroup $\mul_n$ of order $n$. The tower (inductive limit) of the finite subfields  $\F_{p^n}\subset\bar\F_p$ is replaced in the limit case $p=1$ by the inductive limit of commutative monoids
\begin{equation}\label{fonebis}
    \F_{1^n}=\mul_n \cup \{0\}\subset \mul \cup \{0\}
\end{equation}
where the inductive structure is partially ordered by divisibility of the index $n$.

Notice that on $\F_{1^\infty}$ one can still define
$$
0+x=x+0=x\qquad\forall x\in\F_{1^\infty}
$$
and this simple rule suffices  to multiply matrices with coefficients in $\F_{1^\infty}$ whose rows have at most one non-zero element. In fact, the multiplicative formula $C_{ik}=\sum A_{ij}B_{jk}$ for the product $C=AB$ of two matrices only makes sense if one can add $0+x=x+0=x$. Hence, with the exception of the special case where $0$ is one of the two summands, one considers the addition $x+y$ in $\F_{1^\infty}$ to be indeterminate.\vspace{.05in}

By construction, $\F_1=\mul_1\cup\{0\}=\{1\}\cup\{0\}$ is the monoid to which $\F_q$ degenerates when $q=1$, \ie by considering the addition $1+1$ on $\F_1$ to be indeterminate. From this simple point of view the category of commutative algebras over $\F_1$ is simply the category $\Mo$ of commutative monoids with a unit $1$ and an absorbing element $0$. In particular we see that  a (commutative) \dioid is in particular a (commutative) algebra over $\F_1$ in the above sense which is compatible with the absolute property of $\Spec\F_1$ of \eqref{functCFTmap}.

\begin{rem}\label{careful}{\rm
By construction, each finite field $\F_{p^n}$ has the same multiplicative structure as  the  monoid $\F_{1^{p^n-1}}$. However, there exist infinitely many monoidal structures $\F_{1^n}$ which do not correspond to any degeneration of a finite field: $\F_{1^5}=\mul_5\cup\{0\}$ and $\F_{1^9}=\mul_9\cup\{0\}$ are the first two cases. For the sake of clarity, we make it clear that even when $n=p^\ell-1$ for $p$ a prime number, we shall always consider the algebraic structures $\F_{1^n}= \mul_n \cup \{0\}$ only as multiplicative monoids.
}\end{rem}

\section{The functorial approach}

In this chapter we give an overview on the  geometric theory of algebraic schemes over $\F_1$ that we have introduced in our paper \cite{announc3}. The second part of the chapter contains a new result on the geometric realization of an $\Mo$-functor. In fact, our latest development of the study of the algebraic geometry of the $\F_1$-schemes shows that, unlike the case of a $\Z$-scheme, the topology and the structure sheaf of an $\F_1$-scheme can be obtained naturally on the set of its rational $\F_1$-points.

Our original viewpoint in the development of this theory of schemes over $\F_1$ has been an attempt at unifying the theories developed on the one side by C.~Soul\'e  in \cite{Soule} and in our earlier paper \cite{ak} and on the other side by A.~Deitmar in \cite{deit}, \cite{deit1}  (following  N.~Kurokawa,  H.~Ochiai and M.~Wakayama \cite{KOW}), by K.~Kato in \cite{Kato} (with the geometry of logarithmic structures) and by B.~T\"oen and  M.~Vaqui\'e in \cite{TV}. In the earlier \S \ref{degcase}, we have described how to obtain a naive version of $\F_1$ leading naturally to the point of view developed by A.~Deitmar, where $\F_1$-algebras are commutative monoids (with the slight difference that in our setup one also requires the existence of an absorbing element $0$). It is the analysis performed by C.~Soul\'e  in \cite{Soule} of the extension of scalars from $\F_1$ to $\Z$ that lead us to\vspace{.05in}

$\bullet$~Reformulate the notion of a scheme in the sense of  K.~Kato and A.~Deitmar, in functorial terms \ie as a covariant functor from the category of (pointed) monoids to the category of sets.\vspace{.05in}

$\bullet$~Prove a new result on the geometric realization of functors satisfying a suitable local representability condition.\vspace{.05in}

$\bullet$~Refine the notion of an $\F_1$-scheme by allowing more freedom on the choice of the $\Z$-scheme obtained by base change.\vspace{.05in}

In this chapter we shall explain this viewpoint in some details, focussing in particular on the description and the properties of the geometric realization of an $\Mo$-functor that was only briefly sketched in \cite{announc3}.

Everywhere in this chapter we denote by $\Se$, $\Mo$, $\Ab$, $\An$ respectively the categories of sets, commutative monoids\footnote{with a unit and a zero element}, abelian groups and commutative (unital) rings.\vspace{.05in}

\subsection{Schemes as locally representable $\Z$-functors: a review.}\hfill \label{locfunctor}\vspace{.05in}

  In the following three subsections we will shortly review the basic notions of the theory of $\Z$-functors and $\Z$-schemes: we refer to \cite{demgab} (Chapters I, III) for a detailed exposition.\vspace{.05in}

   \begin{defn} \label{Zfunc} A $\Z$-functor is a covariant functor  $\cF:\An\to\Se$.

 \end{defn}

 Morphisms in the category of $\Z$-functors are natural transformations (of functors).\vspace{.05in}

  Schemes over $\Z$ determine a full subcategory $\Sch$  of the category of $\Z$-functors. In fact, a scheme $X$ over $\Spec \Z$ is {\em entirely} characterized by the $\Z$-functor
 \begin{equation}\label{sctofun}
 \underline X: \An\to \Se,\qquad  \underline  X(R)=\Hom_{\Z}(\Spec R, X).
\end{equation}
 To a ring homomorphism $\rho: R_1\to R_2$ one associates the morphism of (affine) schemes $\rho^*: \Spec(R_2)\to \Spec(R_1)$, $\rho^*(\ffp) = \rho^{-1}(\ffp)$ and the  map of sets
 \[
 \rho:\underline X(R_1) \to \underline X(R_2),\qquad \varphi\to \varphi\circ\rho^*.
 \]

  If $\psi: X \to Y$ is a morphism of schemes then one gets for every ring $R$ a map of sets
 \[
 \underline X(R)\to\underline Y(R),\quad\varphi\to \psi\circ\varphi.
 \]
  The   functors of the form $\underline X$, for $X$ a scheme over $\Spec \Z$,
are referred to as {\em local} $\Z$-functors (in the sense that we shall recall in  \S~\ref{localZfunctors}, Definition~\ref{locfuncdef}). These functors are also {\em locally representable} by commutative rings, \ie they have an open cover by representable $\Z$-subfunctors (in the sense explained in \S\S~\ref{opensubsub},~\ref{cov}).\vspace{.05in}

\subsubsection{Local $\Z$-functors}\hfill \label{localZfunctors}\vspace{.05in}

 For any commutative (unital) ring $R$, the geometric space $\Spec(R)$ is the set of prime ideals $\ffp\subset R$. The topology on $\Spec(R)$ is the Jacobson topology \ie the closed subsets are the  sets $V(J) = \{\ffp\in\Spec(R)|\ffp\supset J\}$, where $J\subset R$ runs through the collection of all the ideals of $R$.
The open subsets of $\Spec(R)$ are the complements of the $V(J)$'s \ie they are the sets
\[
D(J) = V(J)^c= \{\ffp\in\Spec(R)|\exists f\in J, f\notin\ffp\}.
\]
 It is well known that the open sets $D(f) = V(fR)^c$, for $f\in R$, form a base of the topology of $\Spec(R)$. For any $f\in R$ one lets $R_f=S^{-1}R$, where $S=\{f^n|n\in\Z_{\ge 0}\}$ denotes the multiplicative system of the non-negative powers of $f$. One has a natural ring homomorphism $R\to R_f$.  Then, for any scheme  $X$ over $\Spec \Z$, the associated functor
$\underline X$ as in  \eqref{sctofun} fulfills the following {\em locality} property.
For any finite cover of $\Spec(R)$ by open sets $D(f_i)$, with $f_i\in R$ ($i\in I=$ finite index set) the following sequence of maps of sets is exact
\begin{equation}\label{locfunc0}
    \underline X(R)\stackrel{u}{\longrightarrow}\prod_{i\in I} \underline X(R_{f_i})\;\xy
{\ar@{->}_{w} (0,-1)*{}; (6,-1)*{}};
{\ar@{->}^{v} (0,1)*{}; (6,1)*{}};
\endxy\prod_{(i,j)\in I\times I}
    \underline X(R_{f_if_j}).
\end{equation}
 This means that $u$ is injective and the range of $u$ is characterized as the set $\{z\in \prod_i \underline X(R_{f_i})|v(z)=w(z)\}$. The exactness of \eqref{locfunc0} is a consequence of the fact that a morphism of schemes is defined by  local conditions. For $\Z$-functors  we have the following definition

\begin{defn} \label{locfuncdef} A $\Z$-functor $\cF$ is {\em local} if for any object $R$ of $\An$ and a partition of unity $\sum_{i\in I} h_if_i=1$ in $R$ ($I=$ finite index set), the following  sequence of sets is exact:
\begin{equation}\label{locfunc}
    \cF(R)\stackrel{u}{\longrightarrow}\prod_i\cF(R_{f_i})\;\xy
{\ar@{->}_{w} (0,-1)*{}; (6,-1)*{}};
{\ar@{->}^{v} (0,1)*{}; (6,1)*{}};
\endxy\prod_{ij}
    \cF(R_{f_if_j}).
\end{equation}
\end{defn}

\begin{example}\label{exproj}{\rm  This example shows that locality is not automatically verified by a general $\Z$-functor.

The Grassmannian $\text{Gr}(k,n)$ of the $k$-dimensional linear spaces in an $n$-dimensional linear space  is defined by the functor which associates to a ring $R$ the set of all complemented submodules of {\em rank} $k$ of the free (right) module $R^n$. Since any  such complemented submodule is projective, by construction we have
\[
{\text Gr}(k,n): \An\to \Se,\qquad \text{Gr}(k,n)(R) = \{E\subset R^n|~E~ \text{projective},~\text{rk}(E) = k\}.
\]
  Let $\rho: R_1\to R_2$ be a homomorphism of rings, the corresponding map of sets is given as follows: for $E_1\in\text{Gr}(k,n)(R_1)$, one lets $E_2 = E_1\otimes_{R_1}R_2\in \text{Gr}(k,n)(R_2)$. If one takes a naive definition of the rank, \ie by just requiring that $E\cong R^k$ as an $R$-module, one does not obtain a local $\Z$-functor. In fact, let us consider the case $k=1$ and $n=2$ which defines the projective line $\P^1= {\text Gr}(1,2)$. To show that locality fails in this case, one takes the algebra $R = C(S^2)$ of continuous functions on the sphere $S^2$ and the partition of unity $f_1+f_2=1$ subordinate to a covering of $S^2$ by two disks $D_j$ ($j=1,2$), so that $\text{Supp}(f_j)\subset D_j$. One then considers the non-trivial line bundle $L$ on $S^2$ arising from the identification $S^2 \simeq \P^1(\C)\subset M_2(\C)$ which determines an idempotent $e\in M_2(C(S^2))$. The range of $e$ defines a finite projective submodule $E\subset R^2$. The localized algebra $R_{f_j}$ is the same as $C(\bar D_j)_{f_j}$ and thus the induced module $E_{f_j}$ is free (of rank one). The modules $E_i = E\otimes_R R_{f_i}$ are free submodules of $R_{f_i}^2$ and the induced modules on $R_{f_if_j}$ are the same. But since $E$ is not free they do not belong to the image of $u$ and the sequence \eqref{locfunc} is not exact in this case.\vspace{.05in}

To obtain a local $\Z$-functor one has to implement a more refined definition of the rank which requires that for any prime ideal $\ffp$ of $R$ the induced module on the residue field of $R$ at $\ffp$ is a vector space of dimension $k$.
}\end{example}\vspace{.05in}

 \subsubsection{Open $\Z$-subfunctors}\hfill \label{opensubsub}\vspace{.05in}

 The $\Z$-functor associated to an affine scheme $\Spec A$ ($A\in{\rm obj}(\An)$) is defined by
\begin{equation}\label{sctofunbis}
\ssp A: \An\to \Se,\quad \ssp A(R)=\Hom_{\Z}(\Spec R, \Spec A)\simeq \Hom_{\Z}(A,R).
\end{equation}
The open sets of an affine scheme are in general {\em not} affine and they provide interesting examples of schemes. The subfunctor of $\ssp  A$
 \begin{equation}\label{defonsubfun}
 \Hom_{\Z}(\Spec R, D(J))\subset \Hom_{\Z}(\Spec R, \Spec A)
\end{equation}
 associated to the open set $D(J)\subset\Spec(A)$ (for a given ideal $J\subset A$) has the following explicit description
\begin{equation}\label{condition}
\underline D(J): \An \to \Se,\qquad \underline D(J)(R)=\{\rho\in\text{Hom}_\Z(A,R)|\rho(J)R = R\}.
\end{equation}
This follows from the fact that $\Sp(\rho)^{-1}(D(J))=D(\rho(J)R)$. In general, we say that $\mathcal U$ is a subfunctor of a functor $\cF: \An\to \Se$ if for each ring $R$, $\mathcal U(R)$ is a subset of $\cF(R)$ (with the natural compatibility for the maps).

\begin{defn} \label{openfunc} Let $\mathcal U$ be a subfunctor of $\cF: \An \to \Se$. One says that $\cU$ is {\em open} if for any ring $A$ and any natural transformation $\phi: \ssp  A\to\cF$, the subfunctor of $\ssp A$ inverse image of $\mathcal U$ by $\phi$
\[
\phi^{-1}(\mathcal U): \An \to \Se,\quad \phi^{-1}(\mathcal U)(R) = \{x\in\ssp A(R)|\phi(x)\in\mathcal U(R)\subset\cF(R)\}
\]
is of the form $\underline D(J)$, for some {\em open} set $D(J)\subset\Spec(A)$.
\end{defn}
Equivalently, using Yoneda's Lemma, the above definition can be expressed by saying that, given any ring $A$ and an element $z\in \cF(A)$, there exists an ideal $J\subset A$ such that, for any $\rho\in \Hom(A,R)$
\begin{equation}\label{criteropen}
\cF(\rho)(z)\in \mathcal U(R)\iff \rho(J)R = R\,.
\end{equation}
For any open subset $Y\subset X$ of a scheme $X$ the subfunctor
$$
\Hom_{\Z}(\Spec R, Y)\subset \Hom_{\Z}(\Spec R, X)
$$
is open, and all open subfunctors of $\underline X$ arise in this way.

\begin{example}\label{exproj1}{\rm We consider the projective line $\P^1$ identified with $X=\text{Gr}(1,2)$. Let  $\mathcal U\subset X$ be the subfunctor described, on a ring $R$, by the collection of all submodules of rank one of $R^2$ which are supplements of the submodule $P=\{(0,y)|y\in R\}\subset R^2$. Let $p_1$ be the   projection on the first copy of $R$, then:
\[
\mathcal U(R) = \{E\in\text{Gr}(1,2)|\;{p_1}_{|E}~\text{isomorphism}\}.
\]
 Proving that $\mathcal U$ is open is equivalent, using \eqref{criteropen}, to find for any ring $A$ and $E\in   X(A)$ an ideal $J\subset A$ such that for any ring $R$ and $\rho: A \to R$ one has
\begin{equation}\label{proj}
{p_1}_{|E\otimes_A R}~\text{isomorphism}\iff R = \rho(J)R.
\end{equation}
It is easy to see that the ideal $J$ given by the annihilator of the cokernel of ${p_1}_{|E}$ satisfies \eqref{proj} (\cf \cite{demgab} Chapter I, Example~3.9).
}\end{example}\vspace{.05in}

  \subsubsection{Covering by $\Z$-subfunctors}\hfill \label{cov}\vspace{.05in}

 To motivate the definition of a covering of a $\Z$-functor, we start by describing the case of an affine scheme. Let $\underline X$  be the $\Z$-functor
\begin{equation}\label{zfun}
\underline X: \An\to \Se, \quad \underline X(R)=\Hom(A,R)
\end{equation}
 that is associated to the affine scheme $\Spec(A)$. We have seen that the open subfunctors $\underline D(I)\subset \underline X$ correspond to ideals $I\subset A$ with
 \[
 \underline D(I)(R) = \{\rho\in\Hom(A,R)|\rho(I)R = R
 \}.
 \]
 The condition that the open sets $D(I_\alpha)$ ($\alpha\in I=$ index set) form a covering of $\Spec(A)$ is expressed algebraically by the equality $\sum_{\alpha\in I} I_\alpha = A$. We want to describe this condition in terms of the open subfunctors $\underline D(I_\alpha)$.

 \begin{lem} \label{coveringsch} Let $\underline X$ be as in \eqref{zfun}. Then $\sum_\alpha I_{\alpha\in I} = A$ if and only if for any field $K$ one has $$\underline X(K) = \bigcup_{\alpha\in I} \underline D(I_\alpha)(K).$$
 \end{lem}

 \proof Assume first that $\sum I_{\alpha\in I} = A$ (a finite sum) \ie $1 = \sum_\alpha a_\alpha$, with $a_\alpha\in I_\alpha$. Let $K$ be a field, then for $\rho\in\Hom(A,K)$, one has $\rho(a_\alpha) \neq 0$ for some $\alpha\in I$. Then $\rho(a_\alpha)K = K$, \ie $\rho\in \underline D(I_\alpha)(K)$ so that the union of all $\underline D(I_\alpha)(K)$ is  $ \underline X(K)$.

 Conversely, assume that $\sum_\alpha I_{\alpha} \neq A$. Then there exists a prime ideal $\ffp\subset A$ containing all $I_\alpha$'s. Let $K$  be the field of fractions of $A/\ffp$ and  let $\rho: A \to K$ be the natural homomorphism. One has $I_\alpha\subset\text{Ker}\rho$, thus $\rho\notin\bigcup_\alpha\underline D(I_\alpha)(K)$.
 \endproof

Notice that  when $R$ is neither a field or a local ring,  the equality $\underline X(R) = \bigcup_\alpha \underline D(I_\alpha)(R)$ {\em cannot} be expected. In fact the range of a morphism   $\rho\in\Hom(\Spec(R),\Sp A)=\underline X(R)$ may  not be contained in a single open set of the covering of $\Sp(A)$ by the $D(I_\alpha)$ so that $\rho$ belongs to none of the $\underline D(I_\alpha)(R)=\Hom(\Spec(R),D(I_\alpha))$.

\begin{defn}\label{defcover} Let $\underline X$ be a $\Z$-functor. Let $\{{\underline X_\alpha}\}_{\alpha\in I}$ be a family of open subfunctors of  $\underline X$. Then, we say that the set $\{{\underline X_\alpha}\}_{\alpha\in I}$ form a covering of $\underline X$ if for any field $K$ one has $$\underline X(K) = \bigcup_{\alpha\in I} \underline X_\alpha(K).$$
\end{defn}

For affine schemes, one recovers the usual notion of an open cover. In fact, any open cover of an affine scheme admits a finite subcover. Indeed, the condition is $\sum_\alpha I_\alpha = A$ and if it holds one gets $1 = \sum_{\alpha\in F}a_\alpha$ for some {\em finite} subset of indices $F\subset I$. For an arbitrary scheme this finiteness condition may {\em not} hold. However, since a scheme is always ``locally affine'', one can say, calling ``quasi-compact'' the above finiteness condition, that any scheme is locally quasi-compact.

To conclude this short review of the basic properties of schemes viewed as $\Z$-functors, we quote the main theorem which allows one to consider a scheme as a local and locally representable $\Z$-functor.

\begin{thm}\label{charsch}
 The   $\Z$-functors of the form $\underline X$, for $X$ a scheme over $\Spec \Z$
are  local and admit an open cover by representable subfunctors.
\end{thm}
\proof \cf\cite{demgab} Chapter I,\S~1, 4.4\endproof\vspace{.05in}

 \subsection{Monoids: the category $\Mo$.}\hfill \label{monoids}\vspace{.05in}

We recall that $\Mo$ denotes the category of commutative `pointed' monoids $M\cup\{0\}$, \ie $M$ is a semigroup with a commutative multiplicative operation `$\cdot$' (for simplicity we shall use the notation $xy$ to denote the product $x\cdot y$ in $M\cup\{0\}$) and an identity element $1$. Moreover, $0$ is an absorbing element in $M$ \ie $0\cdot x =x\cdot 0= 0,~\forall x\in M$.

The morphisms in $\Mo$ are unital homomorphisms of monoids $\varphi: M\to N$ ($\varphi(1_M) = 1_N$) satisfying the condition $\varphi(0) = 0$.

We also recall (\cf \cite{Gilmer} p. 3) that an ideal $I$ of a monoid $M$ is a non-empty subset $I\subset M$ such that  $I\supseteq xI=\{xi|i\in I\}$ for each $x\in M$. An ideal $I\subset M$ is {\em prime} if it is a proper ideal $I\subsetneq M$ and $xy\in I~\implies~x\in I~\text{or}~y\in I~\forall x,y\in M$. Equivalently, a proper ideal $\ffp\subsetneq M$ is prime if and only if $\ffp^c:=M\setminus\ffp$ is a multiplicative subset of $M$, \ie $x\notin \ffp,~y\notin \ffp \implies xy \notin \ffp$.

It is a standard fact that the pre-image of a prime ideal by a morphism of monoids is a prime ideal. Moreover, it is also straightforward to verify that the complement $\ffp_M=(M^\times)^c$ of the set of invertible elements in a monoid $M$ is a prime ideal in $M$ which contains all other prime ideals of the monoid. This interesting fact points out to one of the main properties that characterize monoids, namely monoids are local algebraic structures.

 We recall that the smallest ideal containing a collection of ideals $\{I_\alpha\}$ of a monoid $M$ is the union $I = \cup_\alpha I_\alpha$.

If $I\subset M$ is an ideal of a monoid $M$, the relation $\sim$ on $M$ defined by
\[
x\sim y~\Leftrightarrow~x=y~~\text{or}~~x,y\in I
\]
is an example of a {\em congruence} on $M$, \ie an equivalence relation on $M$ which is compatible with the product, (\cf \cite{Grillet}, \S 1 Proposition 4.6). The quotient monoid $M/I := M/\sim$ (Rees quotient) is identifiable with the pointed monoid $(M\setminus I)\cup\{0\}$, with the operation defined as follows
\[
x\star y = \begin{cases}xy& \text{if $xy\notin I$,}
\\0& \text{if $xy\in I$.}
\end{cases}
\]
Another interesting example of congruence in a monoid is provided by the operation of {\em localization} at a multiplicative subset $S\subset M$. One considers the congruence $\sim$ on the submonoid $M\times S\subset M\times M$ generated by the relation
\[
(m,s)\sim (m',s')\quad \Leftrightarrow\quad\exists~u\in S\quad ums' = um's.
\]
By introducing the symbol $m/s=(m,s)$, one can easily check  that the product
$
(m/s).(m'/s')=mm'/ss'
$
is well-defined on the quotient monoid $S^{-1}M=(M\times S)/\sim$. A particular case of this construction is when, for $f\in M$, one considers the multiplicative set $S=\{f^n; n\in \Z_{\ge 0}\}$: in analogy with rings, the corresponding quotient monoid $S^{-1}M$ is usually denoted by $M_f$.\vspace{.05in}

For an ideal $I\subset M$, the set $D(I)=\{\ffp\subset M|\ffp~\text{prime ideal},~\ffp\nsupseteq I\}$ determines an open set for the natural topology on the set $\Spec(M)=\{\ffp\subset M|\ffp~\text{prime ideal}\}$ (\cf~\cite{Kato}). For $I = \cup_\alpha I_\alpha$ ($\{I_\alpha\}$ a collection of ideals), the corresponding open set $D(I)$  satisfies the property  $D(\cup_\alpha I_\alpha) = \cup_\alpha D(I_\alpha)$.

The following equivalent statements characterize the open subsets $D(I)\subset\Spec(M)$.
\begin{prop}\label{prop} Let $\rho: M \to N$ be a morphism in the category $\Mo$ and let $I\subset M$ be an ideal. Then, the following conditions are equivalent:\vspace{.05in}

$(1)$~$\rho(I)N=N$.\vspace{.05in}

$(2)$~$1\in\rho(I)N$.\vspace{.05in}

$(3)$~$\rho^{-1}((N^\times)^c)$ is a prime ideal belonging to  $D(I)$.\vspace{.05in}

$(4)$~$\rho^{-1}(\ffp)\in D(I)$, for any prime ideal $\ffp\subset N$.
\end{prop}
\proof One has $(1)~\Leftrightarrow~(2)$. Moreover, $\rho^{-1}((N^\times)^c)\not\supset I$ if and only if $\rho(I)\cap N^\times \neq \emptyset$ which is equivalent to $(1)$. Thus $(2)~\Leftrightarrow~(3)$.

If an ideal $J\subset M$ does not contain $I$, then the same holds obviously for all the sub-ideals of $J$. Then $(3)~\Rightarrow~(4)$ since $\ffp_N=(N^\times)^c$ contains all the prime ideals of $N$. Taking $\ffp = \ffp_N$ one gets $(4)~\Rightarrow~(3)$.
\endproof

The proof of the following lemma is straightforward (\cf \cite{Gilmer}).
\begin{lem}\label{radical} Given an ideal $I\subset M$, the intersection  of the prime ideals $\ffp\subset M$, such that $\ffp\supset I$ coincides with the radical of $I$
\[
\bigcap_{\ffp\supset I}\ffp = \sqrt I := \{x\in M|\exists n\in\N, x^n\in I\}.
\]
\end{lem}

Given a commutative group $H$, the following definition determines a pointed monoid in $\Mo$
\[
\F_1[H] := H\cup\{0\}\qquad (0\cdot h = h\cdot 0 = 0,\quad\forall h\in H).
\]
Thus, in $\Mo$ a monoid of the form $\F_1[H]$ corresponds to a field $F$ ($F = F^\times\cup\{0\}$) in the category of commutative rings with unit. It is elementary to verify that the collection of monoids like $\F_1[H]$, for $H$ an abelian group, forms a full subcategory of $\Mo$ isomorphic to the category $\Ab$ of abelian groups.

In view of the fact that monoids are local algebraic structures, one can also introduce a notion which corresponds to that of the residue field for local rings and related homomorphism. For a commutative monoid $M$, one defines the pair $(\F_1[M^\times],\epsilon)$ where $\epsilon$ is the natural homomorphism
\begin{equation}\label{epsilon}
 \epsilon: M\stackrel{}{\to}\F_1[M^\times]\,\qquad \epsilon(y)=\begin{cases}0& \text{if $y\notin M^\times$}
\\y& \text{if $y\in M^\times$.}
\end{cases}
\end{equation}
The non-invertible elements of $M$ form a prime ideal, thus $\epsilon$ is a multiplicative map. The following lemma describes an application of this idea

\begin{lem} \label{eval} Let $M$ be a commutative monoid and $\ffp\subset M$ a prime ideal. Then\vspace{.05in}

a) $\ffp M_\ffp\cap (M_\ffp)^\times = \emptyset$.\vspace{.05in}

b) There exists a unique homomorphism $\epsilon_\ffp\, :M_\ffp\to  \F_1[(M_\ffp)^\times]$ such that\vspace{.05in}

$\bullet$~$\epsilon_\ffp(y)=0,\, \forall y\in \ffp$\vspace{.05in}

$\bullet$~ $\epsilon_\ffp(y)=y,\, \forall y\in (M_\ffp)^\times$.\vspace{.05in}

\end{lem}
\proof a) If $j\in \ffp$ then the image of $j$ in $M_\ffp$ cannot be invertible, since this would imply an equality of the form $sja=sf$ for some $s\notin \ffp$, $f\notin \ffp$ and hence a contradiction.

b) The first statement a) shows that
the corresponding (multiplicative) map $\epsilon_\ffp$ fulfills the two conditions.
To check uniqueness, note that the two conditions suffice to determine $\epsilon_\ffp(y)$ for any $y=a/f$ with $a\in M$ and $f\notin \ffp$.
\endproof\vspace{.05in}

\subsection{Geometric monoidal spaces.}\hfill \vspace{.05in}\label{geomesp}

In this subsection we review the construction of the geometric spaces which generalize, in the categorical context of the commutative monoids, the classical theory of the geometric $\Z$-schemes that we have reviewed in \S\ref{locfunctor}. We refer to  \S 9 in \cite{Kato}, \cite{deit} and \cite{deit1} for further details.\vspace{.05in}

A geometric monoidal space $(X,\cO_X)$ is a topological space $X$ endowed with a sheaf of monoids $\cO_X$ (the structural sheaf).
Unlike the case of the geometric $\Z$-schemes (\cf\cite{demgab}, Chapter I, \S~1 Definition 1.1), there is no need to impose the condition that the stalks of the structural sheaf of a geometric monoidal space are local algebraic structures,  since {\em by construction} any monoid is endowed with this property.

A morphism $\varphi: X \to Y$ of geometric monoidal spaces is a pair $(\varphi,\varphi^\sharp)$ of a continuous map $\varphi: X\to Y$ of topological spaces and a homomorphism of sheaves of monoids $\varphi^\sharp: \varphi^{-1} \cO_Y\to\cO_X$
which satisfies the property of being {\em local}, \ie  $\forall x\in X$ the homomorphisms connecting the stalks $\varphi^\sharp_x: \cO_{Y,\varphi(x)}\stackrel{}{\to}\cO_{X,x}$ are local, \ie they fulfill the following definition (\cf \cite{deit})

\begin{defn} \label{localmorphism} A homomorphism $\rho: M_1\to M_2$ of monoids is local if the following property holds
\begin{equation}\label{locality}
\rho^{-1}(M_2^\times) = M_1^\times.
\end{equation}
\end{defn}
This locality condition can be equivalently phrased by either one of the following statements\vspace{.05in}

$(1)$~ $\rho^{-1}((M_2^\times)^c) = (M_1^\times)^c$\vspace{.05in}

$(2)$~ $\rho((M_1^\times)^c)\subset (M_2^\times)^c$.\vspace{.05in}

The equivalence of \eqref{locality} and $(1)$ is clear since $\rho^{-1}(S^c) = \rho^{-1}(S)^c$, for any subset $S\subset M_2$. The equivalence of \eqref{locality} and $(2)$ follows since \eqref{locality} implies $(2)$. Conversely, if $(2)$ holds one has $\rho^{-1}((M_2^\times)^c) \supset (M_1^\times)^c$ and this latter jointly with $\rho^{-1}(M_2^\times)\supset M_1^\times$ imply $(1)$.\vspace{.05in}

We shall denote by $\GS$ the category of the geometric monoidal spaces.\vspace{.05in}

Notice that by construction the morphism $\epsilon_\ffp$ of Lemma \ref{eval} $b)$ is {\em local}. Thus it is natural to consider, for any geometric monoidal space $(X,\cO_X)$,  the analogue of the residue field at a point
$x\in X$ to be $\kappa(x)= \F_1[\cO_{X,x}^\times]$. Then, the associated evaluation map
\begin{equation}\label{evalmap}
\epsilon_x\, :\; \cO_{X,x}\to \kappa(x)= \F_1[\cO_{X,x}^\times]
\end{equation}
satisfies the properties as in $b)$ of Lemma~\ref{eval}.\vspace{.05in}

For a pointed monoid $M$ in $\Mo$, the set $\Spec(M)$  of the prime ideals $\ffp\subset M$ is called the {\em prime spectrum} of $M$ and is endowed with the topology whose closed subsets are the
\[
V(I) = \{\ffp\in\Spec(M)|\ffp\supset I\}
\]
as $I$ varies among the collection of all ideals of $M$.
Likewise for rings, the subset $V(I)\subset \Spec(M)$ depends only upon the radical of $I$ (\cf \cite{Bourbaki} II, Chpt. 2, \S 3). Equivalently, one can characterize the topology on $\Spec(M)$ in terms of a basis of open sets of the form $D(fM) = \{\ffp\in\Spec(M)|f\notin\ffp\}$, as $f$ varies in $M$.

The sheaf $\mathcal O$ of monoids associated to  $\Sp(M)$ is determined by the following properties:\vspace{.05in}

$\bullet$~The stalk of $\mathcal O$ at $\ffp\in\Spec(M)$ is $\mathcal O_{\ffp}=S^{-1}M$,  with  $S = \ffp^c$.\vspace{.05in}

$\bullet$~Let $U$ be an open set of $\Spec(M)$, then a section $s\in \Gamma(U,\mathcal O)$ is an element of $\prod_{\ffp\in U}\mathcal O_\ffp$ such that its restriction $s_{|D(fM)}$ to any open $D(fM)\subset U$ agrees with an element of $M_f$.\vspace{.05in}

$\bullet$~The homomorphism of monoids $\varphi: M_f \to \Gamma(D(fM),\mathcal O)$
  \[
\varphi(x)(\ffp) = \frac{a}{f^n}\in \mathcal O_\ffp,\quad \forall\ffp\in D(fM)\qqq x = \frac{a}{f^n}\in M_f
\]
is an isomorphism.\vspace{.05in}

For any monoidal geometric space $(X,\cO_X)$ one has a canonical morphism  $\psi_X:X\to \Sp(\Gamma(X,\cO_X))$, $\psi_X(x) = \ffp_x$, such that $s(x)=0$ in $\cO_{X,x},~\forall s\in\ffp_x$. It is easy to verify that a monoidal geometric space $(X,\cO_X)$ is a {\em prime spectrum} if and only if the morphism $\psi_X$ is an isomorphism.

\begin{defn} \label{geometricsch} A monoidal geometric space $(X,\cO_X)$ that admits an open cover by prime spectra $\Sp(M)$ is called a {\em  geometric $\Mo$-scheme}.
\end{defn}

Prime spectra fulfill the following locality property (\cf~\cite{deit}) that will be considered again in  \S \ref{mofunctors} in the functorial definition of an $\Mo$-scheme.

\begin{lem}\label{local2}  Let $M$ be an object in $\Mo$ and let $\{W_\alpha\}_{\alpha\in I}$  be an open cover of  the topological space $X=\Spec(M)$. Then $W_\alpha = \Spec(M)$, for some index $\alpha\in I$.
\end{lem}
\proof The point $\ffp_M = (M^\times)^c\in\Spec(M)$ must belong to some open set $W_\alpha$, hence $\ffp_M \in D(I_\alpha)$ for some index $\alpha\in I$ and this is equivalent to $I_\alpha\cap M^\times \neq\emptyset$, \ie $I_\alpha = M$.
\endproof\vspace{.05in}

\subsection{$\Mo$-schemes}\hfill \vspace{.05in}\label{mofunctors}

In analogy to and generalizing the theory of $\Z$-schemes, one develops the theory of $\Mo$-schemes following both the functorial and the geometrical viewpoint.\vspace{.05in}

\subsubsection{$\Mo$-functors}\hfill \vspace{.05in}\label{mofunctors1}

\begin{defn} \label{Mfunc} An $\Mo$-functor is a covariant functor $\cF: \Mo \to \Se$ from commutative (pointed) monoids to sets.
 \end{defn}

 To a (pointed) monoid $M$  in $\Mo$ one associates the $\Mo$-functor
\begin{equation}\label{specm}
\ssp M: \Mo \to \Se\qquad \ssp M(N)=\Hom_{\Mo}(M,N).
\end{equation}
 Notice that by applying Yoneda's lemma, a morphism of $\Mo$-functors (natural transformation)  such as $\varphi:\ssp M\to \cF$ is completely determined by the element $\varphi(id_M)\in \cF(M)$ and  moreover any such element gives rise to a morphism $\ssp M\to \cF$. By applying this fact to the functor $\cF = \ssp N$, for $N\in\text{obj}(\Mo)$, one obtains an {\em inclusion} of $\Mo$ as a full subcategory of the category of $\Mo$-functors (where morphisms are natural transformations).\vspace{.05in}

 An ideal $I\subset M$ defines the sub-$\Mo$-functor $\underline D(I)\subset\ssp M$
 \begin{equation}\label{di}
 \underline D(I): \Mo \to \Se,\qquad \underline D(I)(N) = \{\rho\in\ssp(M)(N)|\rho(I)N=N\}.
 \end{equation}\vspace{.05in}

 {\em Automatic locality.}\label{autloc}\hspace{.05in}

 We recall (\cf\S\ref{localZfunctors}) that for $\Z$-functors the ``locality'' property is defined by requiring, on coverings of prime spectra $\Spec(R)$ by open sets $D(f_i)$ of a basis, the exactness of sequences such as \eqref{locfunc0} and \eqref{locfunc}.
The next lemma shows  that locality is automatically verified for any $\Mo$-functor.

First of all, notice that for any $\Mo$-functor $\cF$ and any monoid $M$  one has a sequence of maps of sets
\begin{equation}\label{locf}
    \cF(M)\stackrel{u}{\longrightarrow}\prod_{i\in I}\cF(M_{f_i})\;\xy
{\ar@{->}_{w} (0,-1)*{}; (6,-1)*{}};
{\ar@{->}^{v} (0,1)*{}; (6,1)*{}};
\endxy\prod_{i,j\in I}
    \cF(M_{f_if_j})
\end{equation}
that is obtained by using the open covering of $\Spec(M)$ by the open sets $D(f_iM)$ of a basis ($I=$ finite index set) and the natural homomorphisms $M \to M_{f_i}$, $M_{f_i}\to M_{f_if_j}$.
\begin{lem} \label{lemauto} For any $\Mo$-functor $\cF$ the sequence \eqref{locf} is exact.
\end{lem}
\proof By Lemma~\ref{local}, there exists an index $i\in I$ such that $f_i\in M^\times$. We may also assume that $i = 1$. Then, the homomorphism $M\to M_{f_1}$ is invertible, thus $u$ is injective. Let $(x_i)\in \prod_i\cF(M_{f_i})$ be a family, with $x_i\in \cF(M_{f_i})$ and such that $v(x_i)=(x_i)_{f_j} = (x_j)_{f_i}=w(x_i)$, for all $i,j\in I$. This gives in particular the equality between the image of $x_i\in\cF(M_{f_i})$ under the isomorphism  $\cF(\rho_{i1}): \cF(M_{f_i}) \to \cF(M_{f_if_1})$ and $\cF(\rho_{i1})(x_1) \in \cF(M_{f_if_1})$. By writing $x_1 = \rho_1(x)$ one finds that $u(x)$ is equal to the family $(x_i)$.
\endproof\vspace{.05in}

 \subsubsection{Open $\Mo$-subfunctors.}\hfill \vspace{.05in}\label{openMo}

 In exact analogy with the theory of $\Z$-schemes (\cf \S\ref{opensubsub}), we introduce the notion of open subfunctors of $\Mo$-functors.

 \begin{defn} \label{opensubmo} We say that a subfunctor $\mathcal G\subset\cF$ of an $\Mo$-functor $\cF$ is open if for any object $M$ of $\Mo$ and any morphism of $\Mo$-functors $\varphi: \ssp M \to  \cF$ there exists an ideal $I\subset M$ satisfying the following property:\vspace{.05in}

For any object $N$  of $\Mo$ and for any $\rho\in\ssp M(N)=\Hom_{\Mo}(M,N)$ one has
\begin{equation}\label{condition1}
\varphi(\rho)\in \mathcal G(N)\subset\cF(N)~\Leftrightarrow~\rho(I)N = N.
\end{equation}
\end{defn}\vspace{.05in}

To clarify the meaning of the above definition we develop a few  examples.\vspace{.05in}

\begin{example} \label{mulsub} {\rm The functor
\[
 \mathcal G: \Mo\to \Se,\quad \mathcal G(N)=N^\times
 \]
 is an open subfunctor of the (identity) functor $\cD^1$
 \[
  \cD^1: \Mo \to \Se,\quad \cD^1(N)=N.
  \]
  In fact, let $M$ be a monoid, then by Yoneda's lemma a morphism of functors $\varphi: \ssp M \to  \cD^1$  is determined by an element $z\in \cD^1(M) = M$. For any monoid $N$ and $\rho\in\Hom(M,N)$, one has $\varphi(\rho)=\rho(z)\in  \cD^1(N)=N$, thus the condition $\varphi(\rho)\in  \mathcal G(N) = N^\times$ means that $\rho(z)\in N^\times$. One takes  for $I$ the ideal generated by $z$ in $M$: $I = zM$. Then, it is straightforward to check that \eqref{condition1} is fulfilled.}\end{example}

  \begin{example} \label{opensubf} {\rm We start with a monoid $M$ and an ideal $I\subset M$ and define the following sub-$\Mo$-functor of $\ssp(M)$
\[
\underline D(I): \Mo \to \Se,\quad \underline D(I)(N) = \{\rho\in\text{Hom}(M,N)|\rho(I)N = N\}.
\]
This means that for all prime ideals $\ffp\subset N$, one has $\rho^{-1}(\ffp)\not\supset I$ (\cf~Proposition~\ref{prop}). Next, we show that $\underline D(I)$ defines an open subfunctor of $\ssp M$. Indeed, for any object $A$ of $\Mo$ and any natural transformation $\varphi: \ssp A \to \ssp M$ one has $\varphi(id_A) = \eta\in\ssp M(A) = \text{Hom}_{\Mo}(M,A)$; we can take in $A$ the ideal $J = \eta(I)A$. This ideal fulfills the condition \eqref{condition1} for any object $N$ of $\Mo$ and any $\rho\in\text{Hom}_{\Mo}(A,N)$.
In fact, one has $\varphi(\rho) = \rho\circ\eta\in\text{Hom}_{\Mo}(M,N)$ and $\varphi(\rho)\in \underline D(I)(N)$ means that $\rho(\eta(I))N = N$. This holds if and only if $\rho(J)N = N$.
  }\end{example}\vspace{.05in}

 \subsubsection{Open covering by $\Mo$-subfunctors.}\hfill \vspace{.05in}\label{coverMo}

There is a natural generalization of the notion of open cover for $\Mo$-functors.
We recall that the category $\Ab$ of abelian groups embeds as a full subcategory of $\Mo$ by means of the functor $H\to \F_1[H]$.

\begin{defn} \label{defncover} Let $\cF$ be an $\Mo$-functor and let $\{{\cF_\alpha}\}_{\alpha\in I}$ ($I=$ an index set) be a family of open subfunctors of $\cF$. One says that $\{{\cF_\alpha}\}_{\alpha\in I}$ is an open cover of $\cF$ if
\begin{equation}\label{cover}
 \cF(\F_1[H]) = \bigcup_{\alpha\in I} {\cF_\alpha}(\F_1[H]),\quad\forall H\in{\rm Obj}(\Ab).
\end{equation}
\end{defn}
Since commutative groups (with a zero-element) replace fields in $\Mo$, the above definition is the natural generalization of the definition of open covers for $\Z$-functors (\cf\cite{demgab}) within the category of $\Mo$-functors. The following proposition gives a precise characterization of the properties of the open covers of an $\Mo$-functor.

\begin{prop} \label{coverspec} Let $\cF$ be an $\Mo$-functor and let $\{{\cF_\alpha}\}_{\alpha\in I}$ be a family of open subfunctors of $\cF$. Then, the family $\{{\cF_\alpha}\}_{\alpha\in I}$ forms an open cover of $\cF$ if and only if  \[
\cF(M)=\bigcup_{\alpha\in I} \cF_\alpha (M),\quad\forall~M\in{\rm obj}(\Mo).
\]
 \end{prop}

 \proof The condition is obviously sufficient. To show the converse, we assume \eqref{cover}. Let $M$ be a pointed monoid and let $\xi\in \cF(M)$, one needs to show that $\xi\in \cF_\alpha (M)$ for some $\alpha\in I$.
 Let $\varphi: \ssp M \to \cF$ be the morphism of functors with $\varphi(id_M)=\xi$. Since each $\cF_\alpha$ is an open subfunctor of $\cF$, one can find ideals $I_\alpha\subset M$ such that for any object $N$  of $\Mo$ and for any $\rho\in\ssp M(N)=\Hom_{\Mo}(M,N)$ one has
 \begin{equation}\label{condcover}
 \varphi(\rho)\in  \cF_\alpha(N)\subset  \cF(N)~\Leftrightarrow~\rho(I_\alpha)N = N.
\end{equation}
One applies this to the morphism $\epsilon_M:M\to \F_1[M^\times]=\kappa$ as in \eqref{epsilon}.
One has $\epsilon_M\in \ssp(M)(\kappa)$ and $\varphi(\epsilon_M)\in \cF(\kappa)=\bigcup_{\alpha\in I} {\cF_\alpha}(\kappa)$. Thus,
 $\exists \alpha\in I$ such that  $\varphi(\epsilon_M)\in \cF_\alpha(\kappa)$. By \eqref{condcover} one concludes that $\epsilon_M(I_\alpha)\kappa=\kappa$ and
$I_\alpha\cap M^\times\neq \emptyset$ hence $I_\alpha=M$. Applying then \eqref{condcover} to $\rho=id_M$ one obtains $\xi\in \cF_\alpha (M)$ as required. \endproof\vspace{.05in}

For $\cF=\ssp M$ and the $\cF_\alpha=\underline D(I_\alpha)$ open subfunctors corresponding to ideals $I_\alpha\subset M$, the covering condition in Definition~\ref{defncover} is equivalent to state that
$\exists \alpha\in I$ such that $I_\alpha=M$. In fact, one takes $H=M^\times$ and $\epsilon: M\to \F_1[M^\times]=\kappa$. Then,
$\epsilon\in  \ssp M(\F_1[H])$ and $\exists \alpha$, $\epsilon\in\underline D(I_\alpha)(\F_1[H])$ and thus
$I_\alpha\cap M^\times\neq \emptyset$ hence $I_\alpha=M$.\vspace{.05in}

Let $X$ be a geometric monoidal space and let $U_\alpha\subset X$ be (a family of) open subsets ($\alpha\in I$). One introduces the following $\Mo$-functors
\begin{equation}\label{functdef}
\underline X_\alpha: \Mo\to \Se,\qquad \underline X_\alpha(M):=\Hom(\Sp(M),U_\alpha).
\end{equation}
Notice that $\underline X_\alpha(M)\subset \underline X(M):= \Hom(\Sp(M),X)$.

\begin{prop} \label{coverspec1} The following conditions are equivalent\vspace{.05in}

$(1)$~$\underline X = \bigcup_{\alpha\in I}\underline X_\alpha$.\vspace{.05in}

$(2)$~$X = \bigcup_{\alpha\in I}U_\alpha$.\vspace{.05in}

$(3)$~$\underline X(M)=\bigcup_{\alpha\in I} \underline X_\alpha(M)\quad\forall M\in\text{obj}(\Mo)$.
  \end{prop}
\proof (1) $\implies$ (2). Assume that (2) fails and let $x\notin \bigcup_\alpha U_\alpha$. Then the local evaluation
map $\epsilon_x\, :\; \cO_{X,x}\to \kappa(x)$ of \eqref{evalmap} determines a morphism of geometric monoidal spaces $\Sp(\kappa(x))\to X$ and for $M=\kappa(x)$ a corresponding element $\epsilon\in \underline X(M)=\Hom(\Sp(M),X)$.
By applying Definition \ref{defncover}, there exists an index $\alpha$ such that $\epsilon\in \underline X_\alpha(M) $ and this shows that $x\in U_\alpha$ so that the open sets $U_\alpha$'s cover $X$.

(2) $\implies$ (3). Let $\phi\in \underline X(M)=\Hom(\Sp(M),X)$. Then, with $\ffp_M$ the maximal ideal of $M$, one has $\phi(\ffp_M)\in X=\bigcup_\alpha U_\alpha$ hence there exists an index $\alpha$ such that
$\phi(\ffp_M)\in  U_\alpha$. It follows that $\phi^{-1}(U_\alpha)\ni \ffp_M$ is $\Sp(M)$, and one gets $\phi\in \underline X_\alpha(M)=\Hom(\Sp(M),U_\alpha)$.

The implication (3) $\implies$ (1) is straightforward.
\endproof\vspace{.05in}

 \subsubsection{$\Mo$-schemes.}\hfill \vspace{.05in}\label{moschemes}

 In view of the fact that an $\Mo$-functor is local, the definition of an
 $\Mo$-scheme simply involves the local representability.

 \begin{defn}\label{defnmosch} An $\Mo$-scheme is an $\Mo$-functor which admits an open cover by representable subfunctors.
\end{defn}

We shall consider several elementary examples of $\Mo$-schemes\vspace{.05in}

\begin{example}\label{affineex}{\rm  The affine space $\cD^n$. For a fixed $n\in\N$, we consider the following $\Mo$-functor
\[
\cD^n: \Mo \to \Se,\qquad \cD^n(M) = M^n
\]
This functor  is representable since it is described by
\[
\cD^n(M) =  \text{Hom}_{\Mo}(\F_1[T_1,\ldots ,T_n],M),
\]
where
 \begin{equation}\label{F_1T}
\F_1[T_1,\ldots,T_n] := \{0\}\cup\{T_1^{a_1}\cdots T_n^{a_n}|a_j\in\Z_{\ge 0}\}\,
\end{equation}
is the pointed monoid associated to the semi-group generated by the variables $T_j$.
}\end{example}

\begin{example}\label{projex}{\rm The projective line $\underline \P^1$.  We consider the $\Mo$-functor $\underline \P^1$ which associates to an object $M$ of $\Mo$  the set $\underline\P^1(M)$ of complemented submodules $E\subset M^2$ of rank one, where the rank is defined locally. By definition, a complemented submodule is the range of an idempotent matrix $e\in M_2(M)$ (\ie $e^2=e$) whose rows have {\em at most}\footnote{Note that we need the $0$-element to state this condition} one non-zero entry.   To a morphism $\rho:M\to N$ in $\Mo$, one associates the following map
$$
\underline\P^1(\rho)(E)= N\otimes_M E\subset N^2.
$$
In terms of projectors $\underline\P^1(\rho)(e)=\rho(e)\in M_2(N)$. The condition of rank one means that for any prime ideal $\ffp\in \Sp M$ one has $\epsilon_\ffp(e)\notin\{0,1\}$ where $\epsilon_\ffp:M\to\F_1[(M_\ffp)^\times]$ is the morphism of Lemma~\ref{eval} $b)$ (where  $M_\ffp=S^{-1}M$, with $S = \ffp^c$).

Now, we compare $\underline \P^1$ with the $\Mo$-functor
\begin{equation}\label{projlinemo}
\mathcal P(M) = M \cup_{M^\times}M
\end{equation}
where the gluing map is given by $x\to x^{-1}$. In other words, we define on the disjoint union $ M\cup M$ an equivalence relation given by (using the identification $M\times\{1,2\} = M\cup M$)
\[
(x,1)\sim (x^{-1},2)\quad\forall x\in M^\times.
\]
 We define a natural transformation $e:\mathcal P\to\underline \P^1$ by observing that the matrices
$$
e_1(a)=\left(
  \begin{array}{cc}
    1 & 0 \\
    a & 0 \\
  \end{array}
\right)\,, \ \
e_2(b)=\left(
  \begin{array}{cc}
    0 & b \\
    0 & 1 \\
  \end{array}
\right)\,, \ \ a, b\in M
$$
are idempotent ($e^2=e$) and their ranges also fulfill the following property
$$
 {\rm Im}\,e_1(a)={\rm Im}\,e_2(b)\iff ab=1.
$$

\begin{lem}\label{projtwocomp} The natural transformation $e$ is bijective on the objects \ie
$$
 \mathcal P(M)= M \cup_{M^\times}M \cong \underline\P^1(M).
$$
Moreover, $\underline \P^1$ is covered by two copies of representable sub-functors $\cD^1$.
\end{lem}
\proof We refer to \cite{announc3} Lemma 3.13.\endproof}\end{example}

\begin{example}\label{opensubff}{\rm  Let $M$ be a monoid and let $I\subset M$ be an ideal. Consider the $\Mo$-functor $\underline D(I)$ of Example \ref{opensubf}. The next proposition states that this functor is an $\Mo$-scheme.
}
\end{example}

\begin{prop}\label{opensubfff} 1) Let $f\in M$ and $I = fM$. Then the subfunctor $\underline D(I)\subset\ssp M$ is represented by $M_f$.\vspace{.05in}

2) For any ideal $I\subset M$, the $\Mo$-functor
$\underline D(I)$ is an $\Mo$-scheme.
\end{prop}

\proof $1)$ For any monoid $N$ and for $I = fM$, one has
\[
\underline D(I)(N) = \{\rho\in\text{Hom}_{\Mo}(M,N)|\rho(f)\in N^\times\}.
 \]
 The condition $\rho(f)\in N^\times$ means  that $\rho$ extends to a morphism $\tilde\rho\in\text{Hom}(M_f,N)$, by setting
\[
\tilde\rho(\frac{m}{f^n}) = \rho(m)\rho(f)^{-n}\in N.
\]
Thus one has a canonical and functorial isomorphism $\underline D(I)(N) \simeq \text{Hom}_{\Mo}(M_f,N)$ which proves the representability of $\underline D(I)$ by $M_f$.

$2)$ For any $f\in I$, the ideal $fM\subset I$ defines a subfunctor of $\underline D(fM)\subset\underline D(I)$. This sub-functor is open because it is already open in $\ssp M$ and because there are less morphisms of type $\Spec(A)\to D(I)$ than those of type  $\Spec(A)\to\Spec(M)$, as $\underline D(I)(A)\subset\ssp M(A)$. Moreover, by $1)$, this subfunctor is representable.
\endproof\vspace{.05in}

 \subsubsection{Geometric realization.}\hfill \vspace{.05in}\label{geomereal}

In this final subsection we describe the construction of the geometric realization of an $\Mo$-functor (and scheme) following, and generalizing for commutative monoids, the exposition presented in \cite{demgab} \S 1 n.~4 of the geometric realization of a $\Z$-scheme for rings. Although the general development of this construction presents clear analogies with the case of rings, new features also arise which are specifically inherent to the discussion with monoids. The most important one is that the full sub-category $\Ab$ of $\Mo$ which plays the role of ``fields" within monoids, admits a final object. This fact simplifies greatly the description of the geometric realization of an $\Mo$-functor as we show in Proposition \ref{ex1} and in Theorem \ref{explicitgereal}.\vspace{.05in}

An $\Mo$-functor $\cF$ can be reconstructed by an inductive limit $\varinjlim\ssp(M)$ of representable functors, and the geometric realization is then defined by the inductive limit $\varinjlim\Spec(M)$, \ie by trading the $\Mo$-functors $\ssp(M)$ for the geometric monoidal spaces $(\Spec(M),\cO)$.  However, some set-theoretic precautions are also needed since the inductive limits which are taken over pairs $(M,\rho)$, where $M\in\text{obj}(\Mo)$ and $\rho\in\cF(M)$, are indexed over families rather than sets. Another issue arising in the construction is that of seeking that natural transformations of $\Mo$-functors form a set.
We bypass this problem by adopting the same caution as in \cite{demgab} (\cf Conventions g\'en\'erales). Thus, throughout this subsection we shall work with {\em models} of commutative monoids in a fixed suitable {\em universe} ${\bf U}$.
This set up is accomplished by introducing the full subcategory $\Mo'\subset \Mo$ of models of commutative  monoids in $\Mo$, \ie of commutative (pointed) monoids whose underlying set is an element of the fixed universe {\bf U}. To lighten the notations and the statements we adopt the convention that all the geometric monoidal spaces considered here below only involve monoids in $\Mo'$, thus we shall not keep the distinction between the category $\Mo\mathfrak S$ of $\Mo$-functors and the category $\Mo'\mathfrak S$ of (covariant) functors from $\Mo'$ to $\Se$.
 \vspace{.05in}

Given an $\Mo$-functor $\cF: \Mo\to \Se$, we denote by $\mathcal M_\cF$ the category of {\em $\cF$-models}, that is the category whose objects are pairs $(M,\rho)$, where $M\in\text{obj}(\Mo')$ and $\rho\in\cF(M)$. A morphism $(M_1,\rho_1)\to (M_2,\rho_2)$ in $\mathcal M_\cF$ is given by assigning a morphism $\varphi: M_1\to M_2$ in $\Mo'$ such that $\varphi(\rho_1) = \rho_2\in\cF(M_2)$.\vspace{.05in}

For a geometric monoidal space $X$, \ie  an object of the category $\GS$, we  recall the definition of the {\em functor defined by $X$} (\cf\eqref{functdef}). This is the $\Mo$-functor
\[
\underline X: \Mo\to \Se\qquad \underline X(M)=\Hom_{\GS}(\Spec(M),X).
\]
This definition can of course be restricted to the subcategory $\Mo'\subset \Mo$ of models.

One also introduces the functor
\[
d_\cF:\mathcal M_\cF^{op} \to \GS,\qquad d_\cF(M,\rho)=\Spec(M)
\]
\begin{defn} The geometric realization of an $\Mo$-functor $\cF$ is  the geometric monoidal space
defined by the inductive limit
\[
 |\cF| = \displaystyle{\varinjlim_{(M,\rho)\in\mathcal M_\cF}} d_\cF(M,\rho).
\]
The functor $|\cdot|:\cF \to |\cF|$ is called the functor geometric realization.
\end{defn}
For each $(M,\rho)\in\text{obj}(\mathcal M_\cF)$, one has a canonical morphism
\begin{equation}\label{canmorphism}
   i(\rho): d_{\cF}(M,\rho)=\Spec(M)\to \varinjlim d_{\cF}=|\cF|.
\end{equation}

The following result shows that the functor geometric realization $|\cdot|$ is left adjoint to the functor $\underline{~\cdot~}$.

\begin{prop}\label{adjointgeom} For every geometric monoidal space $X$ there is a bijection of sets which is functorial in $X$
\[
\varphi(\cF,X): \Hom_{\GS}(|\cF|,X) \stackrel{\sim}{\to} \Hom_{\Mo\mathfrak S}(\cF,\underline X),\qquad \varphi(\cF,X)(f)=g
\]
where $\forall (M,\rho)\in\rm{obj}(\mathcal M_{\cF})$
\[
g(M): \cF(M)\to \underline X(M),\qquad g(M)(\rho)=f\circ i(\rho)
\]
 and $i(\rho): \Spec(M)\to |\cF|$ is the canonical morphism of \eqref{canmorphism}.\vspace{.05in}

In particular, the functor geometric realization $|\cdot|$ is left adjoint to the functor $\underline{~\cdot~}$.
\end{prop}
\proof The proof is similar to the proof given in \cite{demgab} (\cf \S 1 Proposition~4.1).
\endproof
One obtains, in particular, the following morphisms
\[
\varphi(\cF,|\cF|): \Hom_{\GS}(|\cF|,|\cF|) \stackrel{\sim}{\to}\Hom_{\Mo\mathfrak S}(\cF,\underline{|\cF|}),\quad \varphi(\cF,|\cF|)(id_{|\cF|})=:\Psi(\cF)
\]
\[
\varphi(\underline X,X)^{-1}: \Hom_{\Mo\mathfrak S}(\underline X,\underline X)\stackrel{\sim}{\to}\Hom_{\GS}(|\underline X|,X),\quad \varphi(\underline X,X)^{-1}(id_{\underline X}) =: \Phi(X).
\]

Notice that as a consequence of the fact that the functor $\Spec(\cdot)$ is fully-faithful, $\Psi(\ssp(M))$ is invertible, thus $\Spec(M) \simeq |\ssp(M)|$, \ie $\Phi(\Spec(M))$ is invertible, for any model $M$ in $\Mo'$, \cf also Example~\ref{ex3}.
\medskip

We have already said that the covariant functor
\[
\F_1[~\cdot~]: \Ab\to \Mo\qquad H\mapsto \F_1[H]
\]
embeds the category of abelian groups as a full subcategory of $\Mo$. We shall identify $\Ab$ to this full subcategory of $\Mo$. One has a pair of adjoint functors: $H\mapsto \F_1[H]$ from $\Ab$ to $\Mo$ and  $M\mapsto M^\times$ from $\Mo$ to $\Ab$ \ie one has the following natural isomorphism
$$\Hom_{\Mo}(\F_1[H],M) \cong \Hom_{\Ab}(H,M^\times).$$
Moreover, any $\Mo$-functor $\cF: \Mo\to \Se$ restricts to $\Ab$ and gives rise to a functor: its {\em Weil restriction}
\[
\cF_{|\Ab}: \Ab \to \Se
\]
 taking values into sets.

 Let $X$ be a geometric monoidal space. Then, the Weil restriction $\underline X_{|\Ab}$ of the $\Mo$-functor $\underline X$ to the full sub-category $\Ab\subset \Mo$,   is a direct sum of {\em indecomposable} functors, \ie  a direct sum of  functors that cannot be decomposed further into a disjoint sum of non-empty sub-functors.

\begin{prop} \label{propdec} Let $X$ be a geometric monoidal space. Then the functor
  $$
  \underline X_{|\Ab}: \Ab\to \Se,\qquad \underline X_{|\Ab}(H)=\Hom_{\mathfrak G\mathcal S}(\Sp\F_1[H],X)
  $$
  is the disjoint union of representable functors
  $$
  \underline X_{|\Ab}(\F_1[H]) =  \coprod_{x\in X}(\underline X_{|\Ab})_x(\F_1[H])=\coprod_{x\in X} \Hom_\Ab(\cO_{X,x}^\times,H)
  $$
  where $x$ runs through the points of $X$.
    \end{prop}

\proof Let $\varphi \in  \Hom(\Sp\F_1[H],X)$. The unique point $\ffp\in \Sp\F_1[H]$ corresponds to the ideal $\{0\}$. Let
$\varphi(\ffp)=x\in X$ be its image; there is a corresponding map of the stalks
$$
\varphi^{\#}\,:\cO_{X,\varphi(p)}\to \cO_\ffp=\F_1[H].
$$
This homomorphism is local by hypothesis: this means that the inverse image of  $\{0\}$ by $\varphi^{\#}$ is
the maximal ideal of $\cO_{X,\varphi(p)}=\cO_{X,x}$. Therefore, the map $\varphi^{\#}$ is entirely determined by the group homomorphism $\rho:\cO_{X,x}^\times\to H$ obtained as the restriction of $\varphi^{\#}$. Thus $\varphi \in  \Hom(\Sp\F_1[H],X)$ is entirely specified by a point $x\in X$ and a group homomorphism $\rho \in \Hom_{\Ab} (\cO_{X,x}^\times,H)$.
\endproof

 We recall that the set underlying the geometric realization of a $\Z$-scheme can be obtained by restricting the $\Z$-functor to the full subcategory of fields and passing to a suitable inductive limit (\cf\cite{demgab}, \S 4.5). For $\Mo$-schemes this construction simplifies considerably since the
full subcategory $\Ab$ admits the {\em final object}
\[
\F_1=\F_1[\{1\}]=\{0,1\}.
 \]
  Notice though that while $\F_1$ is a final object for the subcategory $\Ab\subset\Mo$, it is {\em not}  a final object in $\Mo$, in fact the following proposition shows that for any monoid $M\in\text{obj}(\Mo)$, the set $\Hom_\Mo(M,\F_1)$ is  in canonical correspondence with the points of $\Spec M$.\vspace{.05in}

For any geometric monoidal scheme $X$, we  denote by $X^{\mathfrak e}$ its underlying set.

\begin{prop}\label{ex1}
$(1)$~For any object $M$ of $\Mo$, the map
\[
\varphi:\Spec(M)^{\mathfrak e}\to \Hom_\Mo(M,\F_1)\qquad\ffp\mapsto\varphi_\ffp
\]
\begin{equation}\label{phip}
    \varphi_\ffp(x)=0\qqq x\in \ffp\, , \ \ \varphi_\ffp(x)=1\qqq x\notin \ffp
\end{equation}
determines a  natural bijection of sets.

$(2)$~For any $\Mo$-functor $\cF$, there is a canonical isomorphism $~|\cF|^{\mathfrak e}\simeq \cF(\F_1)$.

$(3)$~Let $X$ be a geometric monoidal space, and $\cF = \underline X$. Then, there is a canonical bijection
$\underline X(\F_1)\simeq X$.
\end{prop}
\proof  $(1)$~ The map $\varphi$
is well-defined since the complement of a prime ideal $\ffp$ in a monoid is a multiplicative set. To define the inverse of $\varphi$, one assigns to $\rho\in \Hom_\Mo(M,\F_1)$ its kernel which is a prime ideal of $M$ that uniquely determines $\rho$.

$(2)$~ The set  $|\cF|^{\mathfrak e}$, \ie the set underlying $\varinjlim d_\cF$ is canonically in bijection with
\[
Z=\varinjlim_{(M,\rho)\in\mathcal M_\cF}(\Spec(M))^\mathfrak e
\]
By $(1)$, one has a canonical bijection
\[
\Spec(M)^\mathfrak e\simeq \Hom_\Mo(M,\F_1)
\]
and using the identification $$
\cF=\varinjlim_{(M,\rho)\in\mathcal M_\cF}\ssp(M)
$$
we get the canonical bijection
$$
Z=
\varinjlim_{(M,\rho)\in\mathcal M_\cF}\Hom_\Mo(M,\F_1)
=
\varinjlim_{(M,\rho)\in\mathcal M_\cF}\ssp(M)(\F_1)=\cF(\F_1).
$$
$(3)$~ An element $\rho\in \cF(\F_1)=\underline X(\F_1)$ is a morphism of geometric spaces $\Spec\F_1\to X$ and the image of the unique (closed) point of $\Spec\F_1$ is a point $x\in X$ which uniquely determines $\rho$. Conversely any point $x\in X$ determines  a morphism $(\varphi,\varphi^\sharp)$ of geometric spaces $\Spec\F_1\to X$. The homomorphism of sheaves of monoids $\varphi^\sharp: \varphi^{-1} \cO_X\to\cO_{\F_1}$ is uniquely determined by its locality property, and sends $\cO_{X,x}^\times$ to $1$ and the complement of $\cO_{X,x}^\times$ (\ie the maximal ideal $\mathfrak m_x$) to $0$.
\endproof

Proposition \ref{ex1} shows how to describe concretely the set underlying the geometric realization $|\cF|$ of an $\Mo$-functor $\cF$ in terms of the set $\cF(\F_1)$. To describe the topology of $|\cF|$ directly on $\cF(\F_1)$, we  use the following construction of sub-functors of $\cF$ defined in terms of arbitrary subsets $P\subset \cF(\F_1)$. It corresponds to the construction of \cite{demgab} I, \S 1, 4.10.

\begin{prop}\label{equequ} Let $\cF$ be an $\Mo$-functor and $P\subset \cF(\F_1)$ be a subset of the set $\cF(\F_1)$. For $M\in\text{obj}(\Mo)$ and $\rho\in \cF(M)$ the following two conditions are equivalent:

a) For every homomorphism $\varphi\in\Hom_\Mo(M,\F_1)$ one has
$$\varphi(\rho)\in P\subset \cF(\F_1)$$

b) For every $H\in\text{obj}(\Ab)$ and every homomorphism $\varphi\in\Hom(M,\F_1[H])$ one has
 $$\varphi(\rho)\in \bigcup_{x\in P}({\cF}_{|\Ab})_x(\F_1[H])$$
where $P$ is viewed as a subset of $~\varinjlim {\cF}_{|\Ab} \simeq |\cF|$.\vspace{.05in}

The above equivalent conditions define a sub-functor $~\cF_P\subset\cF$.
\end{prop}

It is clear that $a)$ defines a sub-functor  $\cF_P\subset\cF$. We shall omit the proof of the equivalence between the conditions $a)$ and $b)$ which is only needed to carry on the analogy with the construction in \cite{demgab} I, \S 1, 4.10.

Thus, {\em we adopt $a)$ as the definition of the sub-functor  $\cF_P\subset\cF$}.

Notice that one can reconstruct $P$ from $\cF_P$ using the equality
\[
P = \cF_P(\F_1)\subset \cF(\F_1).
\]

Let $X$ be a geometric monoidal space. If one endows the subset $P\subset X$ with the induced topology from $X$ and with the structural sheaf ${\cO_X}_{|P}$, the functor $\underline X_P$ gets identified with $\underline P\subset\underline X$. Moreover, if $P$ is an open subset in $X$, then $\underline X_P$ is an open sub-functor of $\underline X$. In particular, if $\cF = \ssp(M)$ then $P\subset\Spec(M)=|\cF|$ is open if and only if $\cF_P$ is an open sub-functor of $\ssp(M)$.\vspace{.05in}

More generally, the following results hold for arbitrary  $\Mo$-functors.

\begin{thm}\label{explicitgereal} Let $\cF$ be an $\Mo$-functor. The following facts hold\vspace{.05in}

$(1)$~A subset $P\subset|\cF|= \cF(\F_1)$ is open if and only if for every $(M,\rho)\in\text{obj}(\mathcal M_\cF)$ the set of prime ideals $\ffp\in \Spec M$ such that $\varphi_\ffp(\rho)\in P$ is open in $\Spec M$.\vspace{.05in}

$(2)$~The map $P\to \cF_P$ induces a bijection of sets between the collection of open subsets of $|\cF|=\cF(\F_1)$ and the open sub-functors of $\cF$.\vspace{.05in}

$(3)$~The structure sheaf of $|\cF|=\cF(\F_1)$ is given by
$$
\cO(U)=\Hom_{\Mo\mathfrak S}(\cF_U, \cD)\qqq U\subset|\cF| \ \textit{open}\,,
$$
where $\cD:\Mo\to \Se$, $\cD(M)=M$ is the functor affine line.
\end{thm}

\proof $(1)$ By definition of the topology on the inductive limit $|\cF| = \displaystyle{\varinjlim_{(M,\rho)\in\mathcal M_\cF}} d_\cF$, a subset $P$ of $|\cF|$ is open if and only if its inverse image under the canonical maps $i(\rho)$ of \eqref{canmorphism} is open. Using the identification $|\cF|\simeq \cF(\F_1)$ one obtains the required statement.

$(2)$ We give the simple direct argument (\cf also \cite{demgab} I, \S 1 Proposition~4.12). Let $\cF_1$ be an open sub-functor of $\cF$. By Definition \ref{opensubmo}, given an object $M$ of $\Mo$ and an element $\xi \in \cF(M)$, there exists an ideal $I\subset M$ such that for  any object $N$  of $\Mo$ and for any $\rho\in\ssp M(N)=\Hom_{\Mo}(M,N)$ one has
$$
\cF(\rho)\xi\in \cF_1(N)\subset\cF(N)~\Leftrightarrow~\rho(I)N = N.
$$
Take $N=\F_1$ and $\rho=\varphi_\ffp$ for $\ffp\in \Spec M$, then one gets
$$
\varphi_\ffp(\xi)\in \cF_1(\F_1)~\Leftrightarrow~1\in \varphi_\ffp(I)~\Leftrightarrow~\ffp\in D(I)
$$
which shows that $D(I)$ and hence the radical of $I$ is determined by the subset  $P=\cF_1(\F_1)\subset \cF(\F_1)$. Now take $N=M$ and $\rho=id_M$, then one has
$$
\xi\in \cF_1(M)\subset\cF(M)~\Leftrightarrow~I = M\,,
$$
and this holds if and only if $D(I)=\Spec M$ or equivalently $\varphi_\ffp(\xi)\in \cF_1(\F_1)$ for all
$\ffp\in \Spec M$ which is the definition of the sub-functor associated to the open subset $P=\cF_1(\F_1)\subset \cF(\F_1)$.

$(3)$ By using the above identifications, the proof is the same as in \cite{demgab} I, \S 1 Proposition~4.14. The structure of monoid on $\cO(U)=\Hom_{\Mo\mathfrak S}(\cF_U, \cD)$ is given by
$$
(\alpha\beta)(M)(x)=\alpha(M)(x)\beta(M)(x)\in M \qqq x\in \cF_U(M),  \ \alpha,\beta \in \cO(U).
$$
\endproof
\begin{example}\label{ex3}{\rm Let $M\in\text{obj}(\Mo)$. The map $\ffp\mapsto \varphi_\ffp$ as in  \eqref{phip} defines a natural bijection $\Spec M\simeq|\ssp M|=\Hom(M,\F_1)$ and using the topology defined by (1) of Theorem \ref{explicitgereal}, one can check directly that this map is an homeomorphism.  One can also verify directly that the structure sheaf of $|\ssp M|$ as defined by (3) of Theorem \ref{explicitgereal} coincides with the structure sheaf of $\Spec M$. Thus $\Spec(M) \simeq |\ssp(M)|$, \ie $\Phi(\Spec(M))$ and $\Psi(\ssp(M))$ are invertible, for any object $M$ in $\Mo$.
}\end{example}

\begin{cor}\label{thecor} Let $X$ be a geometric monoidal space. If $\cF\subset\underline X$ is an open sub-functor of $\underline X$, then $|\cF|\subset |\underline X|$ is an open subset of $X$.
\end{cor}

\proof For the proof we refer to \cite{demgab} I, \S 1 Corollary~4.15.\endproof\vspace{.05in}

 The following result describes a property which {\em does not hold} in general for $\Z$-schemes and provides a natural map
$\cF(M)\to |\cF|\simeq \cF(\F_1)$, $\forall M\in\text{obj}(\Mo)$ and any $\Mo$-functor $\cF$.

\begin{prop} \label{coverspecbis}
 Let $\cF$ be an $\Mo$-functor, then the following facts hold

1) For any monoid $M\in\text{obj}(\Mo)$ there exists a canonical map
\begin{equation}\label{pim}
    \pi_M\,:\, \cF(M)\to |\cF|\simeq \cF(\F_1)
\end{equation}
such that
\begin{equation}\label{pim1}
\pi_M(\rho)=\varphi_{\ffp_M}(\rho)\qqq \rho\in \cF(M).
\end{equation}

2) Let $U$ be an open subset of $|\cF|= \cF(\F_1)$ and $\cF_U$ the associated subfunctor of $\cF$, then
 \begin{equation}\label{pim2}
    \cF_U(M)=\pi_M^{-1}(U)\subset \cF(M).
 \end{equation}
  \end{prop}

\proof $1)$ follows from the definition of the map $\varphi_{\ffp_M}$.

$2)$ By Theorem \ref{explicitgereal}, the set $W=\{\ffp\in \Spec M|\varphi_\ffp(\rho)\in U\}$ is open in $\Spec M$. Thus by Lemma \ref{local}, one has $W=\Spec(M)$ if and only if $\ffp_M\in W$ which gives the required statement using the definition of the sub-functor $\cF_U$ as in Proposition \ref{equequ}.
\endproof

 The next result establishes the equivalence between the category of $\Mo$-schemes and the category of  geometric $\Mo$-schemes, using sufficient conditions for the morphisms
 \[
 \Psi(\cF): \cF \to \underline{|\cF|},\qquad \Phi(X): |\underline X|\to X
  \]
  to be invertible. This result generalizes Example \ref{ex3} where this invertibility has been proven in the case $\cF=\ssp M$, $X=\Spec M$ and for objects $M$ of $\Mo$.

\begin{thm}  \label{compare}

$(1)$~Let $X$ be a geometric $\Mo$-scheme. Then $\underline X$ is an $\Mo$-scheme and the morphism $\Phi(X): |\underline X| \to X$ is invertible. \vspace{.05in}

$(2)$~Let $\cF$ be an $\Mo$-scheme then the geometric monoidal space $|\cF|$ is a geometric $\Mo$-scheme and $\Psi(\cF): \cF \to \underline{|\cF|}$ is invertible.

Thus, the two functors $|\cdot|$ and $\underline{~\cdot~}$ induce quasi-inverse equivalences between the category of $\Mo$-schemes and the category of  geometric $\Mo$-schemes.
    \end{thm}

\proof $(1)$ In view of Proposition \ref{ex1}, for any geometric scheme $X$ the map $\Phi(X): |\underline X| \to X$ induces a bijection of the underlying sets. If $\{U_i\}_{i\in {\bf U}}$ is a covering of $X$ by prime spectra, it follows that $\Phi(U_i): |\underline U_i| \to U_i$ is an isomorphism. Moreover, since $U_i$ is an open subset of $X$, it follows from Corollary~\ref{thecor}  that $|\underline U_i|$ is an open subspace of $|\underline X|$. Thus, the topologies on the spaces $|\underline X|$ and $X$ as well as the related structural sheaves get identified by means of $\Phi(X)$. Hence  $\Phi(X)$ is invertible.

$(2)$  Let $\cF$ be an $\Mo$-scheme, we shall show that $|\cF|$ is a geometric $\Mo$-scheme. We let $\{\mathfrak U_i\}_{i\in{\bf U}}$ be an open covering of $\cF$ by (functors) prime spectra such that $\{\mathfrak U_{ij\alpha}\}$ is an open covering of $\mathfrak U_i\cap \mathfrak U_j$ by (functors) prime spectra. Then, $\{|\mathfrak U_i|\}$ is an open covering of $|\cF|$ by prime spectra, thus $|\cF|$ verifies the condition of Definition \ref{geometricsch} and is a  geometric $\Mo$-scheme. In view of the fact that $\Psi(\mathfrak U_i): \mathfrak U_i \to \underline{|\mathfrak U_i|}$ is invertible, the composite of $\Psi(\mathfrak U_i)^{-1}$ and the inclusion $\mathfrak U_i\hookrightarrow\cF$ gives a collection of morphisms $\Psi_i^{-1}(\cF): \underline{|\mathfrak U_i|} \to \cF$.  Moreover one has $\Psi_i^{-1}(\cF)_{|\underline{|\mathfrak U_{ij\alpha}|}} = \Psi_j^{-1}(\cF)_{|\underline{|\mathfrak U_{ij\alpha}|}}$, for all $i,j,\alpha$.  The automatic locality of $\cF$ then shows that the $\Psi_i^{-1}$ assemble to a single morphism $\underline{|\cF|}\to\cF$ which is the inverse of $\Psi(\cF): \cF \to \underline{|\cF|}$. In fact Proposition \ref{coverspec} shows that this inverse is surjective.
\endproof

\section{$\F_1$-schemes and their zeta functions}\label{repfunfunctor}

A general formalism of category theory allows one to glue together two categories connected by a pair of adjoint functors (\cf\cite{announc3} \S~4). This construction is particularly interesting when applied to the functorial formalism that describes the most common examples of schemes (of finite type) over $\F_1$ that we have reviewed so far. In all these cases in fact, the two categories naturally involved \ie $\Mo$ and $\An$ are linked by a pair of natural adjoint functors. These are
\begin{equation}\label{beta1}
 \beta: \Mo \to \An\,, \qquad   \beta(M)=\Z[M]
\end{equation}
which associates to a monoid the convolution ring $\Z[M]$ and the forgetful functor
\begin{equation}\label{units}
 \beta^*: \An \to \Mo\qquad   \beta^*(R)= R.
\end{equation}
which associates to a ring its underlying structure of multiplicative monoid. The resulting glued category $\Mr=\An\cup_{\beta,\beta^*} \Mo$ defines the most natural categorical framework on which one introduces the notions of an $\F_1$-functor and of a scheme over $\F_1$.
\begin{defn}\label{defnfunfunc} An $\F_1$-functor is a covariant functor from the category
$\Mr=\An\cup_{\beta,\beta^*} \Mo$ to the category of sets.
\end{defn}
The conditions imposed in the original definition of a variety over $\F_1$ (\cf\cite{Soule}) are now applied to a covariant functor $\mathcal X: \Mr \to \Se$. Such a functor determines a scheme (of finite type) over $\F_1$ if it also fulfills the following  properties.

\begin{defn}\label{defnfonesch} An $\F_1$-scheme is an $\F_1$-functor $\mathcal X: \Mr\to \mathfrak{Sets}$, such that:\vspace{.05in}

$\bullet$~The restriction $X_\Z$ of $\mathcal X$ to $\An$ is a $\Z$-scheme.\vspace{.05in}

$\bullet$~The restriction $\underline X$ of $\mathcal X$ to $\Mo$ is an $\Mo$-scheme.\vspace{.05in}

$\bullet$~The natural transformation $e: \underline X\circ\beta^* \to X_\Z$ associated to a field is a bijection (of sets).
\end{defn}

Morphisms of $\F_1$-schemes are natural transformations of the corresponding functors. \vspace{.05in}

\subsection{Torsion free Noetherian $\F_1$-schemes}\hfill \vspace{.05in}\label{torsionfree}

The first theorem stated in this subsection shows that for any Noetherian, torsion-free $\F_1$-scheme the function counting the number of rational points (for the associated $\Z$-scheme) over a finite field $\F_q$  is {\em automatically} polynomial. The same theorem also provides a complete description of its zeta function over $\F_1$ (originally introduced upside-down in \cite{Soule}).

By following  \opcit (\S~6 Lemme~1) and implementing the correction concerning the inversion in the formula of the zeta function, the definition of the zeta-function  of an algebraic variety $X =(\underline X, X_\C, e_X)$ over $\F_1$ such that $\# X_\Z(\F_q)=N(q)$, with $N(x)$  a polynomial function and $\#\underline X(\F_{1^n}) = N(q)$ if $n = q-1$, is given as the limit
\begin{equation}\label{zetadefn}
\zeta_X(s):=\lim_{q\to 1}Z(X,q^{-s}) (q-1)^{N(1)}.
\end{equation}
Here, the   function $Z(X,q^{-s})$  (\ie the Hasse-Weil exponential series) is defined by
\begin{equation}\label{zetadefn1}
Z(X,T) := \exp\left(\sum_{r\ge 1}N(q^r)\frac{T^r}{r}\right).
\end{equation}
One obtains
\begin{equation}\label{polynomial}
N(x)=\sum_{k=0}^d a_kx^k ~\Longrightarrow~ \zeta_X(s)=\prod_{k=0}^d(s-k)^{-a_k}.
\end{equation}
For instance, in the case of the projective line  $\P^1$ one has $\zeta_{\P^1}(s)=\displaystyle{\frac{1}{s(s-1)}}$.\vspace{.05in}

We say that an $\F_1$-scheme is Noetherian if the associated $\Mo$ and $\Z$-schemes are Noetherian
(\cf \cite{announc3} Definition~4.12). An $\F_1$-scheme is said to be torsion free if the groups $\cO_{X,x}^\times$ of the invertible elements of the monoids $\cO_{X,x}$ ($X$ is the associated geometric $\Mo$-scheme) are torsion free $\forall x\in X$. The following result is related to Theorem 1 of \cite{deit2}, but applies also to non-toric varieties.

\begin{thm}\label{dthmfonesch}
Let $\mathcal X$ be a torsion free, Noetherian $\F_1$-scheme and let $X$ be the geometric realization of its restriction $\underline X$ to $\Mo$. Then\vspace{.05in}

$(1)$~There exists a polynomial $N(x+1)$ with positive integral coefficients such that
  $$
  \# \,X(\F_{1^n})=N(n+1)\qqq n\in \N.
  $$
$(2)$~For each finite field $\F_q$, the cardinality of the set of points of the $\Z$-scheme
  $X_\Z$ which are rational over $\F_q$ is equal to $N(q)$.\vspace{.05in}

$(3)$~The zeta function of $\mathcal X$  has the following description
  \begin{equation}\label{prodfund}
    \zeta_{\mathcal X}(s)=\prod_{x\in X} \frac{1}{\left(1-\frac 1 s\right)^{\otimes^{n(x)}}},
\end{equation}
where $\otimes$ denotes Kurokawa's tensor product and  $n(x)$ is the local dimension of $X$ at the point $x$.
\end{thm}

In \eqref{prodfund}, we use the convention that when $n(x)=0$ the expression
$$
\left( 1-\frac 1 s\right)^{\otimes^{n(x)}}=s.
$$
We refer to \cite{Ku} and \cite{Manin} for the details of the definition of Kurokawa's tensor products and zeta functions and to \cite{announc3} for the proof of the Theorem.\vspace{.05in}

\subsection{Extension of the counting functions in the torsion case}\hfill \vspace{.05in}\label{torsion}

The remaining part of this section is dedicated to the computation of the zeta function of an arbitrary Noetherian $\F_1$-scheme $\mathcal X$. For simplicity, we will always denote by $X$ the geometric realization (\ie the geometric monoidal scheme) of the restriction $\underline X$ of $\mathcal X$ to $\Mo$ and we will shortly refer to $X$ as to the $\F_1$-scheme. In the general case, we shall prove that the counting function
$$
  \# \,X(\F_{1^n})=N_X(n+1)\qqq n\in \N
  $$
of an $\F_1$-scheme is no longer a polynomial function of $n$ and that its description involves periodic functions. First, we will show that there is  a {\em canonical} extension of the counting function $N_X(n)$ to the complex plane, as an entire function $N_X(z)$, whose growth is well controlled. Then, we will explain how to compute the zeta function \eqref{zetadefn} using $N_X(q)$ for arbitrary real
values of $q\geq 1$.

The  simplest example of a Noetherian $\F_1$-scheme is $X=\Sp(\F_{1^m})$.  The number of rational points of $X$ over $\F_{1^n}$ is the cardinality of the set of group homomorphisms $\Hom(\Z/m\Z,\Z/n\Z)$, \ie
 \begin{equation}\label{genform0}
 \# \,X(\F_{1^n})=\#\,\Hom(\Z/m\Z,\Z/n\Z)=\gcd(n,m).
 \end{equation}
This is a {\it periodic} function of $n$. More generally the following statement holds.

\begin{lem} \label{counttors} Let $X$ be a Noetherian $\F_1$-scheme. The counting function is a finite sum of monomials of the form
\begin{equation}\label{genform}
 \# \,X(\F_{1^n})=\sum_{x\in X} n^{n(x)}\prod_{j} \gcd(n,m_j(x))
\end{equation}
where $n(x)$ denotes the local dimension at $x\in X$ and the natural numbers $m_j(x)$ are the orders of the finite cyclic groups which compose the residue ``field" $\kappa(x)=\F_1[\cO_{X,x}^\times]$ \ie
$$
\cO_{X,x}^\times=\Z^{n(x)}\prod_{j} \Z/m_j(x)\Z.
$$
\end{lem}
\proof It follows from the proof of the first part of Theorem~\ref{dthmfonesch} that
$$
\# \,X(\F_{1^n})=\sum_{x\in X}\Hom_{\Ab} (\cO_{X,x}^\times,\Z/n\Z)
$$
and the result follows from \eqref{genform0}.
\endproof\vspace{.05in}

We continue by explaining how to extend canonically a function such as \eqref{genform} to an entire function on the complex plane.\vspace{.05in}

\subsubsection{Extension of functions from $\Z$ to $\C$}\hfill \vspace{.05in}

Let $f(z)$ be an entire function of $z\in \C$ and let $T(R,f)$  be its {\em characteristic function}:
\begin{equation}\label{jensen1}
    T(R,f)=\frac{1}{2\pi}\int_0^{2\pi}\log^+|f(Re^{i\theta})|d\theta.
\end{equation}
This function may be interpreted as the average magnitude of $\log^+|f(z)|$. One sets
\[
\log^+ x = \begin{cases}\log x &\text{if $x\ge 1$}\\
0 & \text{if $0<x<1$}\end{cases}
\]
and one denotes by
\begin{equation}\label{jensen2}
    N(a, R)=\sum_{f(z)=a\atop |z|<R} \log|\frac R z|
\end{equation}
the sum over the zeros of $(f-a)(z)$ inside the disk of radius $R$.
The First Fundamental Theorem of Nevanlinna theory states as follows

\begin{thm}\label{main1}{\rm(First Fundamental Theorem)}~ If $a\in \C$ then
$$
  N(a, R)\leq T(R,f)-\log|f(0)-a|+\epsilon(a,R)
$$
where
$$
|\epsilon(a,R)|\leq \log^+|a|+\log 2\,.
$$
\end{thm}

From this theorem one deduces that if an entire function $f(z)$ vanishes on $\Z\subset \C$, then one has
$$
N(0,R)\geq \sum_{z\in \Z\backslash 0\atop |z|<R} \log|\frac R z|\sim 2R
$$
and hence that
\begin{equation}\label{mingrowth}
    \underline \lim \frac{T(R,f)}{R}\geq 2.
\end{equation}\vspace{.05in}

We shall now use this inequality to prove the following result.

\begin{prop}\label{propextfun} Let $N(n)$ be a real function on $\Z$ of the form
\begin{equation}\label{data}
    N(n)=\sum_{j\ge 0} T_j(n)\,n^j
\end{equation}
where the $T_j(n)$'s are periodic functions. Then, there exists a unique polynomial $Q(n)$ and a unique entire
function $f:\C\to \C$ such that\vspace{.05in}

$\bullet$~$\overline{f(\bar z)}=f(z)$.\vspace{.05in}

$\bullet$~$\overline\lim~\frac{T(R,f)}{R}< 2$.\vspace{.05in}

$\bullet$~$N(n)-(-1)^nQ(n)=f(n)$\quad $\forall~n\in \Z$.\vspace{.05in}

\end{prop}

\proof Choose $L\in \N$  such that all the $T_j(n)$ fulfill $T_j(m+L)=T_j(m)$, for all $m\in\Z$. For each $j\in\N$, one has a Fourier expansion of the form
$$
T_j(m)=\sum_{\alpha\in I_L} t_{j,\alpha}e^{2\pi i \alpha m} + s_j (-1)^m, \qquad  I_L=\{ \frac kL~|~ |k|<L/2 \}.
$$
We define $Q(n)=\sum_j s_j n^j$ and the function $f(z)$ by
$$
f(z)=\sum_j\sum_{\alpha\in I_L} t_{j,\alpha}e^{2\pi i \alpha z}z^j.
$$
Since $|\alpha|\leq c<\frac 12$, $\forall~\alpha\in I_L$, one gets that  $$\log^+|f(Re^{i\theta})|\leq 2\pi R c|\sin\theta|+o(R)$$
which in turn implies $\overline \lim~\frac{T(R,f)}{R}< 2$. The equality $N(n)-(-1)^nQ(n)=f(n)$ for all $n\in \Z$ holds by construction and it remains to show the uniqueness. This amounts to show that a non-zero entire function $f(z)$ such that $\overline{f(\bar z)}=f(z)$ and $f(n)=(-1)^nQ(n)$ for all $n\in \Z$ and for some polynomial $Q(n)$ automatically fulfills \eqref{mingrowth}. For $x\in \R$ one has $f(x)\in \R$. Let $s_k n^k$ be the leading term of $Q(n)=\sum s_j n^j$. Then, outside a finite
set of indices $n$ the sign of $Q(n)$ is independent of $n$ for $n>0$ and for $n<0$. It follows that the function $f$ changes sign in each interval
$[n, n+1]$ and hence admits at least one zero in such intervals. Moreover, $N(0, R)=\displaystyle{\sum_{f(z)=0\atop |z|<R}} \log|\frac R z|$ fulfills
$$
\underline{\lim}~\frac{N(0,R)}{R}\geq 2
$$
and hence, by Theorem \ref{main1} one gets
$$
    \underline \lim~\frac{T(R,f)}{R}\geq 2.
$$
The uniqueness follows.
\endproof

\begin{defn}\label{canext} Let $N(n)$ be a real function on $\Z$ of the form \eqref{data}. Then the canonical extension $N(z)$, $z\in \C$ is given by
\begin{equation}\label{canextdefn}
    N(z)=f(z)+\cos(\pi z)Q(z)
\end{equation}
where $f$ and $Q$ are uniquely determined by Proposition \ref{propextfun}.
\end{defn}\vspace{.05in}

We shall explicitly compute this extension in the case of the function $N_X(q)$ for $X=\Sp(\F_1[H])$.\vspace{.05in}

\subsubsection{The counting function for $X=\Sp(\F_1[H])$}\hfill \vspace{.05in}

Let $H$ be a finitely generated abelian group. In this subsection we compute the counting function $N_X(n)$ for $X=\Sp(\F_1[H])$. We start by considering the cyclic case $H=\Z/m\Z$.

\begin{lem}\label{fourierexp} Let $X=\Sp(\F_{1^m})$. Then, the canonical extension of the counting function $N_X(n)$
 is given by
\begin{equation}\label{canext1}
  N(z)=  \sum_{d|m}\frac{\varphi(d)}{d}\left(\sum_{|k|<d/2} e^{2i\pi \frac{(z-1)k}{d}}+\epsilon_d\, \cos(\pi(z-1))\right)\,, \ \ z\in\C
\end{equation}
where $\epsilon_d=1$ if $d\in\Z$ is even and $\epsilon_d=0$ otherwise.
\end{lem}

\proof It follows from \eqref{genform0} that $N(n)=\gcd(n-1,m)$. Moreover, one also knows that
$$
\gcd(n,m)=\sum_{d|n\atop d|m}\varphi(d)
$$
where $\varphi$ is the Euler totient function. Thus
\begin{equation}\label{basicform}
\gcd(n,m)=\sum_{d|m}\frac{\varphi(d)}{d}\sum_{} e^{2i\pi \frac{nk}{d}}
\end{equation}
where the sum $\sum e^{2i\pi \frac{nk}{d}}$ is taken over all characters of the
additive group $\Z/d\Z$ evaluated on $n$ mod. $d$. This sum can be written as
$$
\sum e^{2i\pi \frac{nk}{d}}=\sum_{|k|<d/2} e^{2i\pi \frac{nk}{d}}+\epsilon_d\, \cos(\pi n)
$$
which gives the canonical extension of $N(n)$ after a shift of one in the variable.
\endproof\vspace{.05in}
\begin{figure}
\begin{center}
\includegraphics[scale=0.8]{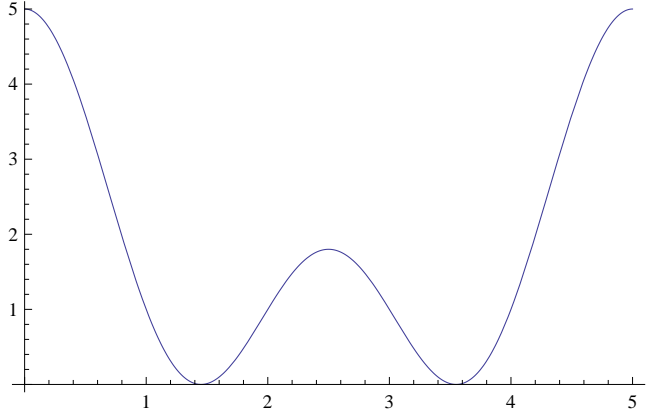}
\end{center}
\caption{The graph of the canonical extension of $N(q)$, for $X=\Sp(\F_1[\Z/5\Z])$. \label{graph1} }
\end{figure}

There is a conceptual explanation for the expansion \eqref{canext1}. Let $C_d$ be the $d$-th cyclotomic polynomial
$$
C_d(x)=\prod_{r|d}(x^{d/r}-1)^{\mu(r)}.
$$
One has the decomposition into irreducible factors over $\Q$
$$
x^m-1=\prod_{d|m}C_d(x)
$$
which determines the decomposition
$$
\Q[\Z/m\Z]=\bigoplus_{d|m}\K_d
$$
where $\K_d=\Q(\xi_d)$ is the cyclotomic field extension of $\Q$ by a primitive $d$-th root of $1$.
This decomposition corresponds to a factorization of the zeta function
$$
\zeta_X(s)=\prod_p \zeta_{X/\F_p}(s)
$$
as a product
$$
\zeta_X(s)=\prod_{d|m} \zeta_{X/\K_d}(s)
$$
where each term of the product is in turns a product of $L$-functions
$$
\zeta_{X/\K_d}(s)=\prod_{\chi\in \hat G(d)}L(s,\chi),\qquad G(d)=(\Z/d\Z)^*
$$
and the order of $G(d)$ is $\varphi(d)$. The decomposition
 $$
\zeta_X(s)=\prod_{d|m} \prod_{\chi\in \hat G(d)}L(s,\chi)
$$
corresponds precisely to the decomposition \eqref{basicform}.\vspace{.05in}

Let now $H=\prod_j \Z/m_j\Z$ be a finite abelian group. Then, the extension of the group ring
$\Z[H]$ over $\Q$ can be equivalently interpreted as the tensor product
$$
\Q[H]=\bigotimes_j(\bigoplus_{d|m_j}\K_d).
$$
In fact, by applying the Galois correspondence,  the decomposition of $H$ into finite cyclic groups corresponds to the disjoint union of orbits for the action of the Galois group $\hat\Z^*$ on the cyclotomic extension $\Q^{{\rm ab}}/\Q$ (\ie the maximal abelian extension of $\Q$). If we denote the group $H$ additively, two elements $x,y\in H$ are on the same Galois orbit if and only if $\{nx\,|\,n\in \Z\}=\{ny\,|\,n\in \Z\}$ \ie if and only if they  generate the same subgroup of $H$. Thus, the set $Z$ of orbits   is the set of cyclic subgroups of $H$.

For a cyclic group $C$ of order $d$ we let, as in \eqref{canext1}
\begin{equation}\label{cyclicgr}
 N(z,C)=  \frac{\varphi(d)}{d}\left(\sum_{|k|<d/2} e^{2i\pi \frac{(z-1)k}{d}}+\epsilon_d\, \cos(\pi(z-1))\right).
\end{equation}

Then, we obtain the following

\begin{lem}\label{fourierexpbis} Let $\Gamma=\Z^k\times H$ be a finitely generated abelian group ($H$ is a finite abelian torsion group) and $X=\Sp(\F_1[\Gamma])$. Then, the canonical extension of the counting function $N_X(q)$
 is given by
\begin{equation}\label{canext2}
  N(z)=  (z-1)^k\sum_{C\subset H} N(z,C),\qquad z\in\C
\end{equation}
where the sum is over the cyclic subgroups $C\subset H$.
\end{lem}

\proof  It is enough to show that the function $N(n+1)$ as in  \eqref{canext2} gives the  number of homomorphisms from $\Gamma$ to $\Z/n\Z$. Indeed, as a function of $z$ it is already in the canonical form of Definition  \ref{canext}. Then, by duality we just need to compare  the number $\iota(n,H)$ of homomorphisms $\rho$ from $\Z/n\Z$ to $H$ with $\sum_{C\subset H} N(n+1,C)$. The cyclic  subgroup $C={\rm Im}\,\rho$ is uniquely determined, thus
$$
\iota(n,H)=\sum_{C\subset H}\mu(n,C)
$$
where $\mu(n,C)$ is the number of surjective homomorphisms from $\Z/n\Z$ to $C$. If $d$ denotes the order of $C$, this number is $0$ unless $d|n$ and in the latter case it is $\varphi(d)$, thus using \eqref{cyclicgr} it coincides in all cases with $N(n+1,C)$.\endproof\vspace{.05in}

\subsection{An integral formula for the logarithmic derivative of $\zeta_N(s)$}\hfill \vspace{.05in}\label{integral}

Let $N(q)$ be a real-valued continuous function on $[1,\infty)$ satisfying a polynomial bound $|N(q)|\leq C q^k$ for some finite positive integer $k$ and a fixed positive integer constant $C$. Then the associated generating function is
$$
Z(q,T)={\rm exp}\left( \sum_{r\geq 1}N(q^r)T^r/r \right)
$$
and one knows that the power series $Z(q,q^{-s})$ converges for $\Re(s)>k$. The zeta function over $\F_1$ associated to $N(q)$ is defined as follows
$$
\zeta_N(s):=\lim_{q\to 1}(q-1)^\chi\,Z(q,q^{-s}) \,, \ \ \chi=N(1).
$$
This definition requires some care to assure its convergence. To eliminate the ambiguity in the extraction of the finite part, one works instead with the logarithmic derivative
\begin{equation}\label{normael}
    \frac{\partial_s\zeta_N(s)}{\zeta_N(s)}=-\lim_{q\to 1} F(q,s)
\end{equation}
where
\begin{equation}\label{fqsdefn}
F(q,s)=-\partial_s \sum_{r\ge 1}N(q^r)\frac{q^{-rs}}{r}\,.
\end{equation}

One then has (\cf \cite{announc3})

\begin{lem} \label{compute1}With the above notations and for $\Re(s)>k$
\begin{equation}\label{lim}
    \lim_{q\to 1} F(q,s) = \int_1^\infty N(u)u^{-s}d^*u\,,\ \ d^*u=du/u
\end{equation}
and
\begin{equation}\label{logzetabis}
    \frac{\partial_s\zeta_N(s)}{\zeta_N(s)}=-\int_1^\infty  N(u)\, u^{-s}d^*u\,.
\end{equation}
\end{lem}\vspace{.05in}

\subsection{Analyticity of the integrals}\hfill \vspace{.05in}\label{analytici}

 The statement of Lemma~\ref{compute1} and the explicit form of the function $N(q)$ as in Lemma~\ref{fourierexp} clearly indicate that one needs to analyze basic integrals of the form
\begin{equation}\label{anfun0}
f(s,a)=\int_1^\infty e^{iau} u^{-s}d^*u.
\end{equation}

\begin{lem}\label{analytic} For $a>0$
\begin{equation}\label{anfun}
    f(s,a)=a^s\int_a^\infty\,e^{iu}u^{-s}d^*u
\end{equation}
defines an entire function of $s\in \C$.
\end{lem}
\proof For $\Re(s)>0$, the integral in \eqref{anfun} is absolutely convergent. Integrating by parts one obtains
$$
\int_a^\infty\,e^{iu}u^{-s-1}du+\int_a^\infty\,\frac 1i e^{iu}(-s-1)u^{-s-2}du=[\frac 1i e^{iu}u^{-s-1}]_a^\infty
$$
which supplies the equation
\begin{equation}\label{ancont}
af(s,a)=-i(s+1)f(s+1,a)+ie^{ia}.
\end{equation}
This shows that $f$ is an entire function of $s$, since after iterating \eqref{ancont} the variable $s$ can be moved to the domain of absolute convergence.\endproof\vspace{.05in}

The function $f(s,a)$ as in \eqref{anfun0} can be expressed in terms of hypergeometric functions in the following way
$$
f(s,a)=e^{-\frac 12 i\pi s}a^s\Gamma(-s)+ i\frac{a}{s-1}H[\frac{1-s}{2}, (\frac 32, \frac{3-s}{2}), -\frac{a^2}{
      4}] +\frac{1}{s}H[{-\frac s2}, (\frac 12, \frac{2-s}{2}), -\frac{a^2}{
      4}]
      $$
where the hypergeometric function $H$ is defined by
$$
H(\alpha,(\beta,\gamma),z)=\sum_{k=0}^\infty\frac{(\alpha)_k}{(\beta)_k(\gamma)_k}\frac{z^k}{k!}.
$$
One sets $\displaystyle{(x)_k=\prod_{j=0}^{k-1} (x+j)}$. Note that the  formulas simplify as follows
$$
\frac{1}{s-1}H[\frac{1-s}{2}, (\frac 32, \frac{3-s}{2}), -\frac{a^2}{
      4}]=\sum_{k=0}^\infty\frac{(-a^2)^{k}}{(2k+1)! (s-(2k+1))}
$$
and
$$
\frac 1 s H[{-\frac s2}, (\frac 12, \frac{2-s}{2}), -\frac{a^2}{
      4}]=\sum_{k=0}^\infty\frac{(-a^2)^{k}}{(2k)! (s-2k)}.
$$
By using \eqref{anfun0} and expanding the exponential in powers of $a$ in the integral
$$
\int_0^1 e^{iau} u^{-s}d^*u=\sum_{n=0}^\infty\int_0^1 \frac{(iau)^n}{n!} u^{-s}d^*u
$$
the  resulting expression for $f(s,a)$ in terms of hypergeometric functions becomes
\begin{equation}\label{hyperexpr}
f(s,a)=e^{-\frac 12 i\pi s}a^s\Gamma(-s)+\sum_{n=0}^\infty\frac{(ia)^{n}}{n! (s-n)}.
\end{equation}

Thus one gets
\begin{lem}\label{analyticbis} For $a>0$, the following identity holds
\begin{equation}\label{anfunbis}
   f(s,a)=e^{-\frac 12 i\pi s}a^s\Gamma(-s)-\sum_{{\rm poles}}\frac{{\rm residue}}{s-n}.
\end{equation}
\end{lem}
\proof The last term of the  expression \eqref{hyperexpr} is a sum of the form $\sum_{{\rm poles} \;g}\frac{{\rm residue}(g)}{s-n}$ where $g(s)=e^{-\frac 12 i\pi s}a^s\Gamma(-s)$.
This suffices to conclude the proof since one knows from lemma~\ref{analytic} that $f(s,a)$ is an entire function of $s$. One can also check directly that the expression for the residues as in \eqref{hyperexpr} is correct by using \eg the formula of complements.
\endproof\vspace{.05in}

\subsection{Zeta function of Noetherian $\F_1$-schemes}\hfill \vspace{.05in}\label{zetasch}

The above discussion has shows that the elementary constituents of the zeta function of a Noetherian $\F_1$-scheme are of the following type
\begin{equation}\label{elemform}
    \xi_d(s)=s^{-\frac{\varphi(d)}{d}}e^{-\frac{\varphi(d)}{d} H_d(s)}
\end{equation}
where $H_d(s)$ is the primitive, vanishing for $\Re(s)\to \infty$, of the entire function
\begin{equation}\label{prim}
    \partial_s H_d(s)=\sum_{|k|\leq d/2\atop k\neq 0} e^{-i\frac {2\pi k}{d}}  f(s,\frac {2\pi k}{d}).
\end{equation}
When $d$ is even, one divides by  $2$ the contribution of $k=\frac d2$.
Note that, by Lemma~\ref{analyticbis}, the function $\partial_s H_d(s)$ is obtained by applying the transformation
$$
F\mapsto F-\sum_{{\rm poles}}\frac{\hbox{ residue}}{s-n}
$$
to the function
$$
  F_d(s)=  \Gamma(-s)\sum_{|k|\leq d/2\atop k\neq 0} \cos(\frac {2\pi k}{d}+\frac{\pi s}{2})(2\pi\frac kd)^s.
$$

\begin{figure}
\begin{center}
\includegraphics[scale=0.8]{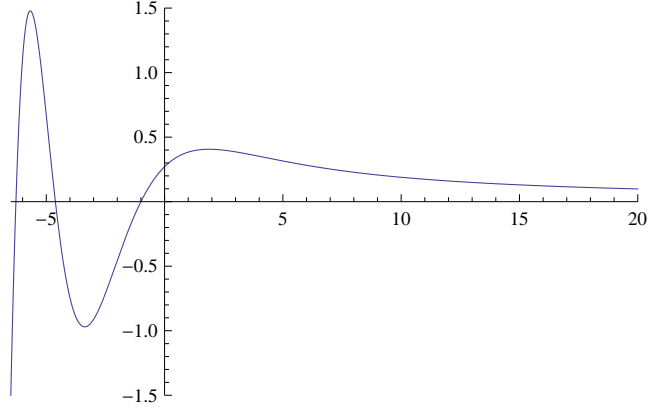}
\end{center}
\caption{Graph of the function $\partial_s H_d(s)$ for $d=3$. \label{graph2} }
\end{figure}

\begin{thm}\label{dthmfonesch1}
Let $X$ be a  Noetherian $\F_1$-scheme. Then the zeta function of $X$ over $\F_1$ is a finite product
\begin{equation}\label{zetaX}
    \zeta_X(s)=\prod_{x\in X}\zeta_{X_x}(s)
\end{equation}
where $X_x$ denotes the constituent of $X$ over the residue field $\kappa(x) = \F_1[\mathcal O_{X,x}^*]$ and $\mathcal O_{X,x}^*=\Gamma$ is a finitely generated abelian group such as $\Gamma=H\times \Z^n$ ($H=$ finite torsion group). The zeta function of each constituent $X_x$ takes the following form
\begin{equation}\label{prodzeta1}
\zeta_{X_x}(s)=\prod_{k\ge 0} \zeta(\F_1[H],s-n+k)^{(-1)^k{n\choose k}}
\end{equation}
where $\zeta(\F_1[H],s)$ is the product
\begin{equation}\label{prodzeta2}
    \zeta(\F_1[H],s)=\prod_{d||H|}\xi_d(s)^{\gamma(H,d)}
\end{equation}
where  $\gamma(H,d)$ denotes the number of cyclic subgroups $C\subset H$ of order $d$.
\end{thm}

\proof The proof of \eqref{zetaX} follows the same lines as that given in Theorem~\ref{dthmfonesch}. The expression \eqref{prodzeta1}
for $\zeta_{X_x}(s)$ in terms of $\zeta(\F_1[H],s)$ follows from \eqref{canext2} and the simple translation in $s$ that the multiplication of the counting function by a power of $q$ generates using Lemma~\ref{compute1}. By \eqref{canext2} and \eqref{cyclicgr} it remains to show that the zeta function $Z_d(s)$ associated to the counting function \eqref{cyclicgr} coincides with $\xi_d(s)$ defined in \eqref{elemform}. By Lemma~\ref{compute1} and \eqref{logzetabis} the logarithmic derivative of $Z_d(s)$ is
$$
\frac{\partial_s Z_d(s)}{Z_d(s)}=-\int_1^\infty    \frac{\varphi(d)}{d}\left(\sum_{|k|<d/2} e^{2i\pi \frac{(u-1)k}{d}}+\epsilon_d\, \cos(\pi(u-1))\right)u^{-s}d^*u.
$$
The term coming from $k=0$ contributes with $-\frac{\varphi(d)}{d}\frac 1s$ and thus corresponds to the fractional power $\displaystyle{s^{-\frac{\varphi(d)}{d}}}$ in \eqref{elemform}. The other terms contribute with integrals of the form
$$
I(s,\alpha)=\int_1^\infty  e^{2\pi i\alpha(u-1)}\, u^{-s}d^*u.
$$
One has
$$
I(s,\alpha)=e^{-2\pi i\alpha}\int_1^\infty  e^{2\pi i\alpha u}\, u^{-s}d^*u
=e^{-2\pi i\alpha}(2\pi \alpha)^s\int_{2\pi \alpha}^\infty  e^{i v}\, v^{-s}d^*v.
$$
Thus with the notation \eqref{anfun} one obtains
\begin{equation}\label{fineform}
I(s,\alpha)=e^{-ia}  f(s,a)\,, \ a=2\pi \alpha.
\end{equation}
The equality $Z_d(s)=\xi_d(s)$ follows using \eqref{prim}.\endproof\vspace{.05in}

\begin{cor}\label{corzeta} Let $X$ be a  Noetherian $\F_1$-scheme. Then the zeta function of $X$
is the product
\begin{equation}\label{zetaXbis}
    \zeta_X(s)=e^{h(s)}\prod_{j=0}^n(s-j)^{\alpha_j}
\end{equation}
of an exponential of an entire function $h(s)$ by a finite product of fractional powers of simple monomials. The exponents $\alpha_j$ are rational numbers given explicitly by
\begin{equation}\label{expo}
    \alpha_j=(-1)^{j+1}\sum_{x\in X}(-1)^{n(x)}{n(x)\choose j}\sum_{d}\frac 1d\nu(d,\cO_{X,x}^\times)
\end{equation}
where $n(x)$ is the local dimension (\ie the rank of $\cO_{X,x}^\times$), and $\nu(d,\cO_{X,x}^\times)$ is the number of  injective homomorphisms from $\Z/d\Z$ to $\cO_{X,x}^\times$.
\end{cor}

\proof By \eqref{elemform} one has
$$
\zeta(\F_1[H],s)=\prod_{d||H|}\xi_d(s)^{\gamma(H,d)}=s^{-\epsilon_H}e^{-\sum_{d||H|}\frac{\varphi(d)\gamma(H,d)}{d} H_d(s)}
$$
where
\begin{equation}\label{expobis}
    \epsilon_H=\sum_{d||H|}\frac{\varphi(d)\gamma(H,d)}{d}.
\end{equation}
Next, by \eqref{prodzeta1} one has
$$
\zeta(\F_1[\Gamma],s)=\prod_{k\ge 0} \zeta(\F_1[H],s-n+k)^{(-1)^k{n\choose k}}
$$
which is a product of the form
$$
\zeta(\F_1[\Gamma],s)=\prod_{k\ge 0} (s-j)^{-(-1)^{n-j}{n\choose j}\epsilon_H}e^{h(s)}
$$
where $h(s)$ is an entire function. Thus one gets that the exponent $\alpha_j$ of $(s-j)$ is
$$
    \alpha_j=(-1)^{j+1}\sum_{x\in X}(-1)^{n(x)}{n(x)\choose j}\sum_{d||H_x|}\frac{\varphi(d)\gamma(H_x,d)}{d}
$$
where $n(x)$ is the local dimension and $H_x$ is the finite group in the decomposition of $\cO_{X,x}^\times$. The above expression  gives \eqref{expo} since $\varphi(d)\gamma(H_x,d)$ is the number of injective homomorphisms from $\Z/d\Z$ to $\cO_{X,x}^\times$.\endproof

\section{Beyond $\F_1$-schemes}\label{ellipticsect}

One may naturally wonder on the adeptness of the method of the computation of the zeta function over $\F_1$ that we have described in the above section, for a possible generalization beyond the case of $\Z$-schemes associated to Noetherian $\F_1$-schemes.
 In \cite{announc3}, we have already given evidence on the effectiveness of the integral formula of Lemma
 \ref{compute1} of \S \ref{integral}, by showing that there exists a {\em uniquely} defined distribution $N(u)$ on $[1,\infty)$ which describes the counting function for the hypothetical curve $C=\overline {\Sp\,\Z}$ over $\F_1$ whose zeta function $\zeta_C(s)$ (over $\F_1$) coincides with the completed Riemann zeta function  $\zeta_\Q(s)=\pi^{-s/2}\Gamma(s/2)\zeta(s)$ over $\C$.\vspace{.05in}

 The goal of this section is to show that as long as one is only interested in the singularities of the zeta function, one may, in the case of Noetherian $\F_1$-schemes, replace the integral which appears in \eqref{logzetabis}  by the corresponding infinite sum \ie
\begin{equation}\label{replace}
    -\int_1^\infty  N(u)\, u^{-s}d^*u~\rightsquigarrow~ -\sum_{n\ge 1} N(n)\, n^{-s-1},
\end{equation}
which is now defined without any need for extrapolation.
This substitution provides one with a natural definition of a modified zeta function that is defined by the formula
\begin{equation}\label{modifiedze}
\frac{\partial_s\zeta_N^{\rm disc}(s)}{\zeta_N^{\rm disc}(s)}= -\sum_{n=1}^\infty  N(n)\, n^{-s-1}.
\end{equation}

In \S \ref{ellipticsubsect} we compute this function for a {\it particular} choice of the extension of the counting function $N(q)$  in the case of elliptic curves over $\Q$.

\subsection{Modified zeta function}\hfill \vspace{.05in}\label{modizeta}

Within the category of the Noetherian $\F_1$-schemes, one has the following result

\begin{prop}\label{replacelem} The replacement \eqref{replace} in the definition of the logarithmic derivative $\frac{\partial_s\zeta_N(s)}{\zeta_N(s)}$ as in \eqref{logzetabis}, does not alter the singularities of the zeta function  associated to a Noetherian $\F_1$-scheme. In fact, the ratio of these two zeta functions (\ie the logarithmic derivatives of $\zeta_N(s)$ and $\zeta_N^{\rm{disc}}(s)$) is the exponential of an entire function.
\end{prop}

\proof It is enough to treat the case of a counting function $N(u)$, of the real variable $u$, of the form
\begin{equation}\label{refsimple}
N(u)= \frac{\varphi(d)}{d} \,u^\ell \left(\sum_{|k|<d/2} e^{2i\pi \frac{(u-1)k}{d}}+\epsilon_d\, \cos(\pi(u-1))\right)
\end{equation}
and show that the  following is an entire function of $s\in\C$
$$
h(s)=\int_1^\infty  N(u)\, u^{-s}d^*u -\sum_{n=1}^\infty  N(n)\, n^{-s-1}.
$$
Notice that the term $u^\ell$ generates a shift in $s$ and thus we can assume that $\ell=0$. Let $G(d)=(\Z/d\Z)^*$ be the multiplicative group of residues modulo $d$ which are prime to $d$. By extending the characters $\chi\in \hat G(d)$ by $0$ on residues modulo $d$ which are not prime to $d$, one obtains the equality
$$
\sum_{\chi\in \hat G(d)} \chi(n)=\begin{cases} \varphi(d)&\text{if}~n=1\\
0&\text{if}~n\neq 1.\end{cases}
 $$
Thus, the restriction of $N(u)$ to the integers agrees with the following sum
$$
N(n)=\sum_{\chi\in \hat G(d)} \chi(n)
$$
and one gets
$$
\sum_{n\ge 1}  N(n)\, n^{-s-1}=\sum_{\chi\in \hat G(d)}L(\chi,s+1).
$$
When $\chi\neq 1$, the function $L(\chi,s+1)$ is entire, thus the only singularity
 arises from the function $L(1,s+1)$ which is known to have a unique pole at $s=0$ with residue $\prod_{p|d}(1-\frac 1p)$. The only singularity of the integral $\int_1^\infty  N(u)\, u^{-s}d^*u$
is due to the contribution of the constant term \ie
$$
\int_1^\infty   \frac{\varphi(d)}{d}\, u^{-s}d^*u= \frac{\varphi(d)}{d}\frac 1s.
$$
Thus the equality $\prod_{p|d}(1-\frac 1p)=\frac{\varphi(d)}{d}$ shows that the function $h(s)$ is entire.
\endproof

We can therefore adopt the following definition
\begin{defn}\label{modifiedzeta} Let $X$ be an $\F_1$-scheme. The modified zeta function $\zeta_X^{\rm disc}(s)$ is defined by
\begin{equation}\label{modified}
   \frac{\partial_s\zeta_X^{\rm disc}(s)}{\zeta_X^{\rm disc}(s)}=-\sum_{n\ge 1}  N(n)\, n^{-s-1}
\end{equation}
where $N(n+1)=\# X(\F_{1^n})$.
\end{defn}
By Proposition~\ref{replacelem}, the singularities of this modified zeta function are the same as the singularities of $\zeta_X(s)$: \cf  Corollary~\ref{corzeta}. The advantage of this more general definition is that it requires no choice of an interpolating function.
This means that one can define $\zeta_X^{\rm disc}(s)$, up to a multiplicative normalization factor, by just requiring some polynomial control on the size of growth  of the finite set $X(\F_{1^n})$.\vspace{.05in}

\begin{rem}\label{superf} {\rm There is a superficial resemblance between the right hand side of  \eqref{modified} that defines the logarithmic derivative of $\zeta_X^{\rm disc}(s)$ and the ``absolute" Igusa zeta function associated in
\cite{abs} to an $\F_1$-scheme. However, there is a shift of $1$ in the argument $n$ which is not easy to eliminate and that introduces a fundamental difference. Thus, in terms of points defined  over $\F_{1^n}$ we obtain
\begin{equation}\label{diffform}
    \sum_{n\ge 1}  N(n)\, n^{-s-1}=\sum_{n\ge 1}  X(\F_{1^{n-1}})\, n^{-s-1}.
\end{equation}
On the other hand, the zeta function associated in
\cite{abs} is given by
\begin{equation}\label{diffform1}
    \zeta^I(s,X)=\sum_{n\ge 1}  X(\F_{1^{n}})\, n^{-s}.
\end{equation}
The shift of $1$ in $s$ is easily accounted for, but not the replacement of $X(\F_{1^{n-1}})$ by $X(\F_{1^{n}})$. In the next subsection we will show how this shift manifests itself in a concrete example.
}\end{rem}\vspace{.05in}

\subsection{Modified zeta function and $\overline {\Sp\,\Z}$ }\hfill \vspace{.05in}\label{riemannzeta}

 In \cite{announc3}, we have shown how to determine the counting distribution $N(q)$, defined for $q\in [1,\infty)$ in such a way that the zeta function as in \eqref{logzetabis} gives the complete Riemann zeta function, \ie so that the following equation holds
\begin{equation}\label{logzetabisnew}
    \frac{\partial_s\zeta_\Q(s)}{\zeta_\Q(s)}=-\int_1^\infty  N(u)\, u^{-s}d^*u\,,
\end{equation}
where
$$
\zeta_\Q(s)=\pi^{-s/2}\Gamma(s/2)\zeta(s).
$$
Moreover, using convergence in the sense of distributions one obtains
 \begin{equation}\label{inghamn1}
N(q)=q- \sum_{\rho\in Z}{\rm order}(\rho)\, q^{\rho} +1.
\end{equation}
This is in perfect analogy with the theory for function fields. Indeed, for a smooth curve $X$ over a finite field $\F_{q_0}$ and for any power $q=q_0^\ell$ of the order $q_0$ of the finite field of constants one has
$$
\#X(\F_q)=N(q)=q- \sum_{\rho\in Z}{\rm order}(\rho)\, q^{\rho} +1.
$$
In the above formula, we have written the eigenvalues of the Frobenius operator acting on the cohomology of the curve
in the form $\alpha=q_0^\rho$.

While these results supply a strong indication on the coherence of the quest for an arithmetic theory over $\F_1$, they also indicate that the similarities with the geometry over finite fields are also elusive and that the $\F_1$-counterpart of the classical theory will be richer in new ingredients. In this section, we will study what type of additional information one can gather about the hypothetical curve $C=\overline {\Sp\,\Z}$ by using the modified zeta function, \ie the simplified formula \eqref{replace}. Since this simplified formula clearly neglects the contribution of the archimedean place, it seems natural to replace the complete zeta function $\zeta_\Q(s)$ by the Riemann zeta function $\zeta(s)$. Thus, the simplified form of \eqref{logzetabisnew} takes the form
\begin{equation}\label{simplerzeta}
   \frac{\partial_s\zeta(s)}{\zeta(s)}=-\sum_1^\infty  N(n)\, n^{-s-1}.
\end{equation}
This supplies for the counting function $N(n)$ the formula
\begin{equation}\label{Nnform}
    N(n)=n\Lambda(n)\,,
\end{equation}
where  $\Lambda(n)$ is the von-Mangoldt function with value $\log p$ for powers $p^\ell$ of  primes and zero otherwise. We now use the relation $N(n+1)=\# X(\F_{1^n})$ of Definition \ref{modifiedzeta} and conclude that the sought for ``scheme over $\F_1$" should fulfill the following condition
\begin{equation}\label{ideal}
   \#(X(\F_{1^n}))=\left\{
                              \begin{array}{ll}
                                0 & \hbox{if}\ n+1 \ \hbox{is not a prime power}\\
                                (n+1)\log p, & \hbox{if}\ n=p^\ell-1,\ p \ \hbox{prime.}
                              \end{array}
                            \right.
\end{equation}
The two obvious difficulties are that $(n+1)\log p$ is not an integer, and that multiples $m=nk$ of numbers of the form $n=p^\ell-1$ are no longer of a similar form, while there are morphisms from $\F_{1^n}$ to $\F_{1^m}$. This shows that functoriality can only hold in a more restrictive manner.
Moreover, at first sight, it seems hopeless to find a natural construction starting from a monoid $M$ which singles out among the monoids $M=\F_1[\Z/n\Z]=\Z/n\Z\cup\{0\}$ those for which $n=p^\ell-1$, for some prime power. However, this is precisely what is achieved by the study of the additive structure in \S \ref{addi}. Indeed, one has
\begin{cor}\label{corfine} Let $M=\F_1[\Z/n\Z]=\Z/n\Z\cup\{0\}$ viewed as a commutative (pointed) monoid. Then,  $n=p^\ell-1$ for a prime power $p^\ell$  if and only if the following set $A(M)$ is non-empty. $A(M)$ is the set of maps $s:M\to M$ such that\vspace{.05in}

 $\bullet$~$s$ commutes with its conjugates by multiplication by elements of $M^\times$\vspace{.05in}

 $\bullet$~$s(0)=1$.
\end{cor}
\proof This would follow from Corollary \ref{uniqueness} if we assumed that $s$ is bijective. So, we need to prove that this assumption is unnecessary. The proof of Theorem \ref{thmsym} shows that, for $s\in A(M)$, $M$ becomes a semi-field, where we use \eqref{mainadddefn} to define the addition. Moreover since $M$ is finite, if it is not a field $M$ contains a (prime) semi-ring  $B(n,i)$ for finite $n$ and $i$. As explained after Proposition \ref{propprime}, this implies $n=2$ and $i=1$ so that $M$ is  an idempotent semi-field. Finally, by Proposition \ref{uniquesemi} one gets that unless $M$ is a field it is equal to $\B$. Thus we obtain that the cardinality of $M$ is a prime power. \endproof

The cardinality of the set $A(M)$ can be explicitly computed. By using Corollary \ref{uniqueness} and the fact that the automorphisms of the finite field $\F_q$ are powers of the Frobenius, one gets for\footnote{one has $\#(A(\F_{1}))=2$.} $n>1$
\begin{equation}\label{cardxm}
     \#(A(\F_{1^n}))=\left\{
                              \begin{array}{ll}
                                0 & \hbox{if}\ n+1 \ \hbox{is not a prime power}\\
                                \varphi(p^\ell-1)/\ell, & \hbox{if}\ n=p^\ell-1,\ p \ \hbox{prime}
                              \end{array}
                            \right.
\end{equation}
where $\varphi$ is the Euler totient function. This produces a result which, although not exactly equal to the (non-integer) value as in \eqref{ideal}, comes rather close. Moreover, one should acknowledge the fact that even though the definition of $A(M)$ is natural it is only functorial for isomorphisms and not for arbitrary morphisms.\vspace{.05in}

\subsection{Elliptic curves}\hfill \vspace{.05in}\label{ellipticsubsect}

Definition~\ref{modifiedzeta} applies any time one has a reasonable guess for the definition of an extension of the counting function of points
from the case $q=p^\ell$  to the case of an arbitrary positive integer. For elliptic curves $E$ over $\Q$, one may use the associated modular form
\begin{equation}\label{mod2}
F_E(\tau)=\sum_{n=1}^\infty a(n)q^n,\qquad q=e^{2\pi i \tau}
\end{equation}
since the sequence $a(n)$ of the coefficients of the above series is defined for all values of $n$. Moreover, these Fourier coefficients  also fulfill equations
\begin{equation}\label{mod1}
N(p,E)=\#E(\F_p) = p + 1 - a(p)
\end{equation}
 at each prime number $p$ of good reduction for the curve $E$. This choice is, however, a bit too naive since \eqref{mod1} does not continue to hold at primes powers.\vspace{.05in}

 We define the function $t(n)$ as the only multiplicative function which agrees with $a(p)$ at each prime and for which \eqref{mod1} continues to hold for prime powers. Then we have the following result

\begin{lem}\label{taste}
Let $t(m)$ be the only multiplicative function such that $t(1)=1$ and
\begin{equation}\label{sumfrob}
    N(q,E)=q+1-t(q)
\end{equation}
for any prime power $q$. Then, the generating function
\begin{equation}\label{betL}
    R(s,E)=\sum_{n\ge 0} t(n)\,n^{-s}
\end{equation}
has the following description
\begin{equation}\label{taste1}
    R(s,E)=\frac{L(s,E)}{\zeta(2s-1)M(s)}.
\end{equation}
$L(s,E)$ is the $L$-function of the elliptic curve, $\zeta(s)$ is the Riemann zeta function, and
\begin{equation}\label{taste5}
    M(s)=\prod_{p\in S}(1-p^{1-2s})
\end{equation}
is a finite product of simple factors over the set $S$ of
primes of  bad reduction for $E$.
\end{lem}
\proof
One has an Euler product for the $L$-function $L(s,E)=\sum a(n)n^{-s}$ of the elliptic curve of the form
\begin{equation}\label{Lfunc}
L(s,E)=\prod_p L_p(s,E).
\end{equation}
For almost all primes, \ie  for $p\notin S$ with $S$ a finite set
\begin{equation}\label{goodreducpr}
L_p(s,E)=\left((1-\alpha_p p^{-s})(1-\bar\alpha_p p^{-s})\right)^{-1}=\sum_{n=0}^\infty a(p^n)p^{-ns}
\end{equation}
where $\alpha_p=\sqrt{p}~e^{i\theta_p}$.
Since the function $t(n)$ is multiplicative one has an Euler product
\begin{equation}\label{localg}
R(s,E)=\sum_{n\ge 0} t(n)\,n^{-s}=\prod_p\left(  \sum_{n\ge 0} t(p^n)p^{-ns}\right)=\prod_p R_p(s,E).
\end{equation}
Let $p\notin S$ be a prime of good reduction. Then for $q=p^\ell$, one has
\begin{equation}\label{frobeigen}
N(q,E)=q+1-\alpha_p^\ell-\bar\alpha_p^\ell.
\end{equation}
We now show that
 \begin{equation}\label{factorRZ}
(1-p^{1-2s})\,    \sum_{n= 0}^\infty a(p^n)p^{-ns}= \sum_{n=0}^\infty t(p^n)p^{-ns}.
\end{equation}
This relation follows  from \eqref{goodreducpr} and the equality
$$
\frac{(1-p \,x^2)}{(1-\alpha_p x)(1-\bar\alpha_p x)}=\frac{1}{(1-\alpha_p x)}+\frac{1}{(1-\bar\alpha_p x)}-1
$$
for $x=p^{-s}$. The above formula appears in the expression of the Poisson kernel (\cf \cite{Rudin} \S\S~5.24,  11.5)
$$
P_r(\theta)=\sum_{-\infty}^\infty r^{|n|}e^{in\theta}=\frac{1-r^2}{1-2r\cos\theta+r^2}.
$$
Using \eqref{factorRZ} one gets
\begin{equation}\label{localf}
R_p(s,E)=(1-p^{1-2s})  \sum_{n=0}^\infty a(p^n)p^{-ns}.
\end{equation}
      When $p\in S$, the local factor $L_p(s,E)$ of the $L$-function takes one of the following descriptions (\cf \cite{Silverman} Appendix C \S 16)\vspace{.05in}

$\bullet$~$L_p(s,E)=\displaystyle{\sum_{n=0}^\infty p^{-ns}}$ when $E$ has split multiplicative reduction at $p$.\vspace{.05in}

$\bullet$~$L_p(s,E)=\displaystyle{\sum_{n=0}^\infty (-1)^np^{-ns}}$ when $E$ has non-split multiplicative reduction at $p$.\vspace{.05in}

$\bullet$~$L_p(s,E)=1$ when $E$ has additive reduction at $p$.\vspace{.05in}

We prove  now that in any of the above cases one has
\begin{equation}\label{checkequ}
R_p(s,E)=L_p(s,E).
\end{equation}

One knows that in the singular case, the elliptic curve has exactly one singular point. This point is already defined over $\F_p$
since if a cubic equation has a double root this root is rational (\cf~Appendix~\ref{sing23}). Thus, for any power
$q=p^\ell$ one has
\begin{equation}\label{plusone}
N(q,E)=1+N_{\rm ns}(q,E)
\end{equation}
where the cardinality $N_{\rm ns}(q,E)$ of the set of the non-singular points is given, for any power
$q=p^\ell$, as follows (\cf~\cite{Silverman} exercise 3.5 page 104)\vspace{.05in}

$\bullet$~$N_{\rm ns}(q,E)=q-1$, if $E$ has split multiplicative reduction over $\F_q$.\vspace{.05in}

$\bullet$~$N_{\rm ns}(q,E)=q+1$, if $E$ has non-split multiplicative reduction over $\F_q$.\vspace{.05in}

$\bullet$~$N_{\rm ns}(q,E)=q$, if $E$ has additive reduction over $\F_q$.\vspace{.05in}

Thus, it follows from the definition \eqref{sumfrob}, that the value $t(p^\ell)$ assumes one of the  following three descriptions\vspace{.05in}

$\bullet$~$t(p^\ell)=1$ for all $\ell$, if $E$ has split multiplicative reduction at $p$.\vspace{.05in}

$\bullet$~$t(p^\ell)= (-1)^\ell$, if $E$ has non-split multiplicative reduction at $p$.\vspace{.05in}

$\bullet$~$t(p^\ell)=0$, if $E$ has additive reduction at $p$.\vspace{.05in}

The second case arises since the split property of $E$ depends upon the parity of $\ell$.
The equality \eqref{checkequ} is checked for all primes of bad reduction $p\in S$. Together with \eqref{localf} and \eqref{localg}, this implies \eqref{taste1}.
\endproof

\begin{rem}\label{positivefun}{\rm
Note that the function $N(n)=n+1-t(n)$ is always positive. Indeed, it is enough to show that for all $n$ one has $|t(n)|\leq n$. Since $t(n)$ is multiplicative, it is enough to show the above inequality when $n=p^k$ is a prime power. In the  case of good reduction, one has $t(p^k)=\alpha_p^k+\bar{\alpha_p}^k$ and since the modulus of $\alpha_p$ is $\sqrt{p}$ one gets
$$
|t(q)|\leq 2 \sqrt{q},\qquad  q=p^k\,.
$$
This proves that $|t(q)|\leq q$ for all $q$ satisfying $2 \sqrt{q}\leq q$ \ie when $q\geq 4$. For $q=2$  one has $2\sqrt 2\cong 2.828<3$ and since $|t(2)|$ is an integer one obtains $|t(2)|\leq 2$. Similarly, for $q=3$ one has $2\sqrt 3\cong 3.4641<4$ and thus $|t(3)|\leq 3$. In the case of bad reduction one has always $|t(q)|\leq 1$.
}\end{rem}\vspace{.05in}

Definition~\ref{modifiedzeta} gives the following equation for the modified zeta function $\zeta_N^{\rm disc}(s)$  associated to the counting function $N(q,E)$
$$
\frac{\partial_s\zeta_N^{\rm disc}(s)}{\zeta_N^{\rm disc}(s)}=-\sum_{n=1}^\infty  (n+1-t(n))\, n^{-s-1}=-\zeta(s+1)-\zeta(s)+R(s+1,E).
$$
Using Lemma~\ref{taste} we then obtain
\begin{thm}\label{tastethm} Let $E$ be an elliptic curve over $\Q$ and let $L(s,E)$ be the associated $L$-function. Then
the modified zeta function $\zeta_N^{\rm disc}(s)$ of $E$  fulfills the following equality
\begin{equation}\label{mod3}
\frac{\partial_s\zeta_N^{\rm disc}(s)}{\zeta_N^{\rm disc}(s)}=-\zeta(s+1)-\zeta(s)+\frac{L(s+1,E)}{\zeta(2s+1)M(s+1)}
\end{equation}
where $\zeta(s)$ is the Riemann zeta function.
\end{thm}

 The Riemann zeta function  has trivial zeros at the points $-2n$ for $n\geq 1$. They generate singularities of $\zeta_E(s)$ at the points $s=-n-\frac 12$, \ie at $s=-\frac 32,-\frac 52,\ldots$.  If RH holds, the non-trivial zeros of the Riemann zeta function generate singularities at points which have real part $-\frac 14$. Note that the poles of the archimedean local factor $N_E^{s/2}(2\pi)^{-s}\Gamma(s)$ ($N_E=$ conductor) of the $L$-function of $E$ determine trivial zeros. These zeros occur for $s\in -\N$ and do not cancel the above singularities. Finally, we have all the zeros of $M(s+1)$. Each factor $(1-p^{-(2s-1)})$ contributes by an arithmetic progression with real part $-\frac 12$
$$
s\in -\frac 12 +\bigcup_{p\in S} \frac{\pi i}{\log p}\Z.
$$
Note that the pole at $-\frac 12$ has order the cardinality of $S$.

For a better understanding of the role played by the primes of bad reduction we consider the following example.

\begin{figure}
\begin{center}
\includegraphics[scale=0.8]{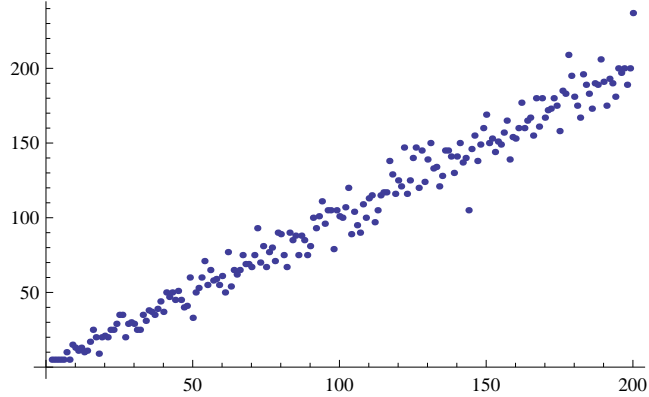}
\end{center}
\caption{Graph of the function $N(n)=n+1-t(n)$ for $n\leq 200$, for the elliptic curve of the Example~\ref{concrete}. \label{counting} }
\end{figure}
\begin{figure}
\begin{center}
\includegraphics[scale=0.5]{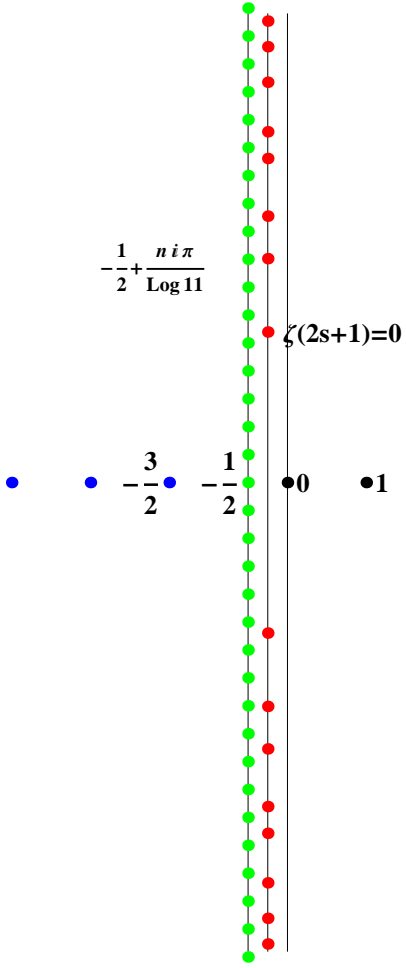}
\end{center}
\caption{Singularities of $\zeta_N^{\rm disc}(s)$. \label{sing} }
\end{figure}

\begin{example}\label{concrete}
{\rm Let consider the elliptic curve in $\P^2(\Q)$ given by the zeros of the homogeneous equation
$$
E :~ Y^2Z + Y Z^2 = X^3 -X^2Z - 10XZ^2 - 20Z^3.
$$
Then the sequence $a(n)$ is given by the coefficients of the modular form
$$
F_E(q)=\sum_{n=1}^\infty
a(n)q^n =q\prod_{n=1}^\infty(1-q^n)^2(1-q^{11n})^2.
$$
These coefficients are easy to compute, for small values of $n$, using the Euler formula
$$
\prod_{n=1}^\infty(1-q^n)=1 + \sum_{n=1}^\infty (-1)^n(q^{(3n^2-n)/2}+q^{(3n^2+n)/2}).
$$
The first values are $$a(0)=0, a(1)=1, a(2)=-2, a(3)=-1, a(4)= 2, a(5)= 1, a(6)= 2, a(7)=-2, \ldots$$
The function $N(n)$ is given by
$N(n)=n+1-t(n)$. It is important to check that it takes positive values. Its graph is depicted in Figure \ref{counting}. In this example, the only prime of bad reduction is $p=11$ and the singularities
of the modified zeta function $\zeta_N^{\rm disc}(s)$ are described in Figure \ref{sing}.}
\end{example}

The challenging question is that to define a natural functor $\underline E$ to finite sets, such that $\#\underline E(\F_{1^n})=N(n+1)$ for all $n$.\vspace{.05in}

\section{Appendix: Primes of bad reduction of an elliptic curve} \label{sing23}\vspace{.05in}

In this appendix we supply more details on the process of counting of points used in the proof of Lemma \ref{taste}.
 In order to
reduce modulo $p$ an elliptic curve $E$ over $\Q$, the first step is to take an equation of the curve which is
in minimal Weierstrass form over the local field $\Q_p$ \ie of the form
\begin{equation}\label{Wegen}
    y^2+a_1 xy+a_3y=x^3+a_2x^2+a_4 x+ a_6
\end{equation}
with $a_j\in\Z_p$ such that the power of $p$ dividing the discriminant
 is minimal. Such a minimal Weierstrass equation can always be found (\cf~\cite{Silverman} Proposition VII 1.3) and is
unique up to a change of variables of the form
$$
x=u^2 x'+r\,, \ \ y=u^3 y'+u^2s x'+t
$$
for $u\in\Z_p^*$ and $r,s,t\in \Z_p$.
One then considers the reduced equation
\begin{equation}\label{Weiminred}
     y^2+\bar a_1 xy+\bar a_3y=x^3+\bar a_2x^2+\bar a_4 x+\bar  a_6
\end{equation}
where $\bar a_j$ is the residue of $a_j$ modulo $p$. This equation is unique up to the standard changes of coordinates on $\F_p$ and  defines a curve over $\F_p$. One then looks for points of this curve \ie solutions of equation \eqref{Weiminred} (together with the point at $\infty$), in the extensions $\F_q$, $q=p^\ell$.
When $p>2$, one can always write the reduced equation in the form
\begin{equation}\label{Weiclass}
     y^2=x^3+\bar a_2x^2+\bar a_4 x+\bar  a_6.
\end{equation}
Thus assuming $p\neq 2$, a singular point has second coordinate $y=0$ while $x$ is the common root of
$x^3+\bar a_2x^2+\bar a_4 x+\bar  a_6=0$, with $3x^2+2\bar a_2x+\bar a_4 =0$. This proves that the elliptic curve has a unique singular point and that both its coordinates belong to $\F_p$. After a translation of $x$ one can then write the equation of the curve in the form $y^2=x^2(x-\beta)$ with $\beta\in \F_p$. There are three cases\vspace{.05in}

$\bullet$~$\beta=0$ is the case of additive reduction: one has a cusp at $(0,0)$.\vspace{.05in}

$\bullet$~$\beta \neq 0,\,-\beta\in \F_p^2$ means split multiplicative reduction.\vspace{.05in}

$\bullet$~$\beta \neq 0,\,-\beta\notin \F_p^2$ is the remaining case of  non-split multiplicative reduction.\vspace{.05in}

Let now $\K$ be a perfect field of characteristic $2$. One starts with a Weierstrass equation in general form
\begin{equation}\label{Weiminredbis}
      y^2+a_1 xy+a_3y=x^3+a_2x^2+a_4 x+ a_6.
\end{equation}
 Then, when the curve is singular it can be written, using affine transformations of $(x,y)$ with coefficients in $\K$, in the form
\begin{equation}\label{char2}
    y^2=x^3+a_4 x+ a_6.
\end{equation}
We refer to \cite{Silverman} (Appendix A, Propositions 1.2 and 1.1) for details. Indeed, the curve is  singular if and only if its discriminant vanishes (\cf~Proposition 1.2 in \opcit). Then, one can  use the description of the curve as in Proposition 1.1 c) of \opcit to deduce that $a_3=0$, since $\Delta=0$. The point at $\infty$ is always non-singular. A point $(x,y)$ is singular if and only if  \eqref{char2} holds together with $3x^2+a_4=0$. This equation is equivalent to $x^2=a_4$ and uniquely determines $x$ since the field is perfect. Then one gets $y^2=a_6$ since $x^3+a_4 x=x(x^2+a_4)=0$. This in turns uniquely determines $y$ since $\K$ is perfect. Thus, for any finite field of characteristic two, one concludes that, if the curve is singular, it has exactly one singular point. In fact, for elliptic curves over $\Q$, there are only $4$ cases to look at in characteristic two. They correspond to the values $a_4\in\{0,1\}$ and $a_6\in\{0,1\}$.

One can always assume that $a_6=0$ by replacing $y$ with $y+1$, if needed.
For $a_4=a_6=0$ one has a cusp at $(0,0)$ and the points over $\F_q$ are labeled  by the $x$ coordinate. Thus there are $q+1$ points over $\F_q$ when $q=2^\ell$.
For $a_4=1$, $a_6=0$, the singular point is $P=(1,0)$. After shifting $x$ by $1$ the equation becomes $y^2=x^2(x+1)$ and thus the singular point is also a cusp. The points over $\F_q$ are labeled as
$(1+t^2,t(1+t^2))$ for $t\in \F_q$ (together with the point at $\infty$). Again, there are $q+1$ points over $\F_q$ when $q=2^\ell$.

\end{document}